\documentclass[11pt]{article}
\usepackage[utf8]{inputenc}
\usepackage{amsmath,amssymb,epsfig,bbm}
\usepackage{stmaryrd,mathabx}
\usepackage{comment}
\usepackage{color}
\usepackage[T1]{fontenc}

\usepackage[textsize=small]{todonotes}
\usepackage{enumitem}
\usepackage{varwidth}
\setlist{nolistsep}
\usepackage{hyperref}

\newtheorem{defi}{Definition}[section]
\newtheorem{prop}[defi]{Proposition}
\newtheorem{theo}[defi]{Theorem}
\newtheorem{conj}[defi]{Conjecture}
\newtheorem{lemm}[defi]{Lemma}
\newtheorem{coro}[defi]{Corollary}
\newtheorem{rema}[defi]{Remark}
\newtheorem{exem}[defi]{Example}
\newtheorem{exems}[defi]{Examples}

\newcommand{\bdefi}{\begin{defi}}
\newcommand{\edefi}{\end{defi}}
\newcommand{\bprop}{\begin{prop}}
\newcommand{\eprop}{\end{prop}}
\newcommand{\btheo}{\begin{theo}}
\newcommand{\etheo}{\end{theo}}
\newcommand{\blemm}{\begin{lemm}}
\newcommand{\brema}{\begin{rema}}
\newcommand{\erema}{\end{rema}}
\newcommand{\bexer}{\begin{exem}}
\newcommand{\eexer}{\end{exem}}
\newcommand{\bexems}{\begin{exems}}
\newcommand{\eexems}{\end{exems}}
\newcommand{\bconj}{\begin{conj}}
\newcommand{\econj}{\end{conj}}
\newcommand{\elemm}{\end{lemm}}
\newcommand{\bcoro}{\begin{coro}}
\newcommand{\ecoro}{\end{coro}}
\newcommand{\dem}{\noindent{\bf Proof. }}
\newcommand{\rem}{\noindent{\bf Remark. }}

\usepackage[english]{babel}
\usepackage{bm, mathtools}

\newcommand{\cA}{\mathcal{A}}

\newcommand{\cJ}{\mathcal{J}}
\newcommand{\cL}{\mathcal{L}}
\newcommand{\cM}{\mathcal{M}}
\newcommand{\cN}{\mathcal{N}}
\newcommand{\cO}{\mathcal{O}}
\newcommand{\cP}{\mathcal{P}}
\newcommand{\cQ}{\mathcal{Q}}
\newcommand{\cX}{\mathcal{X}}
\newcommand{\cY}{\mathcal{Y}}

\newcommand{\bF}{\mathbb{F}}
\newcommand{\bP}{\mathbb{P}}
\newcommand{\bQ}{\mathbb{Q}}
\newcommand{\bR}{\mathbb{R}}
\newcommand{\bT}{\mathbb{T}}
\newcommand{\bZ}{\mathbb{Z}}

\newcommand{\Bad}{\mb{Bad}}
\newcommand{\dg}{\deg}
\newcommand{\dimH}{\dim_{\rm Haus}}
\newcommand*{\rom}[1]{\expandafter\@slowromancap\romannumeral #1@}
\newcommand\idist[1]{\langle\,#1\,\rangle}
\newcommand\mb[1]{\mathbf{#1}}
\newcommand{\Supp}{\operatorname{Supp}}
\newcommand{\ta}{{\mathfrak{a}}}
\newcommand{\ASL}{\operatorname{ASL}}
\newcommand{\inj}{\operatorname{inj}}
\newcommand{\bbr}{{\mb{r}}}

\usepackage{mathrsfs}
\renewcommand\mathcal{\mathscr}

\newcommand{\A}{{\cal A}}
\newcommand{\B}{{\cal B}}
\newcommand{\C}{{\cal C}}

\newcommand{\E}{{\cal E}}

\newcommand{\G}{{\cal G}}

\newcommand{\J}{{\cal J}}

\newcommand{\M}{{\cal M}}

\newcommand{\OOO}{{\cal O}}
\renewcommand{\P}{{\cal P}}
\newcommand{\Q}{{\cal Q}}

\newcommand{\maths}[1]{{\mathbb #1}}  

\newcommand{\FF}{\maths{F}}

\newcommand{\NN}{\maths{N}}

\newcommand{\PP}{\maths{P}}

\newcommand{\RR}{\maths{R}}

\newcommand{\TT}{\maths{T}}

\newcommand{\ra}{\rightarrow}

\newcommand{\ov}[1]{{\overline #1}} 
\newcommand{\wt}[1]{{\widetilde{#1}}}
\newcommand{\wh}[1]{{\widehat{#1}}}

\newcommand{\ga}{\gamma}
\newcommand{\Ga}{\Gamma}

\newcommand{\cqfd}{\hfill$\Box$}

\newcommand{\card}{{\operatorname{Card}}}

\newcommand{\covol}{\operatorname{Covol}}

\newcommand{\diam}{{\operatorname{diam}}}
\newcommand{\diag}{{\operatorname{diag}}}

\newcommand{\id}{\operatorname{id}}

\newcommand{\lcm}{\operatorname{lcm}}

\renewcommand{\ln}{\operatorname{ln}}

\newcommand{\vol}{\operatorname{vol}}

\newcommand{\SL}{\operatorname{SL}}
\newcommand{\GL}{\operatorname{GL}}

\newcounter{fig}

\title{On Hausdorff dimension in inhomogeneous Diophantine
  approximation over global function fields}
\author{Taehyeong Kim \and Seonhee Lim \and Fr\'ed\'eric Paulin}

\begin{document}
\bibliographystyle{../alphanum}
\maketitle
\begin{abstract}
In this paper, we study inhomogeneous Diophantine approximation over
the completion $K_v$ of a global function field $K$ (over a finite
field) for a discrete valuation $v$, with affine algebra $R_v$.  We
obtain an effective upper bound for the Hausdorff dimension of the set
\[
\Bad_A(\epsilon)=\left\{\bm{\theta}\in K_v^{\,m} : \liminf_{(\mb{p},\mb{q})\in
  R_v^{\,m} \times R_v^{\,n}, \|\mb{q}\|\to \infty} \|\mb{q}\|^n
\|A\mb{q}-\bm{\theta}-\mb{p}\|^m \geq \epsilon \right\},
\] 
of $\epsilon$-badly approximable targets $\bm{\theta}\in K_v^{\,m}$
for a fixed matrix $A\in\cM_{m,n}(K_v)$, using an effective version of
entropy rigidity in homogeneous dynamics for an appropriate diagonal
action on the space of $R_v$-grids.  We further characterize matrices
$A$ for which $\Bad_A(\epsilon)$ has full Hausdorff dimension for some
$\epsilon>0$ by a Diophantine condition of singularity on average.
Our methods also work for the approximation using weighted ultrametric
distances.
\end{abstract}

\section{Introduction}
\label{sect:intro}

In the theory of inhomogeneous Diophantine approximation of real numbers
by rational ones (in several variables), one studies the distribution of
the vectors $A\,\mb{x} \in \bR^m$ modulo $\bZ^m$, as $\mb{x}$ varies
over $\bZ^n$, near a vector $b \in \bR^m$ for a $m \times n$ real
matrix $A\in \cM_{m,n}(\bR)$. For instance, if $m,n\geq 1$ and
$\idist{\bm{\xi}} =\displaystyle\inf_{\mb{x}\in \mathbb Z^m}
\|\,\xi-\mb{x}\,\|$ denotes the distance from $\xi \in \bR^m$ to a
nearest integral vector with respect to the Euclidean norm $\|\cdot
\|$ on $\bR^m$, using the inhomogeneous Khintchine-Groshev theorem of
\cite[Theorem1]{Schmidt64}, we have
\[
\liminf_{\mb{x}\in\bZ^n,\;\|\,\mb{x}\,\| \to \infty}
  \|\,\mb{x}\,\|^{n} \idist{A\,\mb{x}-b}^{m}=0
\]
for almost every $(A,b) \in \cM_{m,n}(\bR) \times \bR^m$. 

Let us consider the exceptional set of solutions $(A,b)$ of the above
equation. We call $A$ {\it badly approximable for $b$}
if
\[
\liminf_{\mb{x}\in\bZ^n,\;\|\,\mb{x}\,\| \to \infty}
  \|\,\mb{x}\,\|^{n} \idist{A\,\mb{x}-b}^{m} > 0\;.
\]
If the left hand side is at least $\epsilon$, we say that $A$ is {\it
  $\epsilon$-bad} for $b$.  It is known that given any $b\in \RR^m$,
the set of badly approximable matrices $A\in \cM_{m,n}(\bR)$ has zero
Lebesgue measure but full Hausdorff dimension, see \cite{Schmidt66,
  EinTse11}.

In this paper, replacing $\RR$ by a completion of a function field and
using weighted norms, given an $m\times n$ matrix $A$, under a
necessary assumption of non-singularity on average, we prove that the
set of $m$-dimensional vectors with respect to which $A$ is
$\epsilon$-bad does not have full Hausdorff dimension, and we obtain
an explicit upper bound: there exists $c(A)>0$, depending only on $m$,
$n$ and $A$, such that for every $\epsilon>0$, the Hausdorff dimension
of the set of $m$-dimensional vectors that are $\epsilon$-bad for $A$
is bounded from above by $m-c(A)\,\frac{\epsilon}{\ln(1/\epsilon)}$.
See \cite{KimKimLim21} for analogous results in the case of the real
field. Let us state our main results, referring to Section
\ref{subsec:funcfield} for more precise definitions.

\medskip
Let $K$ be any global function field over a finite field $\bF_q$ of
$q$ elements for a prime power $q$, that is, the function field of a
geometrically connected smooth projective curve $\mb{C}$ over
$\bF_q$. The most studied example in Diophantine approximation in
positive characteristic is the case of the field $K=\FF_q(Z)$ of
rational fractions in one variable $Z$ over $\bF_q$, where
$\mb{C}=\PP^1$ is the projective line, but we emphasize the fact that
our work applies in the general situation above.

We fix a (normalized) discrete valuation $v$ on $K$. Let $K_v$ and
$\cO_v$ be the completion of $K$ with respect to $v$ and its valuation
ring, respectively. We fix a uniformizer $\pi_v \in K$, which
satisfies $v(\pi_v)=1$. Let $k_v = \cO_v/ \pi_v\cO_v$ be the residual
field and let $q_v$ be its cardinality. The (normalized) absolute
value $| \cdot |$ associated with $v$ is defined by $|\,x\,| =
q_v^{-v(x)}$. For every $\sigma\in\bZ_{\geq 1}$, let $\|\;\|:K_v^{\,\sigma}\ra
[0,+\infty[$ be the norm $(\xi_1,\dots,\xi_\sigma)\mapsto\max_{1\leq
i\leq \sigma}|\,\xi_i\,|$. We denote by $\dimH$ the Hausdorff
dimension of the subsets of $K_v^{\,\sigma}$ for this standard
norm. 

The discrete object analogous to the set of integers $\bZ$ in $\bR$ is
the affine algebra $R_v$ of the curve $\mb{C}$ minus the point $v$. If
$K=\bF_{q}(Z)$ and $v=[1:0]$ is the standard point at infinity of
$\mb{C}=\PP^1$, then $R_v = \bF_q[Z]$ is the ring of polynomials in
$Z$ over $\bF_q$.

Let $m,n\in\bZ_{\geq 1}$. Let us fix, throughout the paper, two {\it
  weights} consisting of a $m$-tuple $\mb{r}=(r_1,\cdots,r_m)$ and a
$n$-tuple $\mb{s}=(s_1,\cdots,s_n)$ of positive integers such that we
have $|\mb{r}|=\displaystyle\sum_{1\leq i\leq m}r_i=\displaystyle\sum_{1\leq
  j\leq n}s_j$. The {\it $\mb{r}$-quasinorm} of $\bm{\xi}\in K_v^{\,m}$
and {\it $\bm{s}$-quasinorm} of $\bm{\theta}\in K_v^{\,n}$ are given by
\[
\|\,\bm{\xi}\,\|_{\mb{r}}=
\max_{1\leq i\leq m}|\,\bm{\xi}_i\,|^{\frac{1}{r_i}}
\quad\textrm{and}\quad
\|\,\bm{\theta}\,\|_{\mb{s}}=
\max_{1\leq j\leq n}|\,\bm{\theta}_j\,|^{\frac{1}{s_j}}\;.
\]
We denote by $\idist{\bm{\xi}}_\mb{r} =\displaystyle\inf_{\mb{x}\in
  R_v^{\,m}} \|\,\bm{\xi}-\mb{x}\,\|_{\mathbf{r}}$ the (weighted)
distance from $\bm{\xi}$ to the set $R_v^{\,m}$ of integral vectors in
$K_v^{\,m}$.

Let $\epsilon >0$. A matrix $A\in\cM_{m,n}(K_v)$ is said to be {\it
  $\epsilon$-bad} for a vector $\bm{\theta}\in K_v^{\,m}$ if
\begin{equation}\label{eq:1523}
\liminf_{\mb{x}\in R_v^{\,n},\; \|\,\mb{x}\,\|_{\mb{s}} \to \infty}
\|\,\mb{x}\,\|_{\mb{s}}\;\idist{A\mb{x}-\bm{\theta}}_{\mb{r}} \geq \epsilon \;.
\end{equation}
Denote by $\mb{Bad}_A(\epsilon)$ the set of vectors
$\bm{\theta}\in K_v^{\,m}$ such that $A$ is $\epsilon$-bad for
$\bm{\theta}$.
Given two subsets $U$ and $V$ of a given set, we denote $U-V=\{x\in
U:x\notin V\}$. We say that a matrix $A\in \cM_{m,n}(K_v)$ is
\textit{$(\mb{r},\mb{s})$-singular on average} if for every $\epsilon
>0$, we have
\begin{equation}\label{eq:defisingonaver}
\lim_{N\to\infty}\frac{1}{N}\;\card\{\,\ell\in\{1,\dots,N\}: \exists\;
\mb{y}\in R_v^{\,n} - \{0\},\; \idist{A\,\mb{y}}_{\mb{r}}\leq \epsilon\;
q_v^{-\,\ell},\;
\|\,\mb{y}\,\|_{\mb{s}}\leq q_v^{\ell}\,\}=1\;.
\end{equation}

For the basic example of function field, when $K=\bF_q[Z]$ and
$v=[1:0]$, Bugeaud and Zhang \cite{BugZha19} found a sufficient
condition (and an equivalent one when $n=m=1$) on $A$ for the
Hausdorff dimension of $\mb{Bad}_A(\epsilon)$ to be full.  We first
strenghten and extend their result to general function fields.

\btheo\label{theo:introA1} Let $A\in \cM_{m,n}(K_v)$ be a matrix. The
following assertions are equivalent:
\begin{enumerate}
\item\label{S1_introA1} there exists $\epsilon > 0$ such that the set
  $\mb{Bad}_{A}(\epsilon)$ has full Hausdorff dimension,
\item\label{S2_introA1} the matrix $A$ is $(\mb{r},\mb{s})$-singular
  on average.
\end{enumerate}
\etheo

We also provide an effective upper bound on the Hausdorff dimension in
terms of $\epsilon$, which is a new result even in the basic case
$K=\bF_q[Z]$ and $v=[1:0]$.

\btheo\label{theo:introEff1} For every $A\in \cM_{m,n}(K_v)$ which is not
$(\mb{r},\mb{s})$-singular on average, there exists a constant
$c(A)>0$ depending only on $A$, $\mb{r}$, $\mb{s}$, such that for every
$\epsilon>0$, we have $\dimH \mb{Bad}_{A}(\epsilon)\leq
m-c(A)\frac{\epsilon^{\,|\mb{r}|}}{\ln(1/\epsilon)}$.
\etheo

The proofs of the above main theorems of this paper are largely
divided into two parts. Firstly, assuming the singular on average
property in order to prove the full Hausdorff dimension property, we
give a lower bound on the Hausdorff dimension of appropriately chosen
subsets of $K_v^{\,m}$, using new function fields versions of
classical tools in Diophantine approximation such as geometry of
numbers, transference principle and best approximation vectors (see
for instance \cite{Cassels57, Schmidt80, Kristensen06, KleWei08,
  Cheung11, Chevallier13, GerEvd15, CheChe16, KleShiTom17, GanGho17,
  BugZha19, German20, LiaShiSolTam20, ChoGhoGuaMarSim20,
  BugKimLimRam21}). Secondly, in order to prove the upper bound in
Theorem \ref{theo:introEff1}, we use technics of homogeneous dynamics
of diagonal actions and in particular the entropy method (see for
instance \cite{Kleinbock99,EinLin10,LimSaxSha18,EinLinWar22}). Let us
explain briefly the latter part.

Let $d=m+n$. The dynamical space relevant to inhomogeneous Diophantine
approximation is the space $\cY$ of {\it unimodular grids} $\Lambda+b$
in $K_v^{\,d}$, that is of (Haar-covolume $1$) $R_v$-lattices
$\Lambda$ of $K_v^{\,d}$ translated by vectors $b\in K_v^{\,d}$,
endowed with the affine action of the diagonal subgroup of
$\SL_d(K_v)$. This is in higher dimension more convenient than the
study of the commuting actions of $\SL_d(R_v)$ and of the diagonal
group on the Bruhat-Tits building associated with $\SL_d(K_v)$ (see
\cite[Part III]{BroParPau19} when $d=2$).  Given the above weights
$\mb{r}$ and $\mb{s}$, we consider the affine action on $\cY$ of the
$1$-parameter diagonal subgroup $(\ta^k)_{k\in\bZ}$ where
\[
\ta =\diag(\pi_v^{-r_1 }, \cdots,  \pi_v^{-r_m },
\pi_v^{s_1}, \cdots, \pi_v^{s_n} )\;.
\]
The space $\cY$ of unimodular grids $\ta$-equivariantly projects onto
the space of unimodular $R_v$-lattices $\cX=\SL_d(K_v)/\SL_d(R_v)$ by
the map sending $\Lambda+b$ to $\Lambda$.  We say that a unimodular
$R_v$-lattice $\Lambda$ {\it diverges on average} under the action of
$\ta$ if for every compact subset $Q$ of $\cX$, we have
\[
\liminf_{N\to\infty} \frac{1}{N} \;\card\,\big\{\ell\in \{1, \cdots, N\} :
\ta^\ell\Lambda \notin Q \big\} =1\;.
\]
Following Dani's path, we prove in Section \ref{subsec:Dani} that the
lattice $\Lambda_A=\begin{pmatrix} I_m & A \\ 0& I_n \end{pmatrix}
R_v^{\,d}$ diverges on average under $\ta$ if and only if $A$ is
$(\mb{r},\mb{s})$-singular on average.

As developped in the last Section \ref{sec:upperbound}, the main idea
of the entropy method in our situation is that if the point
$\Lambda_A$ does not diverge on average, then the Hausdorff dimension
of $\mb{Bad}_A(\epsilon)$ provides a lower bound on the conditional
entropy of $\ta$ with respect to a measure $\mu$ constructed by
well-separated sets on the fibers of the projection $\cY \to \cX$. An
effective control of the maximal conditional entropy by a control of
the support of $\mu$ on the thin/thick parts of $\cY$ hence gives an
effective upper bound. See \cite{KimKimLim21} for the real case.

Our paper is organized as follows.  In Section~\ref{sec:prelimI}, we
recall basic facts on the geometry of numbers, define the best
approximation sequence and prove transference principle for the
weighted case, which generalize previous results of Bugeaud-Zhang
\cite{BugZha19}.  In Section~\ref{sect:characsing}, we give a
characterization of the singular on average property with weights in
terms of the best approximation sequence.  In
Section~\ref{sect:fullhausdorff}, we establish the lower bound on the
Hausdorff dimension by constructing a subsequence with controlled
growth of the best approximation sequence for a matrix whose transpose
is singular on average.  In Section~\ref{sec:prelimII}, we recall some
background on homogeneous dynamics and conditional entropy, and prove
an effective and positive characteristic version of the variational
principle for conditional entropy of \cite[\S 7.55]{EinLin10} (see
\cite{KimKimLim21} in the real case). Finally, in Section
\ref{sec:upperbound}, we implement the entropy method.

We remark that in \cite{KimKimLim21}, the first
two authors, with Wooyeon Kim, also show that given any vector 
$b\in\bR^m$, the set of matrices $A\in \cM_{m,n}(\bR)$ 
that are $\epsilon$-bad for $b$ does not have full Hausdorff dimension, 
and estimate an explicit upper bound. Thus it seems very interesting 
to obtain a similar result in the global function field case.

\medskip\noindent{\small {\it Acknowledgements.} Taehyeong Kim and
  Seonhee Lim are supported by the National Research Foundation of
  Korea, Project Number NRF-2020R1A2C1A01011543.  Taehyeong Kim
  is supported by the National Research Foundation of Korea, Project
  Number NRF-2021R1A6A3A13039948. Seonhee Lim is an associate member
  of KIAS.}
 
\section{Background material for the lower bound}
\label{sec:prelimI}

\subsection{On global function fields}
\label{subsec:funcfield}
  
We refer for instance to \cite{Goss96, Rosen02}, as well as \cite[\S
  14.2]{BroParPau19}, for the content of this section. Let $\bF_q$ be
a finite field with $q$ elements, where $q$ is a positive power of a
positive prime. Let $K$ be the function field of a geometrically
connected smooth projective curve $\mb{C}$ over $\bF_q$, or
equivalently an extension of $\bF_q$ with transcendence degree $1$, in
which $\bF_q$ is algebraically closed. We denote by $g$ the genus of
$\mb{C}$. There is a bijection between the set of closed points of
$\mb{C}$ and the set of normalized discrete valuations $v$ of $K$, the
valuation of a given element $f\in K$ being the order of the zero or
the opposite of the order of the pole of $f$ at the given closed
point. We fix such an element $v$ throughout this paper, and use the
notation $K_v$, $\cO_v$, $\pi_v$, $k_v$, $q_v$, $|\cdot|$ defined in
the introduction. We furthermore denote by $\deg v$ the degree of $v$,
so that
$$
q_v=q^{\deg v}\;.
$$

We denote by $\vol_v$ the normalized Haar measure on the locally
compact additive group $K_v$ such that $\vol_v(\cO_v)=1$. For any
positive integer $d$, let $\vol_v^d$ be the normalized Haar measure on
$K_v^{\,d}$ such that $\vol_v^d(\cO_v^d)=1$. Note that for every $g\in
\GL_d(K_v)$ we have
$$
d\vol_v^d (gx)=|\,\det(g)\,| \; d\vol_v^d(x)\;,
$$
where $\det$ is the determinant of a matrix. For every discrete
additive subgroup $\Lambda$ of $K_v^{\,d}$, we again denote by $\vol_v^d$
(and simply $\vol_v$ when $d=1$) the measured induced on
$K_v^{\,d}/\Lambda$ by $\vol_v^d$.

Note that the completion $K_v$ of $K$ for $v$ is the field
$k_{v}((\pi_v))$ of Laurent series $x=\sum_{i\in\bZ}x_i
(\pi_v)^i$ in the variable $\pi_v$ over $k_v$, where $x_i \in k_v$ is
zero for $i\in \bZ$ small enough.  We have
$$
|\,x\,| = q_v ^{-\sup\{j\in\bZ \;:\;  \forall i<j, \;x_i =0\}}\,,
$$
and $\cO_v= k_v[[\pi_v]]$ is the local ring of power series
$x=\sum_{i\in\bZ_{\geq 0}}x_i (\pi_v)^i$ in the variable $\pi_v$ over $k_v$.

Recall that the affine algebra $R_v$ of the affine curve
$\mb{C}-\{v\}$ consists of the elements of $K$ whose only poles are at
the closed point $v$ of $\mb{C}$. Its field of fractions is equal to
$K$, hence we can write elements of $K$ as $x/y$ with $x,y\in R_v$ and
$y\neq 0$. By for instance \cite[Eq.~(14.2)]{BroParPau19}, we have
\begin{equation}\label{eq:RcapO}
  R_v\cap\cO_v=\bF_q\;.
\end{equation}
For every $\xi\in K_v$, we denote by
$$
|\idist{\xi}| = \inf_{x\in R_v} \|\,\xi-x\,\|
$$
the distance in $K_v$ from $\xi$ to the set $R_v$ of integral
points of $K_v$.

For instance, if $\mb{C}$ is the projective line $\bP^1$, if
$\infty=[1:0]$ is its usual point at infinity and if $Z$ is a
variable name, then $g=0$, $K=\bF_q (Z)$, $\pi_{\infty}=Z^{-1}$,
$K_{\infty}=\bF_q ((Z^{-1}))$, $\cO_{\infty}=\bF_q [[Z^{-1}]]$,
$k_{\infty}=\bF_q$, $q_\infty=q$ and $R_{\infty}=\bF_q[Z]$. In this
setting, there are numerous results on Diophantine approximation in
the fields of formal power series, see for instance
\cite{Lasjaunias00}, \cite[Chap.~9]{Bugeaud04}. On the other hand,
little is known about Diophantine approximation over general global
function fields, see for instance \cite{KleShiTom17} (for a single
valuation in positive characteristic) for the ground work on the
geometry of number for function fields.

\subsection{On the geometry of numbers and Dirichlet's theorem}
\label{subsec:geomnumb}

Let $d$ be a positive integer. An {\it $R_v$-lattice} $\Lambda$ in
$K_v^{\,d}$ is a discrete $R_v$-submodule in $K_v^{\,d}$ that generates
$K_v^{\,d}$ as a $K_v$-vector space. The {\it covolume} of $\Lambda$,
denoted by $\covol(\Lambda)$, is defined as the measure of the
(compact) quotient space $K_v^{\,d}/\Lambda$~:
\[
\covol(\Lambda)=\vol_v^d(K_v^{\,d}/\Lambda)\;.
\]
For example, $R_v^{\,d}$ is an $R_v$-lattice in $K_v^{\,d}$, and by for
instance \cite[Lem.~14.4)]{BroParPau19}, we have
\begin{equation}\label{eq:covolRv}
\covol(R_v^{\,d})=q^{(g-1)d}\;.
\end{equation}

Let $\ov{B}(0,r)$ be the closed ball of radius $r$ centered at zero in
$K_v^{\,d}$ with respect to the norm $\|\cdot\|:(\xi_1,\dots,\xi_d)\mapsto
\max_{1\leq i\leq d}|\,\xi_i\,|\,$. For every integer $k\in
\{1,\dots,d\}$, the {\it $k$-th minimum} of an $R_v$-lattice $\Lambda$
is defined by
\[
\lambda_k (\Lambda)=\min\{r>0 : \dim_{K_v}
(\operatorname{span}_{K_v}(\ov{B}(0,r)\cap \Lambda))\geq k\},
\]
where $\operatorname{span}_{K_v}$ denotes the $K_v$-linear span of a
subset of a $K_v$-vector space and $\dim_{K_v}$ is the dimension of a
$K_v$-vector space. Note that $\lambda_1(\Lambda),\dots,
\lambda_d(\Lambda) \in q_v^\bZ$. The next result follows from
\cite[Theo.~4.4]{KleShiTom17} and Equation \eqref{eq:covolRv}.

\btheo\label{theo:Mink2}
{\bf (Minkowski's theorem)}
For every $R_v$-lattice $\Lambda$ in $K_v^{\,d}$, we have
$$
q^{-(g-1)d}\;\covol(\Lambda)\leq
\lambda_1(\Lambda)\dots\lambda_d(\Lambda)\leq q_v^d\; \covol(\Lambda)
\;.\;\;\;\Box
$$
\etheo

Since $\lambda_1(\Lambda)\leq\dots \leq\lambda_d(\Lambda)$, the
following result follows immediately from Minkowski's theorem
\ref{theo:Mink2}.

\bcoro \label{coro:minko}\label{coro:Mink1}
For every $R_v$-lattice $\Lambda$ in $K_v^{\,d}$, we have
$$
\lambda_1(\Lambda)\leq q_v\,\covol(\Lambda)^{\frac{1}{d}}
\;.\;\;\;\Box
$$
\ecoro

\medskip
The following result generalizes \cite[Theo.~2.1]{GanGho17}, which is
proved only when $K=\bF_q(Z)$ and $v=\infty$, to all function fields
$K$ and valuations $v$. See also \cite[Theo.~1.3]{KleWei08} in the
case of the field $\bQ$.

\btheo \label{theo:Dirichlet}{\bf (Dirichlet's theorem)} For every
matrix $A\in\cM_{m,n}(K_v)$ whose rows are denoted by $A_1, \dots$,
$A_m$, for all $(r'_1,\dots,r'_m) \in\bZ_{\geq 0}^{\;m}$ and $(s'_1,\dots,s'_n)
\in\bZ_{\geq 0}^{\;n}$ with
$$
r'_i> 1+\frac{g-1}{\deg v} {\rm ~~~and~~~}
\sum_{i=1}^m r'_i=\sum_{j=1}^n s'_j\;,
$$
there exists an element $\mb{y}=(y_1,\dots, y_n)\in R_v^{\,n}-\{0\}$ such
that, for all $i=1,\dots, m$ and $j=1,\dots, n$, we have
$$
|\langle\, A_i\,\mb{y}\,\rangle |\leq q_v\,q^{g-1}\;{q_v}^{-r'_i}{\rm
  ~~~and~~~} |\,y_j\,|\leq q_v\,q^{g-1}\;{q_v}^{s'_j}\;.
$$
\etheo

\dem
With $A$, $r'_1,\dots,r'_m$ and $s'_1,\dots,s'_n$ as in the statement, we
apply Corollary \ref{coro:minko} with $d=m+n$ to the $R_v$-lattice
$$
\Lambda=\begin{pmatrix}\pi_v^{-r'_1}& & & & & 0\\
& \ddots & & & & \\
& & \pi_v^{-r'_m} & & &\\
& & & \pi_v^{s'_1} & & \\
& & & & \ddots &\\
0 & & & & & \pi_v^{s'_n}\end{pmatrix}
\begin{pmatrix}I_m & A \\ 0 & I_n
\end{pmatrix} \;R_v^{\,d}\;,
$$
where $I_k$ is the $k\times k$ identity matrix. Since the above two
matrices have determinant $1$ by the assumption $\sum_{i=1}^m
r'_i=\sum_{j=1}^n s'_j$, and by Equation \eqref{eq:covolRv}, we have
$\covol(\Lambda)=q^{(g-1)d}$. Corollary \ref{coro:minko}
hence says that there exists $(\mb{x}=(x_1,\dots, x_m),\mb{y}=(y_1,\dots,
y_n))\in R_v^{\,d}-\{0\}$ such that
$$
\max\big\{\max_{i=1,\dots, m}\;|\,\pi_v^{-r'_i}(x_i+A_i\mb{y})\,|,\;
\max_{j=1,\dots, n}|\,\pi_v^{s'_j}y_j\,|\;\big\} \leq
q_v\covol(\Lambda)^{\frac{1}{d}}=q_v\,q^{g-1}\;.
$$
Assume for a contradiction that $\mb{y}=0$. Then for all $i=1,\dots,
m$, since $|\,\pi_v\,|=q_v^{-1}$, we have the inequality
$|\,x_i\,|\leq q_v\,q^{g-1}\,q_v^{-r'_i}$.  Since $r'_i>
1+\frac{g-1}{\deg v}$, this would imply that $|\,x_i\,|<1$. By
Equation \eqref{eq:RcapO}, we have $\{z\in R_v:|\,z\,| < 1\}=\{0\}$.
Since $x_i\in R_v$, we would have that $\mb{x}=0$, contradicting the
fact that $(\mb{x},\mb{y})\neq 0$. Therefore $\mb{y}\neq 0$ and the
result follows.  \cqfd

\medskip
The following corollary is due to \cite[Theo.~1.1]{Kristensen06} (see
also \cite[Theo.~3.2]{BugZha19} where the assumption that $c\,m$ is
divisible by $n$ is implicit) in the special case when $K=\bF_q(Z)$
and $v=\infty$ and without weights.

Let $\min\mb{r}=\min_{1\leq i\leq m} r_i$ and similarly for
$\min\mb{s}$, $\max\mb{r}$ and $\max\mb{s}$. 

\bcoro \label{coro:Dirichlet}
For all $A\in\cM_{m,n}(K_v)$ and $\alpha\in\bZ_{\geq 0}$ with $\alpha >
\frac{1}{\min\mb{r}}+\frac{g-1}{(\min\mb{r})(\deg v)}$, there exists
$\mb{y}\in R_v^{\,n}-\{0\}$ such that
$$
\idist{A\,\mb{y}}_{\mb{r}}\leq q^{\frac{\deg v + g-1}{\min\mb{r}}}
\;q_v^{-\alpha}{\rm  ~~~and~~~}\|\,\mb{y}\,\|_{\mb{s}}\leq
q^{\frac{\deg v + g-1}{\min\mb{s}}}\;q_v^{\alpha}\;.
$$
\ecoro

\dem
We apply Theorem \ref{theo:Dirichlet} with $r'_i=\alpha\, r_i> 1
+\frac{g-1}{\deg v}$ for $i=1,\dots,m$ and $s'_j=\alpha\, s_j$ for $j=
1,\dots,n$, noting that $\sum_{i=1}^m r'_i=\sum_{j=1}^n s'_j$ since
$\sum_{i=1}^m r_i=\sum_{j=1}^n s_j$.
\cqfd

\medskip \rem When $\mb{r}=(n,n,\dots, n)$ and $\mb{s}=(m,m,\dots,m)$,
the above result says that for every integer $\alpha >
\frac{1}{n}+\frac{g-1}{n\deg v}$, there exists $\mb{y}\in
R_v^{\,n}-\{0\}$ such that
$$
\min_{\mb{x}\in R_v^{\,m}} \|A\,\mb{y} -\mb{x}\| \leq q_v\,q^{g-1}
\;q_v^{-\alpha\,n}{\rm  ~~~and~~~}\|\,\mb{y}\,\|\leq
q_v\,q^{g-1}\;q_v^{\alpha\,m}\;,
$$
where $\|\cdot\|$ is the sup norm.

\subsection{Best approximation sequences with weights}
\label{subsec:bestapproxseq}

In this subsection, we construct a version with weights, valid for all
function fields, of the best approximation sequences associated with a
completely irrational matrix by Bugeaud-Zhang \cite{BugZha19}.

\medskip
A matrix $A\in\cM_{m,n}(K_v)$ is said to be {\it completely
  irrational} if $\idist{A\,\mb{y}}_{\mb{r}}\neq 0$ for every
$\mb{y}\in R_v^{\,n}-\{0\}$. Note that this does not depend on the
weight $\mb{r}$, and that the fact that $A$ is completely irrational might
not necessarily imply that $\,^{t}\!A$ is completely irrational.

\brema\label{rem:noncomplirrat} Let $A\in\cM_{m,n}(K_v)$ be such that
$\,^{t}\!A$ is not completely irrational.
\begin{enumerate}
\item[(1)] The matrix $\,^{t}\!A$ is $(\mb{s},\mb{r})$-singular on average.
\item[(2)] For every $\epsilon>0$ small enough, the  set
  $\mb{Bad}_{A}(\epsilon)$ has full Hausdorff dimension.
\end{enumerate}
\erema

\dem By assumption, there exist $\mb{x}\in R_v^{\,n}$ and
$\mb{y}=(y_1,\dots, y_m)\in R_v^{\,m}-\{0\}$ such that
$\,^{t}\!A\;\mb{y} -\mb{x}=0$.

\medskip
(1) For every $\epsilon>0$, if $\ell_0=\lceil \log_{q_v}
\|\,\mb{y}\,\|_{\mb{r}}\rceil$ then for all integers $N\geq \ell_0$
and $\ell\in\{\ell_0,\dots,N\}$, we have
$\idist{\,^{t}\!A\,\mb{y}}_{\mb{s}} =0\leq \epsilon\; q_v^{-\,\ell}$
and $\|\,\mb{y}\,\|_{\mb{r}}\leq q_v^{\ell}$, hence $\,^{t}\!A$ is
$(\mb{s},\mb{r})$-singular on average (see Equation
\eqref{eq:defisingonaver}).

\medskip
(2) For every $\bm{\theta}=(\theta_1,\dots, \theta_m)\in K_v^{\,m}$, let
\[
\mb{y} \cdot \bm{\theta}=\sum_{j=1}^m \;y_i\,\theta_i\;\in K_v\;.
\]
For every $\epsilon\in\;]0,\frac{1}{\|\,\mb{y}\,\|_{\mb{r}}}]$, let
$U_{\mb{y},\epsilon}=\big\{\bm{\theta}\in K_v^{\,m}: |\idist{\mb{y}
  \cdot \bm{\theta}}| \geq (\,\epsilon\; \|\,\mb{y}\,\|_{\mb{r}}
)^{\min\mb{r}}\big\}$.  If $\epsilon$ is small enough, then the set
$U_{\mb{y},\epsilon}$ contains a closed ball of positive radius: For
instance, let $j_0\in\{1,\dots,m\}$ be such that $y_{j_0}\neq 0$~;
define $\theta_{0,j}=0$ if $j\neq j_0$, $\theta_{0,j_0} =\frac{\pi_v}
{y_{j_0}}$ and $\bm{\theta}_{0}= (\theta_{0,1}, \dots,\theta_{0,m})$~;
then it is easy to check using the ultrametric inequality that the
closed ball $\ov{B}(\bm{\theta}_0,\frac{1} {q_v^2\;\|\,\mb{y}\,\|})$ is
contained in $U_{\mb{y},\epsilon}$ if $\epsilon< \,q_v^{-\frac{1}
  {\min\mb{r}}}\, \|\,\mb{y}\,\|_{\mb{r}}^{-1} \,$.

Let us prove that $\mb{Bad}_{A}(\epsilon)$ contains
$U_{\mb{y},\epsilon}$, which implies that $\dimH\big(\mb{Bad}_{A}
(\epsilon)\big)=m$ if $\epsilon$ is small enough.  Let $\bm{\theta}\in
U_{\mb{y},\epsilon}$ and $(\mb{y}',\mb{x}')\in R_v^{\,m}\times
(R_v^{\,n}-\{0\})$.

If $\|\,\mb{y}\,\|_{\mb{r}} \|A\mb{x}'+\mb{y}'-\bm{\theta}\,\|_{\mb{r}}\geq
1$, then since $\mb{x}'\in R_v^{\,n}-\{0\}$ so that
$\|\,\mb{x}'\,\|_{\mb{s}}\geq 1$, we have
$$
\|\,\mb{x}'\,\|_{\mb{s}}\;\|A\mb{x}'+\mb{y'}-\bm{\theta}\,\|_{\mb{r}}
\geq \frac{1}{\|\,\mb{y}\,\|_{\mb{r}}}
\;\|\,\mb{y}\,\|_{\mb{r}} \;\|A\mb{x}'+\mb{y'}-\bm{\theta}\|_{\mb{r}}\geq
\frac{1}{\|\,\mb{y}\,\|_{\mb{r}}}\geq \epsilon\;.
$$

If $\|\,\mb{y}\,\|_{\mb{r}} \;\|A\mb{x'}+\mb{y'}-\bm{\theta}\,\|_{\mb{r}}\leq
1$, then since $\mb{y}\cdot(A\mb{x'}+\mb{y'})=
(\,^{t}\!A\mb{y})\cdot\mb{x'}+\mb{y}\cdot\mb{y'}=
\mb{x}\cdot\mb{x'}+\mb{y}\cdot\mb{y'}\in R_v$, and since
$\bm{\theta}\in U_{\mb{y},\epsilon}$, we have
\begin{align*}
\|\,\mb{x}'\,\|_{\mb{s}}\;\|A\mb{x'}+\mb{y'}-\bm{\theta}\,\|_{\mb{r}}
&\geq \frac{1}{\|\,\mb{y}\,\|_{\mb{r}}}
\;\|\,\mb{y}\,\|_{\mb{r}}\; \|A\mb{x'}+\mb{y'}-\bm{\theta}\,\|_{\mb{r}}\\
&\geq \frac{1}{\|\,\mb{y}\,\|_{\mb{r}}}\;
\big(\max_{1\leq j\leq m}\|\,\mb{y}\,\|_{\mb{r}}^{\,r_j}
\|A\mb{x'}+\mb{y'}-\bm{\theta}\,\|_{\mb{r}}^{\,r_j}\big)^{\frac{1}{\min\mb{r}}}
\\&\geq\frac{1}{\|\,\mb{y}\,\|_{\mb{r}}}\;
\big|\,\mb{y}\cdot(A\mb{x'}+\mb{y'}-\bm{\theta})\big|^{\frac{1}{\min\mb{r}}}
\geq\frac{1}{\|\,\mb{y}\,\|_{\mb{r}}}\;
|\idist{\,\mb{y}\cdot\bm{\theta}}|^{\frac{1}{\min\mb{r}}}
\geq \epsilon\;.
\end{align*}
Therefore $\bm{\theta}\in\mb{Bad}_{A}(\epsilon)$, as wanted.
\cqfd

\bigskip
For every matrix $A\in\cM_{m,n}(K_v)$, a {\it best approximation
  sequence} for $A$ with weights $(\mb{r},\mb{s})$ is a sequence
$(\mb{y}_i)_{i\geq 1}$ in $R_v^{\,n}$ such that, with
$Y_i=\|\,\mb{y}_i\,\|_{\mb{s}}$ and
$M_i=\idist{A\,\mb{y}_i}_{\mb{r}}$,

$\bullet$~ the sequence $(Y_i)_{i\geq 1}$ is positive and
strictly increasing,

$\bullet$~ the sequence $(M_i)_{i\geq 1}$ is positive and
strictly decreasing, and

$\bullet$~ for every $\mb{y}\in R_v^{\,n}-\{0\}$ with
$\|\,\mb{y}\,\|_{\mb{s}}<Y_{i+1}$, we have
$\idist{A\,\mb{y}}_{\mb{r}}\geq M_i$.

\noindent
We denote by $\lcm\mb{r}$ the least common multiple of $r_1,\dots,
r_m$, and similarly for $\lcm\mb{s}$.

\blemm\label{lem:bestapprox} Assume that $A\in\cM_{m,n}(K_v)$ is completely
irrational.

\medskip
\noindent
(1) There exists a best approximation sequence
$(\mb{y}_i)_{i\geq 1}$ for $A$ with weights $(\mb{r},\mb{s})$.

\medskip
\noindent
(2) If $(\mb{y}_i)_{i\geq 1}$ is a best approximation sequence for $A$
with weights $(\mb{r},\mb{s})$, then
\begin{enumerate}
\item[i)] we have $M_i\in q_v^{\frac{1}{\lcm \mb{r}}\bZ}$ and
  $M_i\in q_v^{\frac{1}{\lcm \mb{r}}\bZ_{\leq 0}}$ if $i$ is large enough,
\item[ii)] we have $Y_i\in q_v^{\frac{1}{\lcm \mb{s}}\bZ_{\geq 0}}$
  and $Y_i\geq q_v^{\frac{i-1}{\lcm \mb{s}}}$ for every $i\geq 1$,
\item[iii)] the sequence $\big(M_i\,Y_{i+1}
  \big)_{i\geq 1}$ is uniformly bounded.
\end{enumerate}
\elemm

Note that a best approximation sequence might be not unique (and the
terminology ``best'', though traditional, is not very appropriate).
When $m=n=r_1=s_1=1$, $K=\bF_q(Z)$ and $v=\infty$, then $A\in K_v$ is
completely irrational if and only if $A\in K_v-K$, and with
$\big(\frac{P_k}{Q_k}\big)_{k\geq 0}$ the sequence of convergents of
$A$ (see for instance \cite{Lasjaunias00}), we may take $y_i=Q_{i-1}$
for all $i\geq 1$.

If $A\in\cM_{m,n}(K_v)$ is not completely irrational, a {\it best
  approximation sequence} for $A$ with weights $(\mb{r},\mb{s})$ is a
finite sequence $(\mb{y}_i)_{1\leq i\leq i_0}$ in $R_v^{\, n}$, such that,
with $Y_i=\|\,\mb{y}_i\,\|_{\mb{s}}$ and
$M_i=\idist{A\,\mb{y}_i}_{\mb{r}}$,

$\bullet$~ $1=Y_1<\dots <Y_{i_0}$,

$\bullet$~ $M_1>\dots>M_{i_0}=0$,

$\bullet$~ for all $i\in\{1,\dots,i_0-1\}$ and $\mb{y}\in R_v^{\,n}-\{0\}$
with $\|\,\mb{y}\,\|_{\mb{s}}<Y_{i+1}$, we have
$\idist{A\,\mb{y}}_{\mb{r}}\geq M_i$, and

$\bullet$~ which stops at the first $i_0$ such that there exists
$\mb{z}\in R_v^{\,n}$ with $0<\|\,\mb{z}\,\|_{\mb{s}}\leq Y_{i_0}$ and
$\idist{A\,\mb{z}}_{\mb{r}}=0$.

The proof of Lemma \ref{lem:bestapprox} is similar to the one given
after \cite[Def.~3.3]{BugZha19} in the particular case when
$K=\bF_q(Z)$, $v=\infty$  and without weights.

\medskip
\dem (1) Let us prove by induction on $i\geq 1$ that there exist
$\mb{y}_1,\dots, \mb{y}_i$ in $R_v^{\,n}$ such that, with $Y_j =
\|\,\mb{y}_j\,\|_{\mb{s}}$ and $M_j =\idist{A\,\mb{y}_j}_{\mb{r}}$ for
every $1\leq j\leq i$, we have $1=Y_1<\dots <Y_i$, $M_1>\dots >M_i>0$,
and (using $M_0=+\infty$ by convention)

$(a_i)$~ we have $\idist{A\, \mb{y}}_{\mb{r}}\geq M_{i-1}$ for every
$\mb{y}\in R_v^{\,n}-\{0\}$ with $\|\,\mb{y}\,\|_{\mb{s}}<Y_{i}$,

$(b_i)$~ we have $\idist{A\, \mb{y}}_{\mb{r}}\geq M_{i}$ for every
$\mb{y}\in R_v^{\,n}-\{0\}$ with $\|\,\mb{y}\,\|_{\mb{s}} \leq Y_{i}$.

\medskip
Note that $\{x\in R_v:|\,x\,|\leq 1\}= R_v\cap\cO_v=\bF_q$ by Equation
\eqref{eq:RcapO}. Hence the elements with smallest $\mb{s}$-quasinorm
in $R_v^{\,n}-\{0\}$ are the elements in the finite set
$\bF_q^{\,n}-\{0\}$, which is the set of elements in $R_v^{\,n}$ with
$\mb{s}$-quasinorm $1$. Furthermore, the set
$\{\|\,\mb{y}\,\|_{\mb{s}} : \mb{y}\in R_v^{\,n}-\{0\}\}$ is contained
in $q_v^{\bigcup_{j=1}^n \frac{1}{s_j} \bZ_{\geq 0}}\subset
q_v^{\frac{1}{\lcm \mb{s}} \bZ_{\geq 0}}$.  Similarly, for every
$\mb{x}\in K_v^{\,m}-\{0\}$, we have $\idist{\mb{x}}_{\mb{r}} \in
q_v^{\frac{1}{\lcm \mb{r}}\bZ}$.

Therefore there exists an element $\mb{y}_1\in R_v^{\,n}$ with
$\|\,\mb{y}_1\,\|_{\mb{s}} =1$ such that
$$
\idist{A\,\mb{y}_1}_{\mb{r}} =\min\{\;\idist{A\,\mb{y}}_{\mb{r}} :
\mb{y}\in R_v^{\,n},\;\; \|\,\mb{y}\,\|_{\mb{s}}=1\,\}\;.
$$
We thus have $Y_1=\|\,\mb{y}_1\,\|_{\mb{s}}=1$ and $M_1=
\idist{A\,\mb{y}_1}_{\mb{r}} >0$ since $A$ is completely
irrational. There is no $\mb{y}\in R_v^{\,n}-\{0\}$ with
$\|\,\mb{y}\,\|_{\mb{s}} <Y_{1}$, and if $\|\,\mb{y}\,\|_{\mb{s}}
=Y_{1}$, then $\idist{A\, \mb{y}}_{\mb{r}}\geq M_{1}$, hence the
claims $(a_1)$ and $(b_1)$ are satisfied.

Assume by induction that $\mb{y}_1,\dots, \mb{y}_i$ as above are
constructed.  Let
$$
S=\{\,\mb{y}\in R_v^{\,n}:
\;\|\,\mb{y}\,\|_{\mb{s}}>Y_i,\;\;\idist{A\,\mb{y}}_{\mb{r}}<M_i\,\}\;.
$$
Note that the set $\{\mb{z}\in R_v^{\,n} ,\;\;
0<\|\,\mb{z}\,\|_{\mb{s}} \leq Y_i\}$ is finite by the discreteness of
$R_v^{\,n}$, and $\epsilon_i =\min\{\,\idist{A\,\mb{z}}_{\mb{r}}:
\mb{z}\in R_v^{\,n} ,\;\; 0<\|\,\mb{z}\,\|_{\mb{s}}\leq Y_i\,\}$ is
positive, since $A$ is completely irrational. Corollary
\ref{coro:Dirichlet} of Dirichlet's theorem implies in particular, by
taking in its statement $\alpha$ large enough, that for every
$\epsilon>0$, there exists $\mb{y}\in R_v^{\,n}-\{0\}$ such that
$\idist{A\,\mb{y}}_{\mb{r}} <\epsilon$.  Applying this with $\epsilon
=\min\{M_i,\epsilon_i\}>0$ proves that the set $S$ is nonempty. Hence
the set $S_{\rm min}$ of elements of $S$ with minimal
${\mb{s}}$-quasinorm, which is finite again by the discreteness of
$R_v^{\,n}$, is nonempty.  Therefore there exists $\mb{y}_{i+1}\in
S_{\rm min}$ such that
$$
\idist{A\,\mb{y}_{i+1}}_{\mb{r}}=
\min\{\,\idist{A\,\mb{z}}_{\mb{r}}:\mb{z}\in S_{\rm min}\,\}\;.
$$
Then $Y_{i+1} = \|\,\mb{y}_{i+1}\,\|_{\mb{s}} = \min
\|\,S\,\|_{\mb{s}} > Y_i$ by the definition of the set $S$. We also
have that $M_{i+1} = \idist{A\,\mb{y}_{i+1}}_{\mb{r}} <M_i$ since
$\mb{y}_{i+1}\in S_{\rm min} \subset S$, and again by the definition
of $S$.

\medskip
Let us now prove that $\mb{y}_{i+1}$ satisfies the properties $(a_{i+1})$
and $(b_{i+1})$.

$\bullet$~ Let $\mb{y}\in R_v^{\,n}-\{0\}$ be such that
$\|\,\mb{y}\,\|_{\mb{s}} <Y_{i+1}$.  If $\|\,\mb{y}\,\|_{\mb{s}}\leq
Y_i$, then by the induction hypothesis $(b_i)$, we have
$\idist{A\,\mb{y}}_{\mb{r}}\geq M_{i}$, as wanted for Property
$(a_{i+1})$.  If $\|\,\mb{y}\,\|_{\mb{s}}> Y_i$, then by the
definition of $S$, we have $\idist{A\,\mb{y}}_{\mb{r}}\geq M_{i}$ as
wanted for Property $(a_{i+1})$, otherwise $\mb{y}$ would be an
element of $S$ with ${\mb{s}}$-quasinorm strictly less than the
minimum ${\mb{s}}$-quasinorm of the elements of $S$, a contradiction.

$\bullet$~ Let $\mb{y}\in R_v^{\,n}-\{0\}$ be such that
$\|\,\mb{y}\,\|_{\mb{s}}\leq Y_{i+1}$. Either $\|\,\mb{y}\,\|_{\mb{s}}
< Y_{i+1}$, in which case, as just seen, $\idist{A\,\mb{y}}_{\mb{r}}
\geq M_{i}\geq M_{i+1}$, as wanted for Property $(b_{i+1})$. Or
$\|\,\mb{y}\,\|_{\mb{s}}= Y_{i+1}>Y_i$, in which case either
$\idist{A\,\mb{y}}_{\mb{r}}\geq M_{i}\geq M_{i+1}$, as wanted for
Property $(b_{i+1})$, or $\idist{A\,\mb{y}}_{\mb{r}}< M_{i}$, so that
$\mb{y}$ belongs to $S_{\rm min}$, hence $\idist{A\,\mb{y}}_{\mb{r}}
\geq \min\{\,\idist{A\,\mb{z}}_{\mb{r}} :\mb{z}\in S_{\rm min}\,\}
=M_{i+1}$.

By induction, this proves Assertion (1) of Lemma \ref{lem:bestapprox}.

\medskip
(2) i) This follows from the facts that $M_i\in q_v^{\frac{1}{\lcm
    {\mb r}} \bZ}$ and that $M_{i+1}<M_i$.

\medskip
ii) Since $Y_1=1$, this follows by induction from the facts that
$Y_i\in q_v^{\frac{1}{\lcm {\mb s}}\bZ}$ and that $Y_{i+1}>Y_i$.

\medskip
iii) Let $\alpha=\big\lfloor\log_{q_v}(q^{-\frac{\deg v + g-1}
  {\min\mb{s}}} \, Y_{i+1}) \big\rfloor - 1$, which satisfies $\alpha>
\frac{1}{\min\mb{r}}+\frac{g-1}{(\min\mb{r})(\deg v)}$ if $i$ is
large enough, by Assertion (2) ii). By Corollary \ref{coro:Dirichlet},
there exists $\mb{y}\in R_v^{\,n}-\{0\}$ such that
$$
\|\,\mb{y}\,\|_{\mb{s}}\leq
q^{\frac{\deg v + g-1}{\min\mb{s}}} \,q_v^{\alpha}
<  q^{\frac{\deg v + g-1}{\min\mb{s}}}\,
q_v^{\log_{q_v} (q^{-\frac{\deg v + g-1} {\min\mb{s}}}\,Y_{i+1})} = Y_{i+1}
$$
and
\begin{align*}
\idist{A\,\mb{y}}_{\mb{r}}&
\leq q^{\frac{\deg v + g-1}{\min\mb{r}}}\,q_v^{-\alpha}
\\& \leq q^{\frac{\deg v + g-1}{\min\mb{r}}}\,
q_v^{-\big(\log_{q_v} (q^{-\frac{\deg v + g-1}{\min\mb{s}}}\,Y_{i+1})-2\big)}
=q^{(\deg v + g-1)\big(\frac{1}{\min\mb{r}}+\frac{1}{\min\mb{s}}\big)+2\deg v}
\,(Y_{i+1})^{-1}\;.
\end{align*}
Since $M_i\leq\min\{\; \idist{A\,\mb{y}}_{\mb{r}} : \mb{y}\in R_v^{\,n},\;
0<\|\,\mb{y}\,\|_{\mb{s}} <Y_{i+1}\,\}$ by the definition of a best
approxi\-mation sequence, the result follows.
%
\cqfd

\subsection{Transference theorems with weights}
\label{subsec:transference}

In this section, we will show that a matrix $A\in \cM_{m,n}(K_v)$ is
singular on average if and only if its transpose $^t\!A$ is singular
on average. To do this, following \cite[Chap.~V]{Cassels57}, we prove
a transference principle between two problems of homogeneous
approximations with weights. See also \cite{GerEvd15,German20} in the
disjoint case of the field $\bQ$.

Let $d\in\bZ_{\geq 2}$ be a positive integer at least $2$. For all
$\bm{\xi}= (\xi_1,\dots,\xi_d)$ and
$\bm{\theta}=(\theta_1,\dots,\theta_d)$ in $K_v^{\,d}$, we denote
$$
\bm{\xi}\cdot\bm{\theta}= \sum_{k=1}^d
\xi_k\;\theta_k\;.
$$
Let $\alpha_1,\dots,\alpha_d\in\bZ$ be integers and let $\alpha=
\sum_{k=1}^d\alpha_k$. We consider the parallelepiped
$$
\cP=\big\{\bm{\xi}=(\xi_1,\dots,\xi_d)\in K_v^{\,d}:
\forall\; k=1,\dots,d,\;\;|\,\xi_k\,|\leq q_v^{\alpha_k}\big\}\;.
$$
Following Schmidt's terminology \cite[page 109]{Schmidt80} in the case
of the field $\bQ$ (building on Mahler's compound one), we call the
parallelepiped
$$
\cP^{*}=\big\{\bm{\xi}=(\xi_1,\dots,\xi_d)\in K_v^{\,d}:
\forall\; k=1,\dots,d,\;\;|\,\xi_k\,|\leq
\frac{1}{q_v^{\alpha_k}}\prod_{i=1}^d q_v^{\alpha_i}=q_v^{\alpha-\alpha_k}\big\}
$$
the {\it pseudocompound} of $\cP$. Note that $\cP$ and $\cP^*$ are
preserved by the multiplication of the components of their elements by
elements of $\cO_v$.

\btheo\label{theo:transference}
With $\cP$ and $\cP^{*}$ as above, for every $F\in \SL_{d}(K_v)$, 
$$
\text{\it if}\;\;\;\cP^{*}\cap {^{t}F^{-1}}(R_v^{\,d}) \neq \{0\},
\;\;\;\text{\it then}\;\;\; \pi_{v}^{-\beta_d}\cP\cap F(R_v^{\,d}) \neq\{0\},
$$
where
$$
\beta_d = \Big\lceil \frac{1}{d-1}\Big(d+1 +
\frac{(g-1)d}{\dg v}\Big) \Big\rceil\;.
$$
\etheo

\rem The $R_v$-lattice ${^{t}F^{-1}}(R_v^{\,d})$ is called the {\it dual
  lattice} of the $R_v$-lattice $F(R_v^{\,d})$ since we have
$\mb{z}\cdot\mb{w}\in R_v$ for all $\mb{z}\in {^{t}F^{-1}}(R_v^{\,d})$
and $\mb{w}\in F(R_v^{\,d})$. They have the same covolume as $R_v^{\,d}$,
since $\det(F)=1$.

\medskip
\dem
Let $\mb{z}=(z_1,\dots, z_d)\in\cP^{*}\cap {^{t}F^{-1}}(R_v^{\,d})
-\{0\}$ and $\kappa_0=\max\{k\in\bZ_{\geq 0} : \mb{z}\in \pi_{v}^{k}
\cP^{*}\}$.  Up to permuting the coordinates, we may assume that,
for all $k=2,\dots,d$, we have
\begin{equation}\label{eq:upperbound}
|\,z_1\,| = q_v^{\alpha-\alpha_1-\kappa_0}\;\;\;{\rm
and} \;\;\; |\,z_k\,| \leq q_v^{\alpha-\alpha_k -\kappa_0}\;.
\end{equation}

With $F_k$ the $k$-th row of $F$, let us consider the $R_v$-lattice
$\Lambda = M(R_v^{\,d})$ where
$$
M = \begin{pmatrix} \pi_v^{-1}\sum_{k=1}^{d} z_k F_k \\
\pi_v^{\beta_d +\alpha_2}F_2 \\ \vdots \\
\pi_v^{\beta_d +\alpha_d}F_d \end{pmatrix} \;.
$$
By subtracting to the first row a linear combination of the other
rows, and since $\det F=1$, the determinant of the above matrix $M$ is
equal to $\pi_v^{(d-1)\beta_d +\alpha-\alpha_1-1}\;z_1$. By Equations
\eqref{eq:upperbound} and \eqref{eq:covolRv}, we thus have
$$
\covol(\Lambda)=\det(M)\covol(R_v^{\,d})=
q_v^{1-\kappa_0-(d-1)\beta_d}\;q^{(g-1)d}\;.
$$
Since $d\geq 2$ and $\beta_d \geq \frac{1}{d-1}\big(d+1 +
\frac{(g-1)d}{\deg v} \big)$, Corollary \ref{coro:Mink1} applied to
the $R_v$-lattice $\Lambda$ gives that
$$
\lambda_1(\Lambda)\leq q_v \covol(\Lambda)^{\frac{1}{d}}\leq 1.
$$
Hence, by the definition of the first minimum $\lambda_1(\Lambda)$,
there exists $\mb{w}\in R_v^{\,d}-\{0\}$ such that for every $k=2,\dots,
d$, we have
\begin{equation}\label{eq:boundF}
|\,\mb{z}\cdot F(\mb{w})\,| \leq q_v^{-1}<1\;\;\;{\rm and}
\;\;\;|\,F_k(\mb{w})\,| \leq q_v^{\beta_d +\alpha_k}\;.
\end{equation}
Since $\mb{z}\in{^{t}F^{-1}}(R_v^{\,d})$ and $\mb{w}\in R_v^{\,d}$, we have
$\mb{z}\cdot F(\mb{w})\in R_v$ by the above Remark. The first
inequality of Equation \eqref{eq:boundF} hence implies that
$\mb{z}\cdot F(\mb{w})=0$, which means that
$$
z_1 F_1 (\mb{w})=-\sum_{k=2}^{d} z_k F_k(\mb{w})\;.
$$
By the ultrametric property of $|\cdot|$, by Equations
\eqref{eq:upperbound} and \eqref{eq:boundF}, we have
\begin{align*}
q_v^{\alpha-\alpha_1-\kappa_0} |\,F_1 (\mb{w})\,| &= |\,z_1 F_1 (\mb{w})\,|
\leq \max_{2\leq k\leq d} |\,z_k F_k (\mb{w})\,| \\
&\leq \max_{2\leq k\leq d}q_v^{\alpha-\alpha_k -\kappa_0}q_v^{\beta_d + \alpha_k}
= q_v^{\alpha+\beta_d-\kappa_0 }\;.
\end{align*}
Therefore $|\,F_1 (\mb{w})\,| \leq q_v^{\beta_d+\alpha_1}$ and with the
second inequality of Equation \eqref{eq:boundF}, we conclude that
$F(\mb{w})\in\pi_v^{\beta_d}\cP$. \cqfd

\bcoro\label{coro:TranCor1} There exist $\kappa_1,\kappa_2,\kappa_3,
\kappa_4 \geq 0$ with $\kappa_2>0$, depending only on $m$, $n$, $g$,
$\deg v$, $\mb{r}$ and $\mb{s}$, such that for all $A\in \cM_{m,n}(K_v)$
and $\epsilon\in q_v^{\bZ_{\leq -1}}$, for every large enough $Y\in
q_v^{\bZ_{\geq 1}}$, if there exists $\mb{y}\in R_v^{\,n}-\{0\}$ such that
\begin{equation}\label{eq:Asol}
  \idist{A\mb{y}}_{\mb{r}}\leq \epsilon\; Y^{-1} \quad \text{and}\quad
  \|\,\mb{y}\,\|_{\mb{s}} \leq Y\;,
\end{equation}
then there exists $\mb{x}\in R_v^{\,m}-\{0\}$ such that 
\begin{equation}\label{eq:TranAsol}
  \idist{{\,^{t}\!A}\mb{x}}_{\mb{s}} \leq
  q_v^{\kappa_1}\,\epsilon^{\kappa_2}\,X^{-1} \quad\text{and}\quad
  \|\,\mb{x}\,\|_{\mb{r}}  \leq X\;,
\end{equation}
where $X=q_v^{\kappa_3}\,\epsilon^{-\kappa_4}\,Y$.
\ecoro

\dem
Let $|\mb{s}|=\sum_{j=1}^n s_j$. Denoting $\alpha_\epsilon
=-\log_{q_v}\epsilon\in\bZ_{\geq 1}$ and $\alpha_Y =\log_{q_v}Y
\in\bZ_{\geq 1}$, we define $\delta= q_v^{-\alpha_\delta}$ and
$Z=q_v^{\alpha_Z}\;Y$ where
\begin{equation}\label{eqtwoineq}
\alpha_\delta=\Big\lfloor\frac{\alpha_\epsilon-1}
{|\mb{s}|\big(\frac{1}{\min\mb{r}}+\frac{1}{\min\mb{s}}\big)-1}
\Big\rfloor
{\rm ~~~and~~~} \alpha_Z=\Big\lceil\Big(\frac{|\mb{s}|}
{\min\mb{s}}-1\Big)\alpha_\delta\Big\rceil\;.
\end{equation}
Note that $\alpha_\delta$ is well defined since
$\frac{|\mb{s}|}{\min\mb{s}} \geq 1$, and that $\alpha_\delta$ and
$\alpha_Z$ are nonnegative. We have
\begin{align}
&\;\;\big(|\mb{s}|\big(\frac{1}{\min\mb{r}}+\frac{1}{\min\mb{s}}\big)
-1\big)\,\alpha_\delta\leq \alpha_\epsilon-1\;,\nonumber\\{\rm hence} & \;\;
\big(\frac{|\mb{s}|}{\min\mb{s}}-1\big)\,\alpha_\delta +1 \leq
\alpha_\epsilon-\frac{|\mb{s}|}{\min\mb{r}}\,\alpha_\delta\;,\nonumber
\\{\rm therefore} & \;\;
\big(\frac{|\mb{s}|}{\min\mb{s}}-1\big)\,\alpha_\delta \leq
\alpha_Z\leq \alpha_\epsilon-\frac{|\mb{s}|}{\min\mb{r}}\,\alpha_\delta\;.
\label{eq:controlaphasub}
\end{align}

Let $d=m+n\geq 2$. Let us consider the following parallelepipeds
\begin{align*}
\cQ&=\left\{\mb{\xi}=(\xi_1,\dots,\xi_d)\in K_v^{\,d}:
  \begin{array}{l}\forall\; i=1,\dots,m,\;\;
    |\,\xi_i\,| \leq \epsilon^{r_i}\,Y^{-r_i}\\\forall\;j=1,\dots,n,\;\;
    |\,\xi_{m+j}\,| \leq Y^{s_j}
  \end{array} \right\}\;,\\
\cP&=\left\{\mb{\xi}=(\xi_1,\dots,\xi_d)\in K_v^{\,d}:
\begin{array}{l} \forall\; i=1,\dots,m,\;\;|\,\xi_i\,| \leq Z^{r_i} \\
  \forall\;j=1,\dots,n,\;\; |\,\xi_{m+j}\,| \leq \delta^{s_j}\,Z^{-s_j}
\end{array} \right\}\;.
\end{align*}
Since $\sum_{i=1}^m r_i=\sum_{j=1}^n s_j$, the pseudocompound $\cP^*$ of
$\cP$ is equal to
$$
\cP^{*}=\left\{\mb{\xi}=(\xi_1,\dots,\xi_d)\in K_v^{\,d}:
\begin{array}{l}\forall\; i=1,\dots,m,\;\;|\,\xi_i\,| \leq
  \delta^{|\mb{s}|} \,Z^{-r_i} \\\forall\;j=1,\dots,n,\;\;|\,\xi_{m+j}\,| \leq
  \delta^{|\mb{s}|-s_j}Z^{s_j} \end{array}  \right\}\;.
$$
By the right inequality of Equation \eqref{eq:controlaphasub}, for
every $i=1,\dots,m$, we have
$$
\delta^{|\mb{s}|} Z^{-r_i}=q_v^{-|\mb{s}|\alpha_\delta -r_i\alpha_Z}\;Y^{-r_i} \geq
q_v^{-r_i(\alpha_Z+\frac{|\mb{s}|}{\min \mb{r}}\alpha_\delta)}\;Y^{-r_i}
\geq \epsilon^{r_i}\,Y^{-r_i}\;.
$$
By the left inequality of Equation \eqref{eq:controlaphasub}, for
every $j=1,\dots,n$, we have
$$
\delta^{|\mb{s}|-s_j} Z^{s_j}=
q_v^{-(|\mb{s}|-s_j)\alpha_\delta +s_j\alpha_Z}\;Y^{s_j} \geq
q_v^{s_j(\alpha_Z-(\frac{|\mb{s}|}{\min \mb{s}}-1)\alpha_\delta)}\;Y^{s_j}
\geq Y^{s_j}\;.
$$
Therefore $\cQ$ is contained in $\cP^{*}$. 

Now, by the assumption of Corollary \ref{coro:TranCor1}, let
$\mb{y}\in R_v^{\,n}-\{0\}$ be such that the inequalities \eqref{eq:Asol}
are satisfied. Then there exists $(\mb{x}',\mb{y})\in R_v^{\,m}\times
(R_v^{\,n}-\{0\})$ such that
$$
\|A\mb{y}-\mb{x}'\|_{\mb{r}}\leq \epsilon \;Y^{-1} \quad \text{and}\quad
  \|\,\mb{y}\,\|_{\mb{s}} \leq Y\;.
$$
Therefore
$$
\cQ\cap \left(\begin{matrix} I_m & A \\ 0 & I_n \\
\end{matrix}\right) R_v^{\,d} \neq \{0\}\;.
$$
Since $\cQ\subset\cP^{*}$, this implies that
$$
\cP^{*}\cap \left(\begin{matrix} I_m & A \\ 0 & I_n \\
\end{matrix}\right) R_v^{\,d} \neq \{0\}\;.
$$
By Theorem \ref{theo:transference}, we have
$$
\pi_{v}^{-\beta_d}\cP\cap \left(\begin{matrix} I_m & 0
  \\ -{\,^{t}\!A} & I_n \\ \end{matrix}\right) R_v^{\,d} \neq \{0\}\;.
$$
Then there exists $(\mb{x},\mb{y}')\in (R_v^{\,m}\times
R_v^{\,n})-\{0\}$ such that
\begin{equation}\label{eq:controlmbx}
\|\pi_v^{\beta_d}\,\mb{x}\|_{\mb{r}} \leq Z\quad \text{and}\quad
\|\pi_v^{\beta_d}(-{\,^{t}\!A}\mb{x}-\mb{y}')\|_{\mb{s}}\leq
\delta \,Z^{-1} \;.
\end{equation}

The above inequality on the left-hand side and the two equalities of
Equation \eqref{eqtwoineq} give
\begin{align*}
\|\,\mb{x}\,\|_{\mb{r}} &\leq q_v^{\frac{\beta_d}{\min\mb{r}}}\,Z
=q_v^{\frac{\beta_d}{\min\mb{r}}+\alpha_Z}\,Y\leq
q_v^{\frac{\beta_d}{\min\mb{r}}+1+(\frac{|\mb{s}|}
{\min\mb{s}}-1)\frac{\alpha_\epsilon-1}
{|\mb{s}|(\frac{1}{\min\mb{r}}+\frac{1}{\min\mb{s}})-1}}\;Y\\ &
\leq q_v^{\frac{\beta_d}{\min \mb{r}}+1}\;
\epsilon^{-\frac{\frac{|\mb{s}|}{\min\mb{s}}-1}{
    |\mb{s}|(\frac{1}{\min\mb{r}}+\frac{1}{\min\mb{s}})-1}}\; Y\;.
\end{align*}
If $\kappa_3=\frac{\beta_d}{\min\mb{r}}+1> 0$ and $\kappa_4=
\frac{\frac{|\mb{s}|}{\min\mb{s}}-1}{ |\mb{s}| (\frac{1}{\min\mb{r}}
  +\frac{1}{\min\mb{s}})-1}\geq 0$, this proves the right inequality in
Equation \eqref{eq:TranAsol} with $X=q_v^{\kappa_3}\,
\epsilon^{-\kappa_4}\,Y$.

The right inequality in Equation \eqref{eq:controlmbx}, since
$\beta_d\geq 0$ and by using the left inequality in Equation
\eqref{eq:controlaphasub} and the definition \eqref{eqtwoineq} of
$\alpha_\delta$, gives
\begin{align*}
  \idist{{^{t}\!A}\mb{x}}_{\mb{s}} &\leq
  q_v^{\frac{\beta_d}{\min\mb{s}}}\,\delta \,Z^{-1}
  =q_v^{\frac{\beta_d}{\min\mb{s}}-\alpha_\delta-\alpha_Z} \,Y^{-1}
  \leq q_v^{\frac{\beta_d}{\min\mb{s}}-\frac{|\mb{s}|}{\min\mb{s}}
    \alpha_\delta +\kappa_3} \,\epsilon^{-\kappa_4}\,X^{-1}\\ & \leq
  q_v^{\frac{\beta_d}{\min\mb{s}}-\frac{|\mb{s}|}{\min\mb{s}}
  \big(\frac{\alpha_\epsilon-1}
      {|\mb{s}|\big(\frac{1}{\min\mb{r}}+\frac{1}{\min\mb{s}}\big)-1}
      -1\big) +\kappa_3 +\frac{\frac{|\mb{s}|}{\min\mb{s}}-1}{
        |\mb{s}|(\frac{1}{\min\mb{r}}+\frac{1}{\min\mb{s}})-1}
      \alpha_\epsilon}\,X^{-1}\\ &
  =q_v^{\frac{\beta_d}{\min\mb{s}}+\frac{|\mb{s}|}{\min\mb{s}}
  \big(\frac{1}
  {|\mb{s}|\big(\frac{1}{\min\mb{r}}+\frac{1}{\min\mb{s}}\big)-1}
  +1\big) +\kappa_3}\,\epsilon^{\frac{1}{
  |\mb{s}|(\frac{1}{\min\mb{r}}+\frac{1}{\min\mb{s}})-1}}\,X^{-1}\;.
\end{align*}
This proves the left inequality in Equation \eqref{eq:TranAsol} for
appropriate positive constants $\kappa_1$ and $\kappa_2$.

If $\mb{x}=0$, then we have $\mb{y}'\neq 0$ and $\|\,\mb{y}'\,\|_{\mb{s}}
\leq q_v^{\kappa_1-\kappa_3}\, \epsilon^{\kappa_2+\kappa_4} \,Y^{-1}$,
which contradicts the fact that $\mb{y}'\in R_v^{\,n}$ if $Y$ is large
enough.  This concludes the proof of Corollary \ref{coro:TranCor1}.
\cqfd

\bcoro \label{coro:AsingifftAsing} Let $m,n$ be positive integers and
$A\in \cM_{m,n}(K_v)$. Then $A$ is $(\mb{r},\mb{s})$-singular on
average if and only if $\,^{t}\!A$ is $(\mb{s},\mb{r})$-singular on
average.
\ecoro

\dem This follows from Corollary \ref{coro:TranCor1}.
\cqfd

\medskip
It follows from this corollary and from Remark \ref{rem:noncomplirrat}
that if $A\in\cM_{m,n}(K_v)$ is such that $\,^{t}\!A$ is not
completely irrational, then $A$ is $(\mb{r},\mb{s})$-singular on
average.

\section{Characterisation of singular on average property}
\label{sect:characsing}

In this section, we give a characterisation of the singular on average
property with weights in terms of an asymptotic property in average of
the best approximation sequence with weights. In the real case, the
relation between the singular property and the best approximation
sequence has been studied in \cite{Cheung11, Chevallier13, CheChe16,
  LiaShiSolTam20}.  Also in the real case, and with weights, the
relation (similar to the one below) between the singular on average
property and the best approximation sequence has been studied in
\cite[Prop.~6.7]{KimKimLim21}.

For the sake of later applications, we work with transposes of
matrices.

\btheo\label{theo:sing} Let $A\in \cM_{m,n}(K_v)$ and let
$(\mb{y}_i)_{i\geq 1}$ be a best approximation sequence in $K_v^{\,m}$
for $\;^{t}\!A$ with weights $(\mb{s},\mb{r})$. The following
statements are equivalent.
\begin{enumerate}
\item\label{S1} For all $a>1$ and $\epsilon
>0$, we have
$$
\lim_{N\to\infty}\frac{1}{N}\;\card\{\,\ell\in\{1,\dots,N\}: \exists\;
\mb{y}\in R_v^{\,m} - \{0\},\; \idist{^{t}\!A\,\mb{y}}_{\mb{s}}\leq \epsilon\;
a^{-\,\ell},\;
\|\,\mb{y}\,\|_{\mb{r}}\leq a^{\ell}\,\}=1\;.
$$
\item\label{S2}
  The matrix $\;^{t}\!A$ is $(\mb{s},\mb{r})$-singular on average.
\item\label{S3} There exists $a>1$ such that for every $\epsilon >0$,
  we have
$$
\lim_{N\to\infty}\frac{1}{N}\;\card\{\,\ell\in\{1,\dots,N\}: \exists\;
\mb{y}\in R_v^{\,m} - \{0\},\; \idist{{^{t}\!A}\,\mb{y}}_{\mb{s}}\leq \epsilon\;
a^{-\,\ell},\;
 \|\,\mb{y}\,\|_{\mb{r}}\leq a^{\ell}\,\}=1\;.
$$
\item\label{S4} For every $\epsilon'>0$, we have
$$
  \lim_{k\to\infty}\;\frac{1}{\log_{q_v} Y_{k}}\;\card\big\{i\leq k:
  M_{i}\,Y_{i+1}>\epsilon'\big\}=0\;.
$$
\end{enumerate}
\etheo

\dem Since Assertion \eqref{S2} is Assertion \eqref{S1} for $a=q_v>1$,
it is immediate that \eqref{S1} implies \eqref{S2} implies \eqref{S3}.

Let us first  prove that Assertion \eqref{S3} implies Assertion
\eqref{S4}. Let $a>1$ be as in Assertion \eqref{S3} and let
$\epsilon'\in\;]0,1[\,$. Let $\epsilon=\frac{\epsilon'}{a}>0$.

We may assume that the set $I= \{i\in\bZ_{\geq 1}: M_{i}Y_{i+1}
>\epsilon'\} $ is infinite, otherwise Assertion \eqref{S4} is clear
since $\lim_{k\to\infty}Y_k=+\infty$. We consider the increasing
sequence $(i_j)_{j\in\bZ_{\geq 1}}$ of positive integers such that
$I=\{i_j: j\geq 1\}$.  For every $j\geq 1$, by taking the logarithm in
base $a$, we thus have $\log_{a}{\epsilon'} -\log_{a}{M_{i_j}}
<\log_{a} Y_{i_j+1}$, hence
\begin{equation}\label{eq:isubk}
\log_{a}{\epsilon}-\log_{a}{M_{i_j}}<\log_{a} Y_{i_j+1}-1\;.
\end{equation}

Note that for every $i\geq 1$ and $X\in [Y_{i}, Y_{i+1}[\,$, the
system of inequalities
\begin{equation}\label{eq:inequal}
  \idist{{^{t}\!A}\,\mb{y}}_{\mb{s}} \leq \epsilon\, X^{-1} \quad
  \text{and}\quad 0 < \|\,\mb{y}\,\|_{\mb{r}} \leq X
\end{equation}
has a solution $\mb{y}\in R_v^{\,m}$ if and only if $M_i \leq
\epsilon\, X^{-1}$. Indeed, if the latter inequality is satisfied, then
$\mb{y}_i$ is a solution of the system \eqref{eq:inequal} since
$M_i=\idist{{^{t}\!A}\,\mb{y}_i}_{\mb{s}}$ and $X\geq Y_i=
\|\,\mb{y}_i\,\|_{\mb{r}}$. Conversely, if this system has a solution,
then since
$$
M_i\leq\min\{\; \idist{\,^{t}\!A\,\mb{y}}_{\mb{s}} :
\mb{y}\in R_v^{\,m},\; 0< \|\,\mb{y}\,\|_{\mb{r}}<Y_{i+1}\,\}
$$
by the definition of a best approximation sequence, the inequality
$M_i \leq \epsilon\, X^{-1}$ holds since $X< Y_{i+1}$.  Hence, for
every integer $\ell\in [\log_{a}{Y_i}, \log_{a}{Y_{i+1}}[\,$, the
system of inequalities \eqref{eq:inequal} has no integral
solutions for $X=a^\ell$ if and only if
\begin{equation}\label{eq:NoSol}
\log_{a}{\epsilon}-\log_{a}{M_i} <\ell<\log_{a}{Y_{i+1}}\;. 
\end{equation}

There exists an integer $j_0\geq 1$ such that for every integer $j\geq
j_0$, we have $\log_{a} Y_{i_j+1} \geq 2$ by Lemma
\ref{lem:bestapprox} (2) ii). If $\ell$ is the integer in the interval
$[\log_{a} Y_{i_j+1}-1,\log_{a} Y_{i_j+1}[$ (which is half-open and
has length $1$, hence does contain one and only one integer), then
$\ell\geq 1$ and by Equations \eqref{eq:isubk} and \eqref{eq:NoSol},
the system \eqref{eq:inequal} has no integral solutions for $X=a^\ell$.

Let $u=\lceil (\lcm \mb{r})(\log_{q_v} a)\rceil$, which belongs to
$\bZ_{\geq 1}$. By Lemma \ref{lem:bestapprox} (2) ii), for every $k\in
\bZ_{\geq 1}$, since the sequence $(i_j)_{j\in \bZ_{\geq 1}}$ is
increasing, we have
$$
Y_{i_{k+u}+1}\geq q_v^{\frac{u}{\lcm \mb{r}}}\;Y_{i_{k}+1}
\geq a \;Y_{i_{k}+1}\;.
$$
The intervals $[\log_{a}{Y_{i_{uj}+1}} -1, \log_{a} {Y_{i_{uj}+1}}[\,$ 
and $[\log_{a}{Y_{i_{u(j+1)}+1}} -1,\log_{a} {Y_{i_{u(j+1)}+1}}[\,$ are
hence disjoint for every $j\in \bZ_{\geq 1}$. Thus, if $j$ is large
enough, with $N_j=\lceil \log_{a} {Y_{i_{uj} +1}}\rceil$, the number
$n(N_j)$ of integers $\ell\in\{1,\dots, N_j\}$ such that the system of
inequalities \eqref{eq:inequal} has no integral solutions for
$X=a^\ell$ is at least $j-j_0$.  Therefore $\frac{j-j_0}{\lceil \log_{a}
  Y_{i_{uj}+1}\rceil} \leq \frac{n(N_j)}{N_j}$ tends to $0$ as $j\ra
+\infty$, by Assertion \eqref{S3}.  This implies that
$\frac{j}{\log_{a}{Y_{i_j}}}$ tends to $0$ as $j\ra+\infty$.

For every integer $k\geq 1$, let $j(k)\geq 1$ be the unique positive
integer such that we have $i_{j(k)}\leq k <i_{j(k)+1}$, so that
$j(k)=\card\{i\leq k : M_iY_{i+1}>\epsilon'\}$. Hence, since
$(Y_{i})_{i\geq 1}$ is increasing, we have
$$
  \lim_{k\to\infty}\;\frac{1}{\log_{q_v} Y_{k}}\;\card\big\{i\leq k:
  M_{i}Y_{i+1}>\epsilon'\big\}\leq\frac{\ln q_v}{\ln a}\;
  \lim_{k\to\infty}\;\frac{j(k)}{\log_{a}{Y_{i_{j(k)}}}} =0\;,
$$
which proves Assertion \eqref{S4}.

\medskip
Let us now prove that Assertion \eqref{S4} implies Assertion
\eqref{S1}. Let $a>1$ and $\epsilon\in\;]0,1[\,$.  By Lemma
\ref{lem:bestapprox} (2) iii), let $c\geq 1$ be such that for every
$i\geq 1$, we have $M_iY_{i+1}\leq a^c$. By Equation \eqref{eq:NoSol},
since the number of integer points in an open interval is at most
equal to its length, for every $i\geq 1$, the number of integers
$\ell\in [\log_{a} {Y_{i}}, \log_{a} {Y_{i+1}}[\,$ such that the
system of inequalities \eqref{eq:inequal} has no integral
solutions for $X= a^\ell$ is at most
$$
\log_{a}Y_{i+1}-(\log_{a}\epsilon-\log_{a}M_i)=
\big(\log_{a}M_iY_{i+1}-\log_{a}\epsilon\big)\;.
$$
For every $N\geq 1$ large enough, let $k_N\geq 1$ be such that $N\in
[\log_{a} {Y_{k_N}}, \log_{a} {Y_{k_N+1}}[\,$ and let $n'(N)$ be the
number of integers $\ell\in\{1,\dots, N\}$ such that the system of
inequalities \eqref{eq:inequal} has no integral solutions for
$X=a^\ell$. Then
\begin{align*}
\frac{n'(N)}{N}&\leq \frac{1}{N}\sum_{i=1}^{k_N}
\max\big\{0,\log_{a}M_iY_{i+1}-\log_{a}\epsilon\big\}\\&
\leq\big(c-\log_{a}\epsilon\big)\;\frac{1}{\log_{a} {Y_{k_N}}}
\;\card\big\{i\leq k_N : M_{i}Y_{i+1}>\epsilon\big\}\;.
\end{align*}
This last term tends to $0$ as $N\ra+\infty$ by Assertion \eqref{S4}
applied with $\epsilon'=\epsilon$.  Therefore $\lim_{N\ra+\infty}
\frac{n'(N)}{N}=0$, thus proving Assertion \eqref{S1}.
\cqfd

\section{Full Hausdorff dimension for singular on average matrices}
\label{sect:fullhausdorff}

\subsection{Modified Bugeaud-Zhang sequences}
\label{subsect:modBJseq}

In this subsection, we construct a subsequence with controlled growth
of the best approximation sequence with weights of a matrix, assuming
that its transpose is singular on average for those weights. We use as
inspiration \cite[page 470]{BugZha19} in the special case of
$K=\bF_q(Z)$ and $v=v_\infty$, and the first claim of the proof of
\cite[Theo.~2.2]{BugKimLimRam21} in the case of the field
$\bQ$ (with characteristic zero).

\bprop\label{prop:modBZseq} Let $A\in \cM_{m,n}(K_v)$ be such that
$\;^{t}\!A$ is completely irrational and $(\mb{s},\mb{r})$-singular on
average.  Let $(\mb{y}_i)_{i\in\bZ_{\geq 1}}$ be a best approximation
sequence in $K_v^{\,m}$ for $\;^{t}\!A$ with weights
$(\mb{s},\mb{r})$, and let $c>0$ be such that $M_{i}Y_{i+1}\leq q_v^c$
for every $i\in\bZ_{\geq 1}$. For all $a>b>0$, there exists an
increasing map $\varphi: \bZ_{\geq 1} \ra \bZ_{\geq 1}$ such that
\begin{enumerate}
\item[(1)]\label{assertion1:modBZseq} for every $i\in\bZ_{\geq 1}$, we
  have
\begin{equation}\label{eq:modBZseq1}
    Y_{\varphi(i+1)}\geq q_v^b\;Y_{\varphi(i)}\quad{\rm and}\quad
    M_{\varphi(i)}\;Y_{\varphi(i+1)}\leq q_v^{b+c}\;,
\end{equation}
\item[(2)]\label{assertion2:modBZseq} we have
\begin{equation}\label{eq:modBZseq2}
\limsup_{k\to\infty}\;\frac{k}{\log_{q_v} Y_{\varphi(k)}}\leq \frac{1}{a}\;.
\end{equation}
\end{enumerate}
\eprop

\dem Let $A$, $(\mb{y}_i)_{i\in\bZ_{\geq 1}}$ and $a,b$ be as in the
statement. We start by proving a particular case, that will be useful
in two of the four cases below.

\blemm\label{lem:partcase1and3} If furthermore we have $\lim_{k\to\infty}\;
Y_{k}^{\frac{1}{k}} =+\infty$, then there exists an increasing map
$\varphi: \bZ_{\geq 1} \ra\bZ_{\geq 1}$ such that Equations
\eqref{eq:modBZseq1} and \eqref{eq:modBZseq2} are satisfied.
\elemm

\dem
The fact that $\lim_{k\to\infty}\; Y_{k}^{\frac{1}{k}} =+\infty$
implies that the set
$$
\cJ_0=\{j\in\bZ_{\geq 1}:Y_{j+1}\geq q_v^b\,Y_{j}\}
$$
is infinite. We construct the increasing sequence
$(\varphi(i))_{i\in\bZ_{\geq 1}}$ of positive integers by stacks
$\{\varphi(i_k+1),\dots ,\varphi(i_{k+1})\}$ with $i_{k+1}>i_k$, by
induction on $k\in\bZ_{\geq 0}$. For $k=0$, let $i_0=0$, let $i_1=1$
and let $\varphi(1)$ be the smallest element of $\cJ_0$.

For $k\in \bZ_{\geq 0}$, assume that $i_k$ and $\varphi(i_k)$ are
constructed such that $\varphi(i_k)\in\cJ_0$ and Equation
\eqref{eq:modBZseq1} holds for every $i\leq i_k-1$. Let us construct
$i_{k+1}$ and $\varphi(i_k+1), \dots ,\varphi(i_{k+1})$ such that
$\varphi(i_{k+1}) \in\cJ_0$ and Equation \eqref{eq:modBZseq1} holds
for every $i\leq i_{k+1}-1$. Let $j_0$ be the smallest element of
$\cJ_0$ greater than $\varphi(i_k)$.  Let $r'=0$ if the set
$\{j>\varphi(i_k) : Y_{j_{0}}\geq q_v^b\,Y_{j}\}$ is empty. Otherwise,
let $r'\in\bZ_{\geq 1}$ be the maximal integer such that by induction
there exist $j_1,j_2,\dots, j_{r'}\in\bZ_{\geq 1}$ such that for
$\ell=1,\dots, r'$, the set $\{j>\varphi(i_k) : Y_{j_{\ell-1}}\geq
q_v^b\,Y_{j}\}$ is nonempty and for $\ell=1,\dots,r'+1$ the integer
$j_\ell$ is its largest element. Since the sequence
$(Y_i)_{i\in\bZ_{\geq 1}}$ is increasing, this in particular implies
that $j_{\ell-1}>j_\ell$ for $\ell=1,\dots,r'+1$, which itself ensures
the finiteness of $r'$.  Now we define $i_{k+1}=i_k+r'+1$ and
$$
\varphi(i_k+1)=j_{r'},\;\varphi(i_k+2)=j_{r'-1},\;\dots,\;
\varphi(i_k+r')=j_1,\;\varphi(i_{k+1})=j_0\;.
$$
By construction, for $\ell=1,\dots,r'$, we have
$$
Y_{\varphi(i_k+\ell+1)}=Y_{j_{r'-\ell}}\geq q_v^b\;Y_{j_{r'-\ell+1}}=
q_v^b\;Y_{\varphi(i_k+\ell)}\;.
$$
As $\varphi(i_k+1)=j_{r'}>\varphi(i_k)$, we have $Y_{\varphi(i_k+1)} \geq
Y_{\varphi(i_k)+1}\geq q_v^b\; Y_{\varphi(i_k)}$ since $\varphi(i_k)
\in\cJ_0$. Note that $\varphi(i_{k+1})=j_0 \in\cJ_0$. This proves the
claim on the left hand side of Equation \eqref{eq:modBZseq1} for
$i\leq i_{k+1}-1$.

By the maximality property of $j_{r'-\ell}$ in the above construction,
for every $\ell=1,\dots,r'$, we have $ Y_{\varphi(i_k+\ell+1)}=
Y_{j_{r'-\ell}}<q_v^b\;Y_{j_{r'-\ell+1}+1}=q_v^b\;Y_{\varphi(i_k+\ell)+1}$. By
the maximality of $r'$ in the above construction, we have
$Y_{\varphi(i_k+1)} < q_v^b\; Y_{\varphi(i_k)+1}$. Hence, by the
definition of $c$, for every $\ell=0,\dots, r'$, we have
$$
M_{\varphi(i_k+\ell)}Y_{\varphi(i_k+\ell+1)}
\leq M_{\varphi(i_k+\ell)}\; Y_{\varphi(i_k+\ell)+1}
\;q_v^{b} \leq q_v^{b+c}\;.
$$
This proves the claim on the right hand side of Equation
\eqref{eq:modBZseq1} for $i\leq i_{k+1}-1$.

Since $\lim_{k\to\infty}\; \frac{k}{\log_{q_v} Y_{k}} =0$, Equation
\eqref{eq:modBZseq2} is satisfied for $\varphi$, and this concludes the
proof of Lemma \ref{lem:partcase1and3}.  \cqfd

\medskip
Now in what follows, we will discuss four cases on the
configuration in $\bZ_{\geq 1}$ of the set
$$
\cJ=\{j\in\bZ_{\geq 1}: M_{j}\;Y_{j+1}\leq
q_v^{b+c-3a}\}\;.
$$
By Theorem \ref{theo:sing} \eqref{S4} applied with
$\epsilon'=q_v^{b+c-3a}$, we have
\begin{equation}\label{eq:limq}
  \lim_{k\to\infty}\;\frac{1}{\log_{q_v} Y_{k}}\;\card\big\{i\leq k:
  i\in \;^c\!\!\cJ\big\}=0\;.
\end{equation}

\noindent {\bf Case 1. } Assume first that $\cJ$ is finite.

By Equation \eqref{eq:limq}, we then have $\lim_{k\to\infty}\;
\frac{k}{\log_{q_v} Y_{k}} =0$, hence Proposition \ref{prop:modBZseq}
follows from Lemma \ref{lem:partcase1and3}.

\medskip\noindent {\bf Case 2. } Let us now assume that there exists
$j_*\in\bZ_{\geq 1}$ such that $j\in \cJ$ for every $j\geq j_*$.

Let us consider the auxiliary increasing sequence
$(\psi(i))_{i\in\bZ_{\geq 1}}$ of positive integers defined by
induction by setting $\psi(1)=\min\{j_*\in\bZ_{\geq 1} :\;\forall
\;j\geq j_*,\;j\in \cJ\}$ and, for every $i\geq 1$,
$$
\psi(i+1)=\min\{j\in\bZ_{\geq 1}:q_v^a \;Y_{\psi(i)}\leq Y_j\}\;.
$$
Since the sequence $(Y_i)_{i\in\bZ_{\geq 1}}$ is increasing and
converges to $+\infty$,
this is well defined, and $\psi$ is increasing, hence takes value in
$\J$ by the assumption of Case 2.
Let us now define the sequence $(\varphi(i))_{i\in\bZ_{\geq 1}}$ by, for
every $i\in\bZ_{\geq 1}$,
$$
\varphi(i)=\left\{\begin{array}{ll}
\psi(i) & {\rm if~~} M_{\psi(i)}\; Y_{\psi(i+1)}\leq
q_v^{b+c-a}\;,
\\\psi(i+1)-1 & {\rm otherwise.}
\end{array} \right.
$$
Note that the sequence $(\varphi(i))_{i\in\bZ_{\geq 1}}$ is increasing
with $\varphi\geq \psi$.

Let $i\in\bZ_{\geq 1}$. Let us prove that
\begin{equation}\label{eq:modBZseq3}
Y_{\varphi(i+1)}\geq q_v^a\,Y_{\varphi(i)}\quad{\rm and}\quad
M_{\varphi(i)}Y_{\varphi(i+1)}\leq q_v^{b+c}\;,
\end{equation}
by discussing on the values of $\varphi(i)$ and $\varphi(i+1)$. This
implies that Equation \eqref{eq:modBZseq1} is satisfied since $a\geq
b$, and that Equation \eqref{eq:modBZseq2} is satisfied since by induction
$Y_{\varphi(k)}\geq q_v^{a\,(k-1)}\;Y_{\varphi(1)}$ for every $k\in\bZ_{\geq 1}$.

\medskip
$\bullet$~ Assume that $\varphi(i)=\psi(i)$ and $\varphi(i+1)
=\psi(i+1)$. By the definition of $\psi(i+1)$, we have
$$
Y_{\varphi(i+1)}=Y_{\psi(i+1)}\geq q_v^a \;Y_{\psi(i)}=
q_v^a \;Y_{\varphi(i)}\;.
$$
If $\psi(i)\neq \psi(i+1)-1$, then by the definition of $\varphi(i)$,
we have
$$
M_{\varphi(i)}\; Y_{\varphi(i+1)}=
M_{\psi(i)}\; Y_{\psi(i+1)}\leq q_v^{b+c-a}\leq q_v^{b+c}\;.
$$
If $\psi(i)= \psi(i+1)-1$, then $\varphi(i+1)= \varphi(i)+1$ and by
the definition of $c$, we have
$$
M_{\varphi(i)}\; Y_{\varphi(i+1)}=
M_{\varphi(i)}\; Y_{\varphi(i)+1}\leq q_v^{c}\leq q_v^{b+c}\;.
$$
This proves Equation \eqref{eq:modBZseq3}.

\medskip
$\bullet$~ Assume that $\varphi(i)=\psi(i)$ and $\varphi(i+1)=
\psi(i+2)-1$. Since the sequence $(Y_i)_{i\in\bZ_{\geq 1}}$ is increasing
and by the definition of $\psi(i+1)$, we have
$$
Y_{\varphi(i+1)}=Y_{\psi(i+2)-1}\geq Y_{\psi(i+1)}\geq q_v^a
\;Y_{\psi(i)}= q_v^a \;Y_{\varphi(i)}\;.
$$
We have $q_v^a \;Y_{\psi(i+1)}> Y_{\psi(i+2)-1}$ by the minimality
property of $\psi(i+2)$. If $\psi(i+1)>\psi(i)+1$, then $M_{\psi(i)}\;
Y_{\psi(i+1)}\leq q_v^{b+c-a}$ by the dichotomy in the definition of
$\varphi(i)$. Hence
$$
M_{\varphi(i)}\; Y_{\varphi(i+1)}= M_{\psi(i)}\; Y_{\psi(i+2)-1}\leq
M_{\psi(i)}\; Y_{\psi(i+1)}\,q_v^{a}\leq q_v^{b+c-a}q_v^{a}= q_v^{b+c}\;.
$$
If $\psi(i+1)=\psi(i)+1$, then $M_{\psi(i)}\; Y_{\psi(i)+1} \leq
q_v^{b+c-3a}$ since $\psi(i)\in\J$. Hence
$$
M_{\varphi(i)}\; Y_{\varphi(i+1)}= M_{\psi(i)}\; Y_{\psi(i+2)-1}\leq
M_{\psi(i)}\; Y_{\psi(i)+1}\,q_v^{a}\leq
q_v^{b+c-3a}q_v^{a}\leq q_v^{b+c}\;.
$$
This proves Equation \eqref{eq:modBZseq3}.

\medskip
$\bullet$~ Assume that $\varphi(i)=\psi(i+1)-1$ and $\varphi(i+1)=
\psi(i+1)$. Since $\psi(i+1)-1\in\cJ$, we have
$$
M_{\varphi(i)}Y_{\varphi(i+1)}=
M_{\psi(i+1)-1}Y_{\psi(i+1)}\leq
q_v^{b+c-3a}\leq q_v^{b+c}\;.
$$
If $\psi(i+1)-1=\psi(i)$, then by the definition of $\psi(i+1)$, we
have
$$
\frac{Y_{\varphi(i+1)}}{Y_{\varphi(i)}}=\frac{Y_{\psi(i+1)}}{Y_{\psi(i+1)-1}}
=\frac{Y_{\psi(i+1)}}{Y_{\psi(i)}}\geq q_v^a\;.
$$
If $\psi(i+1)-1>\psi(i)$, then we have $M_{\psi(i)}\;
Y_{\psi(i+1)}>q_v^{b+c-a}$ by the dichotomy in the definition of
$\varphi(i)$, we have $Y_{\psi(i+1)-1} <q_v^a \;Y_{\psi(i)}\leq q_v^a
\;Y_{\psi(i)+1}$ by the minimality property of $\psi(i+1)$, and we
have $M_{\psi(i)}\; Y_{\psi(i)+1} \leq q_v^{b+c-3a}$ since
$\psi(i)\in\cJ$. Therefore
$$
\frac{Y_{\varphi(i+1)}}{Y_{\varphi(i)}} =
\frac{Y_{\psi(i+1)}}{Y_{\psi(i+1)-1}}=
\frac{M_{\psi(i)}\;Y_{\psi(i+1)}}{M_{\psi(i)}\;Y_{\psi(i+1)-1}}
\geq\frac{q_v^{b+c-a}}{M_{\psi(i)}\; Y_{\psi(i)+1}\;q_v^{a}}
\geq\frac{q_v^{b+c-a}}{q_v^{b+c-3a}\; q_v^{a}}=q_v^a\;.
$$
This proves Equation \eqref{eq:modBZseq3}.

\medskip
$\bullet$~ Assume that $\varphi(i)=\psi(i+1)-1$ and $\varphi(i+1)=
\psi(i+2)-1$.
By the previous case computations, we have
$$
\frac{Y_{\varphi(i+1)}}{Y_{\varphi(i)}}
=\frac{Y_{\psi(i+2)-1}}{Y_{\psi(i+1)-1}}
\geq\frac{Y_{\psi(i+1)}}{Y_{\psi(i+1)-1}}\geq q_v^a\;.
$$
We have $q_v^a \;Y_{\psi(i+1)}> Y_{\psi(i+2)-1}$ by the minimality
property of $\psi(i+2)$.  Hence since $\psi(i+1)-1\in\cJ$, we have
\begin{align*}
M_{\varphi(i)}Y_{\varphi(i+1)}&=
M_{\psi(i+1)-1}Y_{\psi(i+2)-1}
=M_{\psi(i+1)-1}Y_{\psi(i+1)}
\Big(\frac{Y_{\psi(i+2)-1}}{Y_{\psi(i+1)}}\Big)\\ &\leq
q_v^{b+c-3a}\;q_v^{a}\leq q_v^{b+c}\;.
\end{align*}
This proves Equation \eqref{eq:modBZseq3} and concludes the proof
of Case 2.

\medskip\noindent {\bf Case 3. } Let us now assume that $\cJ$ and
$^c\!\!\cJ$ are both infinite, and that the number of sequences of
consecutive elements of $\cJ$ with length at least $3a$ is finite.

Let $j_0=\min\cJ$. Let us write the set $\bZ_{\geq j_0}=
\bigcup_{i\in\bZ_{\geq 1}} C_i\cup D_i $ as the disjoint union of
nonempty finite sequences $C_i$ of consecutive integers in $\cJ$ and
finite nonempty sequences $D_i$ of consecutive integers in $^c\!\!\cJ$
with $\max C_i< \min D_i\leq \max D_i< \min C_{i+1}$ for all $i\in
\bZ_{\geq 1}$.  Under the assumption of Case 3, let $i_0\in\bZ_{\geq
  1}$ be such that $\card\; C_i <3a$ for every $i\geq i_0$. Let
$k_0=\min C_{i_0}$.

Then there exists an element of $^c\!\!\cJ$ in any finite sequence of
$3\lceil a \rceil +1$ consecutive integers at least $k_0$, so that for
every $k\in\bZ_{\geq 1}$ we have
$$
\frac{k}{\log_{q_v} Y_{k}}\leq
\frac{k_0+(3\lceil a\rceil +1)\;\card\big\{i\leq k:
  i\in \;^c\!\!\cJ\big\}}{\log_{q_v} Y_{k}}\;,
$$
which converges to $0$ as $k\ra+\infty$ by Equation \eqref{eq:limq}
and since $\lim_{k\to\infty}\; Y_{k} =+\infty$.  Therefore
$\lim_{k\to\infty}\; Y_{k}^{\frac{1}{k}} =+\infty$, and Lemma
\ref{lem:partcase1and3} implies Proposition \ref{prop:modBZseq}.

\medskip\noindent {\bf Case 4. } Let us finally assume that $\cJ$ and
$^c\!\!\cJ$ are both infinite, and that there are infinitely many
sequences of consecutive elements of $\cJ$ with length at least $3a$.

With the notation $(C_i)_{i\in\bZ_{\geq 1}}$ and $(D_i)_{i\in\bZ_{\geq
    1}}$ of the beginning of Case 3, let $(i_k)_{k\in\bZ_{\geq 1}}$ be
the increasing sequence of positive integers such that
$\{i\in\bZ_{\geq 1}:\card\; C_i \geq 3a\}=\{i_k :k\in\bZ_{\geq 1}\}$.

For every $k\in\bZ_{\geq 1}$, let us define an increasing finite
sequence $(\psi_k(i))_{1\leq i\leq m_k+1}$ of positive integers by
setting $\psi_k(1)=\min C_{i_k}$ and by induction
$$
\psi_k(i+1)=\min \{\,j\in C_{i_k}:q_v^a\;Y_{\psi_k(i)}\leq Y_j\,\}\;,
$$
as long as this set is nonempty. Since $C_{i_k}$ is a finite sequence
of consecutive positive integers with length at least $3a$ and
$Y_{i+1}\geq q_v^{\frac{1}{\min \mb{r}}}\, Y_i$ for every
$i\in\bZ_{\geq 1}$, there exists $m_k \in\bZ_{\geq 2}$ such that
$\psi_k(i)$ is defined for $i=1,\dots, m_k+1$.  Note that $\psi_k(i)$
belongs to $\cJ$ for $i=1,\dots, m_{k+1}$ since $C_{i_k}\subset \cJ$.

As in Case 2, let us define an increasing finite sequence
$(\varphi_k(i))_{1\leq i\leq m_k}$ of positive integers by
$$
\varphi_k(i)=\left\{\begin{array}{ll}
\psi_k(i) & {\rm if~~} M_{\psi_k(i)}\; Y_{\psi_k(i+1)}\leq
q_v^{b+c-a}\;,
\\\psi_k(i+1)-1 & {\rm otherwise.}
\end{array} \right.
$$
As in the proof of Case 2, since for $i=1,\dots, m_k$, the integers
$\psi_k(i)$, $\psi_k(i+1)$ as well as $\psi_k(i+1)-1$ belong to $\cJ$,
we have, for every $i=1,\dots, m_k-1$,
\begin{equation}\label{eq:controlvarphik}
  Y_{\varphi_k(i+1)}\geq q_v^a\;Y_{\varphi_k(i)}\quad{\rm and}\quad
    M_{\varphi_k(i)}\;Y_{\varphi_k(i+1)}\leq q_v^{b+c}\;.
\end{equation}

Since $\varphi_k(m_k)\in C_{i_k}$ and $\varphi_{k+1}(1)\in
C_{i_{k+1}}$, we have $\varphi_k(m_k)<\varphi_{k+1}(1)$. Let us define
an increasing finite sequence $(\varphi'_k(i))_{1\leq i\leq r_k+1}$ of
positive integers that will allow us to interpolate between
$\varphi_k(m_k)$ and $\varphi_{k+1}(1)$. Let $j_0=\varphi_{k+1}(1)$.
If $\{j\in\bZ_{\geq \varphi_k(m_k)}: Y_{j_0}\geq q_v^b\;Y_j\}$ is
empty, let $r'_k=0$ and $\varphi'_k(1) =j_0=\varphi_{k+1}(1)$.
Otherwise, by decreasing induction, let $r'_k\in\bZ_{\geq 1}$ be the
maximal positive integer such that there exist $j_1,\dots, j_{r'_k}\in
\bZ_{\geq 1}$ such that for $\ell=1 ,\dots, r'_k$, the set
$\{j\in\bZ_{\geq \varphi_k(m_k)} : Y_{j_{\ell-1}}\geq q_v^b\;Y_j\}$ is
nonempty and for $\ell=1,\dots, r'_k+1$, the integer $j_\ell$ is its
largest element. As in the part of the proof of Case 1 that does not
need some belonging to $\cJ_0$, the sequence $(\varphi'_k(i) =
j_{r'_k+1-i})_{1\leq i\leq r'_k+1}$ is well defined, it is contained
in $[\varphi_k(m_k),\varphi_{k+1}(1)]$, and for $i=1,\dots, r'_k$, we
have
\begin{equation}\label{eq:controlvarphiprimek}
Y_{\varphi'_k(i+1)}\geq q_v^b\;Y_{\varphi'_k(i)}\quad{\rm and}\quad
M_{\varphi'_k(i)}\;Y_{\varphi'_k(i+1)}\leq q_v^{b+c}\;.
\end{equation}

Putting alternatively together the sequences $(\varphi_k(i))_{1\leq
  i\leq m_k-1}$ and $(\varphi'_k(i))_{1\leq i\leq r'_k}$ as $k$ ranges
over $\bZ_{\geq 1}$, we now define (with the standard convention that
an empty sum is zero) $N_k=\sum_{\ell=1}^{k-1}(m_\ell-1+r'_\ell)$ and
$$
\varphi(i)=\left\{\begin{array}{ll}
\varphi_k\big(i-N_k\big)&
       {\rm if~~} 1+N_k\leq i\leq
       m_k-1+N_k\\
\varphi'_k\big(i+1-m_k-N_k\big)&
       {\rm if~~} m_k+N_k\leq i\leq
       r'_k-1+m_k+N_k\;.\\
\end{array}\right.
$$
By Equation \eqref{eq:controlvarphik} for $i=1,\dots, m_k-2$, by
Equation \eqref{eq:controlvarphiprimek} for $i=1,\dots, r'_k$, and
since $\varphi'_k(r'_k+1)=\varphi_{k+1}(1)$, in order to prove that
the map $\varphi$ satisfies Equation \eqref{eq:modBZseq1}, hence
Assertion (1) of Proposition \ref{prop:modBZseq}, we only have to
prove the following lemma.

\blemm\label{lem:chainonmanquant} For every $k\in\bZ_{\geq 1}$, we have
\begin{equation}\label{eq:controlvarphikkplusun}
Y_{\varphi'_k(1)}\geq q_v^b\;Y_{\varphi_k(m_k-1)}\quad{\rm and}\quad
M_{\varphi_k(m_k-1)}\;Y_{\varphi'_k(1)}\leq q_v^{b+c}\;.
\end{equation}
\elemm

\dem
Since $\varphi'_k(1)\geq \varphi_k(m_k)$, hence $Y_{\varphi'_k(1)}
\geq Y_{\varphi_k(m_k)}$, the left hand side of Equation
\eqref{eq:controlvarphikkplusun} follows from the left hand side of
Equation \eqref{eq:controlvarphik} with $i=m_k-1$. If $\varphi'_k(1)=
\varphi_k(m_k)$, then the right hand side of Equation
\eqref{eq:controlvarphikkplusun} follows from the right hand side of
Equation \eqref{eq:controlvarphik} with $i=m_k-1$.

Let us hence assume that $\varphi'_k(1)>\varphi_k(m_k)$, so that
\begin{equation}\label{eq:hautpage20}
Y_{\varphi'_k(1)}\leq q_v^b\;Y_{\varphi_k(m_k)}\leq q_v^a\;Y_{\varphi_k(m_k)}
\end{equation}
by the maximality of $r'_k$.  Let us prove that $\varphi_k(m_k)=
\psi_k(m_k)$.  For a contradiction, assume otherwise that
$\varphi_k(m_k)= \psi_k(m_k+1)-1 >\psi_k(m_k)$. As in the third
subcase of Case 2, we have $M_{\psi_k(m_k)}\; Y_{\psi_k(m_k+1)}
>q_v^{b+c-a}$ by the dichotomy in the definition of $\varphi_k(m_k)$,
we have $Y_{\psi_k(m_k+1)-1} <q_v^a \;Y_{\psi_k(m_k)} \leq q_v^a
\;Y_{\psi_k(m_k)+1}$ by the minimality property of $\psi_k(m_k+1)$,
and we have $M_{\psi_k(m_k)}\; Y_{\psi_k(m_k)+1} \leq q_v^{b+c-3a}$
since $\psi_k(m_k)\in\cJ$. Therefore, as in the third subcase of Case
2, we have
$$
\frac{Y_{\psi_k(m_k+1)}}{Y_{\psi_k(m_k+1)-1}}=
\frac{M_{\psi_k(m_k)}\;Y_{\psi_k(m_k+1)}}{M_{\psi_k(m_k)}\;Y_{\psi_k(m_k+1)-1}}
\geq q_v^a\;.
$$
Hence by the construction of $\varphi'_k(1)$, we have
$\varphi'_k(1)=\varphi_k(m_k)$, a contradiction to our assumption
that $\varphi'_k(1)>\varphi_k(m_k)$.  We now discuss on the two
possible values of $\varphi_k(m_k-1)$.

\medskip
First assume that $\varphi_k(m_k-1)=\psi_k(m_k-1)$. If $\psi_k(m_k-1)
\neq \psi_k(m_k)-1$ then $M_{\psi_k(m_k-1)}\; Y_{\psi_k(m_k)} \leq
q_v^{b+c-a}$ by the dichotomy in the definition of $\varphi_k(m_k-1)$.
If on the contrary $\psi_k(m_k-1)= \psi_k(m_k)-1$ then $M_{\psi_k(m_k-1)}\;
Y_{\psi_k(m_k)} \leq q_v^{b+c-3a}\leq q_v^{b+c-a}$ since the integer
$\psi_k(m_k)-1$ belong to $\cJ$ as $m_k\geq 2$. Since $\varphi_k(m_k)
=\psi_k(m_k)$ by Equation \eqref{eq:hautpage20}, we have
$$
M_{\varphi_k(m_k-1)}\;Y_{\varphi'_k(1)} =M_{\psi_k(m_k-1)}\;Y_{\psi_k(m_k)}
\Big(\frac{Y_{\varphi'_k(1)}}{Y_{\varphi_k(m_k)}}\Big)\leq
q_v^{b+c-a}\;q_v^a=q_v^{b+c}\;.
$$
This proves the right hand side of Equation
\eqref{eq:controlvarphikkplusun}.

\medskip
Now assume that $\varphi_k(m_k-1)=\psi_k(m_k)-1$.  Again since
$\varphi_k(m_k)=\psi_k(m_k)$, since the integer $\psi_k(m_k)-1$
belongs to $\cJ$ as $m_k\geq 2$, and by Equation
\eqref{eq:hautpage20}, we have
$$
M_{\varphi_k(m_k-1)}\;Y_{\varphi'_k(1)} =M_{\psi_k(m_k)-1}\;Y_{\psi_k(m_k)}
\Big(\frac{Y_{\varphi'_k(1)}}{Y_{\varphi_k(m_k)}}\Big)\leq
q_v^{b+c-3a}\;q_v^a\leq q_v^{b+c}\;.
$$
This proves the right hand side of Equation
\eqref{eq:controlvarphikkplusun}, and concludes the proof of Lemma
\ref{lem:chainonmanquant}. \cqfd

\medskip
Finally, let us prove Assertion (2) of Proposition
\ref{prop:modBZseq}.  Since there exists an element of $^c\!\!\cJ$ in
any finite sequence of $3\lceil a \rceil +1$ consecutive integers in
the complement of $\bigcup_{k\in\bZ_{\geq1}} C_{i_k}$, there exists
$c_0\geq 0$ such that, for every $k\in\bZ_{\geq 1}$, we have
$$
  \frac{\card\{j\leq \varphi(k):j\notin \bigcup_{k\in\bZ_{\geq1}} C_{i_k}\}}
       {\log_{q_v} Y_{\varphi(k)}}\leq
\frac{c_0+(3\lceil a\rceil +1)\;\card\big\{j\leq \varphi(k):
  j\in \;^c\!\!\cJ\big\}}{\log_{q_v} Y_{\varphi(k)}}\;,
$$
which converges to $0$ as $k\ra+\infty$ as seen at the end of the
proof of Case 3.  Let us define $n(k)=\card\{i\leq k: Y_{\varphi(i)}
\geq q_v^a\; Y_{\varphi(i+1)}\}$.  For every $\ell\in\bZ_{\geq 1}$,
since $Y_{j+1}\geq q_v^{\frac{1}{\min\mb{r}}}\,Y_j$ for every
$j\in\bZ_{\geq 1}$, and by the maximality of $m_\ell$ in the
construction of $\big(\varphi_\ell(i)\big)_{1\leq i\leq m_\ell}$, we
have $\card\{j\in C_{i_\ell}:j\geq \varphi_\ell(m_\ell)\}\leq
2\,\lceil a\rceil\min\mb{r}$. If $\varphi(i)$ belongs to $C_{i_\ell}$
but $\varphi(i+1)$ does not, then $\varphi(i) \geq
\varphi_\ell(m_\ell)$. Since when $\varphi(i)$ and $\varphi(i+1)$
belong to $C_{i_\ell}$ for some $\ell\in\bZ_{\geq 1}$, then $\varphi$
and $\varphi_\ell$ coincide on $i$ and $i+1$, and since Equation
\eqref{eq:controlvarphik} holds, we hence have
$$
k-n(k)=\card\{i\leq k:Y_{\varphi(i)}<
q_v^a\;Y_{\varphi(i+1)}\}\leq 2\,\lceil a\rceil\;\min\mb{r}\;
\card\{j\leq \varphi(k): j\notin \bigcup_{k\in\bZ_{\geq1}} C_{i_k}\}\;.
$$
Hence
\begin{align*}
\limsup_{k\ra+\infty}\frac{k}{\log_{q_v} Y_{\varphi(k)}}&=
\limsup_{k\ra+\infty}\frac{n(k)+k-n(k)}{\log_{q_v} Y_{\varphi(k)}}
=\limsup_{k\ra+\infty}\frac{n(k)}{\log_{q_v} Y_{\varphi(k)}}\\  &\leq
\limsup_{k\ra+\infty}\frac{n(k)}{\log_{q_v} q_v^{a(n(k)-1)}Y_{\varphi(1)}}
=\frac{1}{a}\;.
\end{align*}
This proves Equation \eqref{eq:modBZseq2} and concludes the proof of
Proposition \ref{prop:modBZseq}.
\cqfd

\subsection{Lower bound on the Hausdorff dimension of $\Bad_A(\epsilon)$}
\label{subsec:lowerbound}

In this subsection, we use the scheme of proof in the real case of
\cite[Theo.~6.1]{ChoGhoGuaMarSim20}, which is a weighted version of
\cite[Theo.~5.1]{BugKimLimRam21}, in order to estimate the lower bound
on the Hausdorff dimension of the $\epsilon$-bad sets of
$(\mb{r},\mb{s})$-singular in average matrices.

For a given sequence $(\mb{y}_i)_{i\geq 1}$ in $R_v^{\,m}-\{0\}$ and
for every $\delta >0$, let
$$
\Bad_{(\mb{y}_i)_{i\geq 1}}^{\delta}=\{\bm{\theta}\in(\pi_v\cO_{\!v})^{m}:\;
\forall\;i\geq 1,\;\; |\idist{\bm{\theta}\cdot\mb{y}_i}|\geq \delta\,\}\;.
$$

\bprop\label{prop:SProp} Let $A\in \cM_{m,n}(K_v)$ be such that
$\,^{t}\!A$ is completely irrational and let $(\mb{y}_i)_{i\geq 1}$ be
a best approximation sequence in $K_v^{\,m}$ for $\,^{t}\!A$ with
weights $(\mb{s},\mb{r})$. Suppose that there exist $b,c > 0$ and an
increasing function $\varphi:\bZ_{\geq 1}\to \bZ_{\geq 1}$ such that
$$
\forall\;i\in\bZ_{\geq 1},\;\;\;M_{\varphi(i)}Y_{\varphi(i+1)}\leq q_v^{b+c}\;.
$$
Then for every $\delta\in\;]0,1]$, if $\epsilon=\delta^{\frac{1}{\min
    \mb{r}} +\frac{1}{\min \mb{s}}}\;q_v^{-b-c}$, then the set
$\Bad_{(\mb{y}_{\varphi(i)}) _{i\geq 1}}^{\delta}$ is contained in the
set $\Bad_A(\epsilon)$.
\eprop

\dem Fix $\delta\in\; ]0, 1]$ and $\bm{\theta}\in
\Bad_{(\mb{y}_{\varphi(i)})_{i\geq 1}} ^{\delta}$.  Let $\epsilon_1 =
\delta^{\frac{1}{\min \mb{s}}} \;q_{v}^{-b-c}$. For every
$(\mb{y'},\mb{x'})$ in $R_v^{\,m} \times R_v^{\,n}$ such that
$\|\,\mb{x'}\,\|_{\mb{s}}\geq\epsilon_1Y_{\varphi(1)}$, let $k$ be the
unique element of $\bZ_{\geq 1}$ for which
\[
Y_{\varphi(k)}\leq \epsilon_1^{-1} \|\,\mb{x'}\,\|_{\mb{s}} < Y_{\varphi(k+1)}\;,
\]
which exists since $\|\,\mb{x'}\,\|_{\mb{s}}\geq \epsilon_1
Y_{\varphi(1)}$ and since the sequence $(Y_{\varphi(i)})_{i\geq 1}$ is
increasing, converging to $+\infty$.  Let $\mb{x}_{\varphi(k)}\in
R_v^{\,n}$ be such that $M_{\varphi(k)}= \|{\,^{t}\!A}
\mb{y}_{\varphi(k)} -\mb{x}_{\varphi(k)}\|_{\mb{s}}$.  Then by the
ultrametric inequality, the assumption of the proposition, the fact
that $\epsilon_1\;q_v^{b+c}=\delta^{\frac{1}{\min \mb{s}}}\leq 1$ and
the definition of $\Bad_{(\mb{y}_{\varphi(i)})_{i\geq 1}}^{\delta}$,
we have
\begin{align}
|({\,^{t}\!A}\mb{y}_{\varphi(k)}-\mb{x}_{\varphi(k)})\cdot \mb{x'}| &\leq
\max_{1\leq i\leq n} M_{\varphi(k)}^{\,s_i} \|\,\mb{x'}\,\|_{\mb{s}}^{\,s_i} <
\max_{1\leq i\leq n} (\epsilon_1M_{\varphi(k)}Y_{\varphi(k+1)})^{s_i}
\nonumber\\ &\leq
(\epsilon_1\;q_v^{b+c})^{\min \mb{s}}= \delta\leq \min_{\ell'\in R_v}
|\,\mb{y}_{\varphi(k)}\cdot \bm{\theta}-\ell'\,|\;.\label{eq:stricineqdelta}
\end{align}
Observe that 
\begin{align*}
  \mb{y}_{\varphi(k)}\cdot \bm{\theta}&=
  \mb{y}_{\varphi(k)}\cdot (A\,\mb{x'})+\mb{y}_{\varphi(k)}\cdot\mb{y'}-
  \mb{y}_{\varphi(k)}\cdot (A\,\mb{x'} +\mb{y'}- \bm{\theta})\\
  &= ({\,^{t}\!A}\,\mb{y}_{\varphi(k)})\cdot\mb{x'}-\mb{x}_{\varphi(k)}\cdot\mb{x'}
  +\ell -\mb{y}_{\varphi(k)}\cdot (A\mb{x'} +\mb{y'}- \bm{\theta}),
  \end{align*}
where $\ell = \mb{x}_{\varphi(k)}\cdot\mb{x'}+\mb{y}_{\varphi(k)}\cdot
\mb{y'} \in R_v$.  Thus we have, using the equality case of the
ultrametric inequality for the second equality below with the strict
inequality in Equation \eqref{eq:stricineqdelta}, and again the definition of
$\Bad_{(\mb{y}_{\varphi(i)})_{i\geq 1}}^{\delta}$ for the last inequality below,
\begin{align*}
  |\mb{y}_{\varphi(k)}\cdot (A\mb{x'} +\mb{y'}- \bm{\theta})| &=
  |({\,^{t}\!A}\mb{y}_{\varphi(k)}-\mb{x}_{\varphi(k)})\cdot\mb{x'}-
  \mb{y}_{\varphi(k)}\cdot \bm{\theta}+\ell\,|\\ &= \max
\left\{ |({\,^{t}\!A}\mb{y}_{\varphi(k)}-\mb{x}_{\varphi(k)})\cdot\mb{x'}|,
  |\,\mb{y}_{\varphi(k)}\cdot \bm{\theta}-\ell\,| \right\}\\
  &=|\,\mb{y}_{\varphi(k)}\cdot \bm{\theta}-\ell\,|\geq
  |\idist{\mb{y}_{\varphi(k)}\cdot \bm{\theta}}|\geq \delta.
\end{align*}
Hence, we have 
$$
\delta \leq |\mb{y}_{\varphi(k)}\cdot (A\mb{x'} + \mb{y'} - \bm{\theta})|
\leq \max_{1\leq j\leq m} Y_{\varphi(k)}^{\,r_j}\;\|A\mb{x'} + \mb{y'} -
\bm{\theta}\|_{\mb{r}}^{\,r_j},
$$
which implies, since $\delta\leq 1$, that
\[
Y_{\varphi(k)}\|A\mb{x'} +\mb{y'}- \bm{\theta}\|_{\mb{r}} \geq
\min_{1\leq j\leq m}\delta^{\frac{1}{r_j}}= \delta^{\frac{1}{\min \mb{r}}}.
\]
Finally, for every $(\mb{y'},\mb{x'})$ in $R_v^{\,m} \times R_v^{\,n}$ such that
$\|\,\mb{x'}\,\|_{\mb{s}}\geq\epsilon_1Y_{\varphi(1)}$, we have
$$
\|\,\mb{x'}\,\|_{\mb{s}}\;\|A\mb{x'} +\mb{y'}- \bm{\theta}\|_{\mb{r}} \geq
\epsilon_1 \;Y_{\varphi(k)}\,\|A\mb{x'} +\mb{y'}- \bm{\theta}\|_{\mb{r}}
\geq \delta^{\frac{1}{\min \mb{r}}+\frac{1}{\min \mb{s}}}\;q_v^{-b-c}\;.
$$
By Equation \eqref{eq:1523}, this implies that $\bm{\theta}\in
\Bad_A(\epsilon)$ for $\epsilon=\delta^{\frac{1}{\min \mb{r}}
  +\frac{1}{\min \mb{s}}}\; q_v^{-b-c}$.
\cqfd

\bprop\label{prop:DProp} For every $\delta\in \;]0,\frac{1}{q_v^{\,3m}}[\,$,
there exist $b=b(\delta)>0$ and $C=C(\delta)>0$ such
that for every sequence $(\mb{y}_i)_{i\in\bZ_{\geq 1}}$ in $R_v^{\,m}
- \{0\}$ satisfying $\|\,\mb{y}_{i+1}\,\|_{\mb{r}}\geq q_{v}^{b}\;
\|\,\mb{y}_{i}\,\|_{\mb{r}}$ for all $i \in \bZ_{\geq 1}$, we have
$$
\dimH\;\Bad_{(\mb{y}_i)_{i\geq 1}}^{\delta} \geq
m - C \limsup_{k\to\infty}\frac{k}{\log_{q_v}\|\,\mb{y}_{k}\,\|_{\mb{r}}}\;.
$$
\eprop

\dem 
Fix $\delta\in\;]0,\frac{1}{q_v^{\,3m}}[\,$. Let
\begin{equation}\label{eq:defibbdelta}
b=b(\delta)=\frac{-\log_{q_v} \delta}{\min\mb{r}} \;,
\end{equation}
which is positive since $\delta<1$. By the mass distribution principle
(see for instance \cite[page 60]{Falconer14}), it is enough to prove
that there exist a (Borel, positive) measure $\mu$, supported on
$\Bad_{(\mb{y}_i)_{i\geq 1}}^{\delta}$, and constants $C, C_0, r_0
>0$, with $C$ depending only on $\delta$, such that, for every closed
ball $B$ of radius $r<r_0$, we have
$$
\mu(B)\leq C_0\; r^{m-C\limsup_{k\to \infty}\frac{k}{\log_{q_v}\|\,\mb{y}_k\,\|_\mb{r}}}\;.
$$
We adapt by modifying it quite a lot the measure construction in the
proof of \cite[Theo.~6.1]{ChoGhoGuaMarSim20}.

By convention, let $Y_0=1$ and $n_{0,j}=0$ for $j=1,\dots,m$.  For
every $k\in\bZ_{\geq 1}$, define $Y_k = \|\,\mb{y}_k\,\|_{\mb{r}}$, which
is at least $1$ since $\mb{y}_k\in R_v^{\,m} - \{0\}$, and for every
$j=1,\dots,m$, let $n_{k,j}\in \bZ_{\geq 0}$ be such that
\begin{equation}\label{eq:defnkj}
q_v^{\,-n_{k,j}}\leq Y_{k}^{\,-r_j}< q_v^{\,-n_{k,j}+1}\;.
\end{equation}
Note that the sequence $(n_{k,j})_{k\in\bZ_{\geq 0}}$ is
nondecreasing, for all $j=1,\dots,m$.

For every $k\in\bZ_{\geq 0}$, let us consider the polydisc 
$$
\Pi(Y_k)=\ov{B}(0,\frac{1}{q_v}\,Y_{k}^{-r_1})\times\cdots\times
\ov{B}(0,\frac{1}{q_v}\,Y_{k}^{-r_m})=
\ov{B}(0,q_{v}^{-n_{k,1}-1})\times\cdots\times \ov{B}(0,q_{v}^{-n_{k,m}-1})\;,
$$
where $\ov{B}(0,r')$ is the closed ball of radius $r'>0$ and center $0$ in
$K_v$.  Note that $\Pi(Y_0)=(\pi_v\cO_v)^{\,m}$ is the open unit ball
of $K_v^{\,m}$ and that $\Pi(Y_k)$ is an additive subgroup of
$K_v^{\,m}$. Since the residual field $k_v= \cO_v/\pi_v\cO_v$ lifts as a
subfield of order $q_v$ of $K_v$, for every $\ell\in\bZ_{\geq 0}$, we
have a disjoint union
$$
\ov{B}(0,q_v^{\,-\ell})=
\bigsqcup_{a\in k_v}\big(\,a\;\pi_v^{\,\ell}+\ov{B}(0,q_v^{\,-\ell-1})\,\big)\;.
$$
Hence by induction, the polydisc $\Pi(Y_{k})$ is the disjoint union of
$$
\Delta_{k+1}= \prod_{1\leq j\leq m}q_v^{n_{k+1,j}-n_{k,j}}
$$
translates of the polydisc $\Pi(Y_{k+1})$. Note that
\begin{equation}\label{eq:minorDeltak}
\Delta_{k+1}\geq \prod_{1\leq j\leq m} Y_{k+1}^{\,r_j}Y_{k}^{\,-r_j}q_v^{-1}
=q_v^{\,-m}\,\big(Y_{k+1}Y_{k}^{\,-1}\big)^{|\mb{r}|}\;.
\end{equation}
For every $k\in\bZ_{\geq 0}$, let us fix some elements
$\theta_{1,k+1},\dots, \theta_{\Delta_{k+1},k+1}$ in
$(\pi_v\cO_v)^{\,m}$ (which are not unique in the ultrametric space
$K_v^{\,m}$) such that
$$
\Pi(Y_{k})=\bigsqcup_{i=1}^{\Delta_{k+1}}
\big(\,\theta_{i,k+1}+\Pi(Y_{k+1})\,\big)\;.
$$

By convention, let us define $Z_{0,\delta}=\emptyset$ and $I_0 =
\{\Pi(Y_0)\}$. For every $k\in\bZ_{\geq 1}$, let us define
$$
Z_{k,\delta}= \{\bm{\theta}\in(\pi_v\cO_v)^{\,m}:
|\idist{\mb{y}_k\cdot\bm{\theta}}|<\delta\}
$$
and
$$
I_k=\big\{\theta_{i_1,1}+\dots+\theta_{i_k,k}+\Pi(Y_{k}):
\forall\;j\in\{1,\dots, k\},\;\;1\leq i_j\leq \Delta_{j}\big\}\;.
$$

\blemm\label{lem:Delvol} For every $k\in\bZ_{\geq 1}$, we have
\begin{enumerate}
\item[(1)] for every $I'\in I_{k+1}$, if $I'\cap Z_{k,\delta}\neq
  \emptyset$ then $I'\subset Z_{k,\delta}$,
\item[(2)] for every $I\in I_{k}$,  we have
  $\vol_v^m(I\cap Z_{k,\delta})\leq \delta\; Y_{k}^{-|\mb{r}|}$.
\end{enumerate}
\elemm

\dem (1) If $I'\in I_{k+1}$ and $I'\cap Z_{k,\delta}\neq\emptyset$,
let $\bm{\theta}\in I'\cap Z_{k,\delta}$.  Then for every $\bm{\theta}'
\in I'$, if $x,x'\in R_v$ are such that $|\idist{\mb{y}_k\cdot
  \bm{\theta}}| = |(\mb{y}_k \cdot \bm{\theta})-x|$ and
$|\idist{\mb{y}_k\cdot (\bm{\theta}' -\bm{\theta})}| =|(\mb{y}_k\cdot
(\bm{\theta}'-\bm{\theta})) -x'|$, then by the ultrametric inequality,
since $\bm{\theta}\in Z_{k,\delta}$ and $\bm{\theta}'- \bm{\theta}
\in\Pi(Y_{k+1})$, by the assumption of Proposition \ref{prop:DProp},
and by the definition of $b$, we have
\begin{align*}
  |\idist{\mb{y}_k \cdot \bm{\theta}'}|& \leq
  |\mb{y}_k \cdot (\bm{\theta}+(\bm{\theta}'-\bm{\theta}))- (x+x')|
  \leq \max \left\{|(\mb{y}_k \cdot \bm{\theta})-x|,\;|(\mb{y}_k\cdot
  (\bm{\theta}'-\bm{\theta})) -x'|\right\}\\&=
  \max \left\{|\idist{\mb{y}_k\cdot\bm{\theta}}|,\;
  |\idist{\mb{y}_k\cdot (\bm{\theta}'-\bm{\theta})}|\right\}\\
  &\leq \max\big\{\delta,
  \max_{1\leq j\leq m} Y_{k}^{r_j} \frac{1}{q_v}\;Y_{k+1}^{-r_j} \big\}
  \leq \max\big\{\delta, q_v^{-1-b\min \mb{r}}\big\} = \delta\;.
\end{align*}
This inequality $|\idist{\mb{y}_k \cdot \bm{\theta}'}|\leq \delta$ is
actually strict, since $|\idist{\mb{y}_k\cdot \bm{\theta}}| <\delta$
and by Equation \eqref{eq:defibbdelta}, we have $q_v^{-1-b\min
  \mb{r}}= q_v^{-1}\delta <\delta$. Since $I'$ is contained in
$\Pi(Y_0)=(\pi_v\cO_v)^{\,m}$, we thus have that $\bm{\theta}'\in
Z_{k,\delta}$ and this proves Assertion (1).

\medskip
(2) Let $j_0\in\{1,\dots, m\}$ be such that $Y_{k} =
|y_{k,j_0}|^{1/r_{j_{0}}}$ where $\mb{y}_{k} = (y_{k,1}, \dots,
y_{k,m})$. In particular, $y_{k,j_0}$ is nonzero. For every $z\in
R_v$, let
$$
L_k(z)=\{\bm{\theta}\in K_v^{\,m} : \mb{y}_k\cdot\bm{\theta}=z\}\;,
$$
which is an affine hyperplane of $K_v^{\,m}$ transverse to the
$j_0$-axis, and let
$$
\cN(k,z)=\{\bm{\theta}'\in(\pi_v\cO_v)^m:
\exists\;\mb{u}'\in L_k(z),\;\;
|\theta_{j_0}'-u_{j_0}'|\leq \delta \;Y_k^{-r_{j_0}}
\;\;{\rm and}\;\;\forall\;j\neq j_0,\; \theta_{j}'=u_{j}'\}\;,
$$
which is the intersection with the open unit ball in $K_v^{\,m}$ of
the $(\delta \;Y_k^{-r_{j_0}})$-thickening along the $j_0$-axis of the
affine hyperplane $L_k(z)$.

Fix $I\in I_k$. Since $\vol_v(\ov{B})(0,r'))=q_v^{\lfloor\log_{q_v}r'\rfloor}
\leq r'$ for all $r'>0$, and by Fubini's theorem, we have
\begin{equation}\label{eq:volthickening}
  \vol_v^m(I\cap \cN(k,z))\leq
  \delta \;Y_k^{-r_{j_0}} \;\prod_{j\neq j_0} Y_k^{-r_{j}}
  =\delta \;Y_k^{-|\mb{r}|}\;.
\end{equation}

\medskip\noindent {\bf Claim 1. } Let us prove that the set
$Z_{k,\delta}$ is contained in the union of the sets
$\cN(k,z)$ for $z\in R_v$.

\medskip
\dem Let $\bm{\theta}=(\theta_1,\dots,\theta_m)\in Z_{k,\delta}$ and
let $z\in R_v$ be such that $|\idist{\mb{y}_k \cdot \bm{\theta}}| =
|\,\mb{y}_k \cdot \bm{\theta} - z\,|$.  Let us define $u_j=\theta_j$ if
$j\neq j_0$,
$$
u_{j_0}= \frac{z-\sum_{j\neq j_0}y_{k,j}\theta_{j}}{y_{k,j_0}}
$$
and $\mb{u}=(u_1,\dots, u_m)$, which is the projection of
$\bm{\theta}$ on the affine hyperplane $L_k(z)$ along the
$j_0$-axis. Then, since $\bm{\theta}\in Z_{k,\delta}$, we have
$$
|\theta_{j_0}-u_{j_0}| =
\frac{|\,\mb{y}_k \cdot \bm{\theta} -z\, |}{|y_{k,j_0}|} =
\frac{|\idist{\mb{y}_k \cdot \bm{\theta}}|}{|y_{k,j_0}|} \leq
\delta \;Y_k^{-r_{j_0}}
\;.
$$
Since $Z_{k,\delta}$ is contained in $(\pi_v\cO_v)^m$, this proves
Claim 1.
\cqfd

\medskip\noindent {\bf Claim 2. } Let us prove that
there exists a unique $z\in R_v$ such that $I\cap Z_{k,\delta}$ is
contained in $I\cap \cN(k,z)$.

\medskip
\dem By Claim 1, the set $I\cap Z_{k,\delta}$ is contained in
$\bigcup_{z\in R_v}I\cap \cN(k,z)$. Assume for a contradiction that
there exist two distinct elements $z,z'$ in $R_v$ such that there
exist $\bm{\theta} \in I\cap\cN(k,z)$ and $\bm{\theta}'\in
I\cap\cN(k,z')$. Let $\mb{u}\in L_k(z)$ and $\mb{u'}\in L_k(z')$ be
the projections of $\bm{\theta}$ and $\bm{\theta}'$ along the
$j_0$-axis on $L_k(z)$ and $L_k(z')$ respectively.

Let $j\in \{1,\dots, m\}$. Note that $\bm{\theta} -\bm{\theta}'
\in\Pi(Y_k)$ since $I\in I_k$. If $j\neq j_0$, then
$$
|u_j-u'_j|=|\theta_j-\theta_j'|\leq \frac{1}{q_v}\;Y_k^{-r_j}\;.
$$
Furthermore, by the ultrametric inequality, since $\bm{\theta}$
(respectively $\bm{\theta}'$) is contained in the $(\delta
\;Y_k^{-r_{j_0}})$-thickening along the $j_0$-axis of $L_k(z)$
(respectively $L_k(z')$), and since $\delta\leq \frac{1}{q_v}$, we
have
\begin{align*}
|u_{j_0}-u'_{j_0}|&= |(u_{j_0}-\theta_{j_0})
+(\theta_{j_0}-\theta_{j_0}')+(\theta_{j_0}'-u'_{j_0})|
\\ &\leq \max\{|u_{j_0}-\theta_{j_0}|,\;|\theta_{j_0}-\theta_{j_0}'|,
\;|\theta_{j_0}'-u'_{j_0}|\}
\\&\leq \max\{\delta\;Y_k^{-r_{j_0}},\;\frac{1}{q_v}\;Y_k^{-r_{j_0}}\}=
\frac{1}{q_v}\;Y_k^{-r_{j_0}}\;.
\end{align*}
This implies since $\mb{u}\in L_k(z)$ and $u'\in L_k(z)$ that
\begin{align*}
1 \leq |z-z'| = |\mb{y}_k \cdot \mb{u} -\mb{y}_k \cdot\mb{u}'|
\leq \max_{1\leq j \leq m} |y_{k,j}|\,|u_j-u'_j|\leq
\max_{1\leq j \leq m} Y_k^{r_j}\,\frac{1}{q_v}\;Y_k^{-r_j}=\frac{1}{q_v} \;, 
\end{align*}
which is a contradiction since $q_v>1$. This proves Claim 2.
\cqfd

\medskip
By Equation \eqref{eq:volthickening}, Claim 2 concludes the proof of
Assertion (2) of Lemma \ref{lem:Delvol}.
\cqfd

\medskip
Since every element $I'$ of $I_{k+1}$ is a translate of
$\Pi(Y_{k+1})$, and by Equation \eqref{eq:defnkj}, we have
$$
\vol_v^m (I')= \vol_v^m (\Pi(Y_{k+1})) =\prod_{j=1}^mq_{v}^{-n_{k+1,\,j}-1}
\geq q_v^{-2m}\;Y_{k+1}^{-|\mb{r}|}\;.
$$
For every $I\in I_k$, there are $\Delta_{k+1}$ elements $I'\in
I_{k+1}$ contained in $I$, they are pairwise disjoint and they have
the same volume $\vol_{v}^{m} (\Pi(Y_{k+1}))$. Among them, those who
meet $Z_{k,\delta}$ are actually contained in $I\cap Z_{k,\delta}$ by
Lemma \ref{lem:Delvol} (1), thus their number is at most
$\frac{\vol_{v}^{m}(I\cap Z_{k,\delta})} {\vol_{v}^{m}(\Pi(Y_{k+1}))}$.
Therefore, by Equation \eqref{eq:minorDeltak} and Lemma
\ref{lem:Delvol} (2), we have
\begin{align}
\card\;\{I'\in I_{k+1}: I'\subset I,\;\; I'\cap Z_{k,\delta}=\emptyset\}
&\geq \Delta_{k+1}-
\frac{\vol_{v}^{m}(I\cap Z_{k,\delta})}{\vol_{v}^{m}(\Pi(Y_{k+1}))}
\nonumber\\ & \geq q_v^{-m}(Y_{k+1}Y_{k}^{-1})^{|\mb{r}|}-
\frac{\delta\; Y_{k}^{-|\mb{r}|}}{q_v^{-2m}Y_{k+1}^{-|\mb{r}|}}\nonumber
\\ & = c_1\;(Y_{k+1}Y_{k}^{-1})^{|\mb{r}|}\;, \label{eq:cardstep}
\end{align}
where $c_1=q_v^{-m}-q_v^{\,2m} \delta$ belongs to $]0,1[$ by the
assumption on $\delta$.

\medskip
Now, let us define by induction $J_0=I_0$ and for every $k\in\bZ_{\geq 0}$,
$$
J_{k+1}=
\bigcup_{J\in J_k}\{I\in I_{k+1}:I\subset J,\;\; I\cap Z_{k,\delta}=\emptyset\}.
$$
By Equation \eqref{eq:cardstep} and by induction, we have
\begin{equation}\label{eq:minocardJk}
\card \;J_{k+1}\geq\prod_{j=1}^k \;c_1\;(Y_{j+1}Y_{j}^{-1})^{|\mb{r}|}
=c_1^{\,k}\,(Y_{k+1}Y_{1}^{-1})^{|\mb{r}|}\;.
\end{equation}
By Lemma \ref{lem:Delvol} (1) and by induction, we have
$$
J_{k+1}=\{J\in I_{k+1}:\forall\;j\in\{1,\dots, k\},\;\;
J\cap Z_{j,\delta}=\emptyset\}=\{J\in I_{k+1}:
J\subset \bigcap_{j=1}^k \;^c Z_{j,\delta}\;\}\;,
$$
where $^c$ denotes the complement in $(\pi_v\cO_v)^m$.  Hence
$\big(\bigcup J_k\big)_{k\geq 1}$ is a decreasing sequence of compact
subsets of $(\pi_v\cO_v)^m$, whose intersection is contained in
$\bigcap_{k\geq 1}\; ^cZ_{k,\delta}=\Bad_{(\mb{y}_i)_{i\geq
    1}}^{\delta}$.

\medskip
For every $k\in\bZ_{\geq 0}$, let us define a measure
$$
\mu_k = \big(\vol_v^m (\Pi(Y_k))\;\card\; J_k\big)^{-1}
\sum_{J\in J_k} \vol_v^m|_{J}\;,
$$
which is a probability measure with support $\bigcup J_k$. By the
compactness of $(\pi_v\cO_v)^m$, any weakstar accumulation point $\mu$
of the sequence $(\mu_k)_{k\geq 1}$ is a probability measure with
support in $\Bad_{(\mb{y}_i)_{i\geq 1}}^{\delta}$.

For every closed ball $B$ in $(\pi_v\cO_v)^m$ with radius $r'
\in\;]0,r_0= Y_1^{-\min\mb{r}}]$, let $k\in\bZ_{\geq1}$ be such that
\begin{equation}\label{eq:defkpourrprim}
Y_{k+1}^{-\min\mb{r}} < r' \leq Y_{k}^{-\min\mb{r}}\;.
\end{equation}
Note that $\lceil t\rceil\leq t+1\leq q_v\,t$ if $t\geq 1$, and that $r'\,
q_v^{n_{k+1,j}+1} \geq Y_{k+1}^{-\min\mb{r}}\; Y_{k+1}^{r_j} \;q_v\geq 1$
for every $j=1,\dots, m$, by Equation \eqref{eq:defnkj}.  Then $B$ can
be covered by a subset of $I_{k+1}$ with cardinality at most
$$
\prod_{j=1}^m \;\big\lceil r'\;q_v^{n_{k+1,j}+1} \big\rceil
\leq (r')^m \;q_v^{\,3m}\; Y_{k+1}^{\;|\mb{r}|}\;.
$$
Let $C= \frac{-\log_{q_v} c_1}{\min\mb{r}}>0$, which depends (besides
on $m$, $q_v$ and $\mb{r}$) only on $\delta$. Defining $C_0 =
q_v^{\,3m} \;Y_1^{\;|\mb{r}|}$, by Equations \eqref{eq:minocardJk} and
\eqref{eq:defkpourrprim}, we thus have
\begin{align*}
  \mu_{k+1}(B)&\leq
  q_v^{\,3m} \;(r')^m \;Y_{k+1}^{\;|\mb{r}|} \;(\card\; J_{k+1})^{-1}
  \leq q_v^{\,3m}\; (r')^m\;c_1^{-k}\;Y_1^{\;|\mb{r}|} \\ &
  \leq C_0 \;(r')^{m-C\frac{k}{\log_{q_v}Y_k}}\;.
\end{align*}
Therefore, since the ball $B$ is closed and open and since $r'\leq
r_0\leq 1$, we have
$$
\mu(B) \leq \limsup_{k\to\infty} C_0 \;(r')^{m-C\frac{k}{\log_{q_v}Y_k}}=
C_0 \;(r')^{m-C\limsup_{k\to\infty}\frac{k}{\log_{q_v}Y_k}},
$$
which concludes the proof of Proposition \ref{prop:DProp}.
\cqfd

\subsection{Proof that Assertion \eqref{S2_introA1} implies
  Assertion \eqref{S1_introA1} in Theorem \ref{theo:introA1}}
\label{subsec:S2A1impliesS1A1}

Suppose that $A$ is $(\mb{r},\mb{s})$-singular on average.  Then by
Corollary \ref{coro:AsingifftAsing}, the matrix $\,^{t}\!A$ is also
$(\mb{s},\mb{r})$-singular on average. By Remark
\ref{rem:noncomplirrat} (2), in order to prove that there exists
$\epsilon > 0$ such that $\mb{Bad}_{A}(\epsilon)$ has full Hausdorff
dimension, we may assume that the matrix $\,^{t}\!A$ is completely
irrational.

By Lemma \ref{lem:bestapprox}, let $(\mb{y}_k)_{k\in\bZ_{\geq 1}}$ be
a best approximation sequence in $K_v^{\,m}$ for the matrix
$\,^{t}\!A$ with weights $(\mb{s},\mb{r})$, and let $c>0$ be such that
$M_{i}Y_{i+1}\leq q_v^c$ for every $i\in\bZ_{\geq 1}$.  Fix some
$\delta\in\, \big]0, \;\frac{1}{q_v^{\,3m}}\big[$ and let
$b=b(\delta)>0$ and $C=C(\delta)>0$ as in Proposition
\ref{prop:DProp}. By Proposition \ref{prop:modBZseq}, for every $a >
b$, we have a subsequence $(\mb{y}_{\varphi(k)})_{k\geq 1}$ such that
the properties \eqref{eq:modBZseq1} and \eqref{eq:modBZseq2} are
satisfied. Proposition \ref{prop:SProp}, whose assumption is satisfied
by the second inequality in Equation \eqref{eq:modBZseq1} and where
$\epsilon=\delta^{\frac{1}{\min \mb{r}} + \frac{1}{\min
    \mb{s}}}\;q_v^{-b-c}$, gives that $\Bad_A(\epsilon)$ contains
$\Bad_{(\mb{y}_{\varphi(i)})_{i\geq 1}}^{\delta}$. Therefore, using
Proposition \ref{prop:DProp} applied to the sequence
$(\mb{y}_{\varphi(i)})_{i\geq 1}$, whose assumption is satisfied by
the first inequality in Equation \eqref{eq:modBZseq1}, and using
Equation \eqref{eq:modBZseq2} for the last inequality, we have
\begin{align*}
  \dimH \Bad_A(\epsilon)&\geq
  \dimH \Bad_{(\mb{y}_{\varphi(i)})_{i\geq 1}}^{\delta}
  \geq m-C\limsup_{k\to\infty}\frac{k}{\log_{q_v}Y_{\varphi(k)}}
  \geq m-\frac{C}{a}\;.
\end{align*}
Letting $a$ tend to $+\infty$, this concludes the proof that
Assertion \eqref{S2_introA1} implies Assertion \eqref{S1_introA1} in
Theorem \ref{theo:introA1}.  \cqfd

\section{Background material for the upper bound}
\label{sec:prelimII}

\subsection{Homogeneous dynamics}
\label{subsec:homogdyn}

Let $K_v,\OOO_v,\pi_v,R_v,q_v$ be as in Subsection
\ref{subsec:funcfield}. Let $m,n\in\NN-\{0\}$ and $d=m+n$. We fix some
weights $\mb{r}=(r_1,\dots,r_m) $ and $\mb{s}=(s_1,\dots,s_n)$ as in
the introduction.  In this subsection, we introduce the space of
unimodular grids $\cY$ in $K_v^{\,d}$ and the diagonal flow
$(\ta^\ell)_{\ell\in\bZ}$ acting on this space. Let
\[
G_0=\SL_d(K_v) \quad\text{and}\quad
G=\ASL_d(K_v)=\SL_d(K_v)\ltimes K_v^{\,d},
\]
and let
\[
\Ga_0=\SL_d(R_v) \quad\text{and}\quad
\Ga=\ASL_d(R_v)=\SL_d(R_v)\ltimes R_v^{\,d}.
\] 
The product in $G$ is given by
\begin{equation}\label{eq:prodG}
  (g,u)\cdot (g',u') = (gg',u+gu')
\end{equation}
for all $g,g'\in G_0$ and $u,u'\in K_v^{\,d}$. We also view $G$ as a
subgroup of $\SL_{d+1}(K_v)$ by
\[
G=\left\{\begin{pmatrix}
g & u \\
0 & 1 \\
\end{pmatrix} : g\in \SL_d(K_v),\; u\in K_v^{\,d}  \right\}\;.
\]

We shall identify $G_0$ with the corresponding subgroup of $G$.  We
consider the one-parameter diagonal subgroup $(\ta^\ell)_{\ell\in\bZ}$
of $G_0$, where $\ta= \diag(\ta_-,\ta_+)$ and
\[
\ta_-= \diag(\pi_v^{-r_1},\dots,\pi_v^{-r_m})\in\GL_m(K_v)
\quad {\rm and} \quad
\ta_+ = \diag(\pi_v^{s_1},\dots,\pi_v^{s_n})\in\GL_n(K_v)\;.
\]
Note that for all $\bm{\theta}\in K_v^{\,m}$, $\bm{\xi}\in K_v^{\,n}$ and
$\ell\in\bZ$, we have
\begin{equation}\label{eq:dilatweight}
  \|\,\ta_-^{\,\ell}\,\bm{\theta}\,\|_{\mb{r}}=
  q_v^{\ell}\;\|\,\bm{\theta}\,\|_{\mb{r}}
  \quad{\rm and} \quad  \
  \|\,\ta_+^{\,\ell}\,\bm{\xi}\,\|_{\mb{s}}=
  q_v^{-\ell}\;\|\,\bm{\xi}\,\|_{\mb{s}}\;.
\end{equation}

We denote by $G^+$ the unstable horospherical subgroup for $\ta$
in $G$ and by $U$ the unipotent radical of $G$,
that is,
\[
G^+=\{g\in G : \lim_{\ell\to -\infty} \ta^\ell g \,\ta^{-\ell}= I_{d+1}\}
\;\;\;{\rm and}\;\;\;U=\Big\{\begin{pmatrix} I_d & u \\
0 & 1 \end{pmatrix} : u\in K_v^{\,d} \Big\}\;.
\]
\noindent Let $U^{+}=G^+ \cap U=
\!\Big\{\!\begin{pmatrix} I_m & 0 & w \\ 0& I_n & 0 \\ 0 & 0 &
1 \end{pmatrix} : w\in K_v^{\,m} \Big\}$, which is a closed subgroup
in $G^+$ normalized by $\ta$.

Let us define
\[
\cX = G_0 /\Ga_0\quad{\rm and}\quad \cY=G/\Ga\;.
\]
Even though we have $\covol(R_v^{\,d})=q^{(g-1)d}$ by Equation
\eqref{eq:covolRv}, we say that an $R_v$-lattice $\Lambda$ in
$K_v^{\,d}$ is \textit{unimodular} if $\covol(\Lambda)
=\covol(R_v^{\,d})$.  A translate in the affine space $K_v^{\,d}$ of
an unimodular lattice is called an \textit{unimodular grid}.  We
identify the homogeneous space $\cX= \SL_d(K_v)/\SL_d(R_v)$ with the
space of unimodular lattices in $K_v^{\,d}$ by the equivariant
homeomorphism
\[
x=g\,\Ga_0 \mapsto\Lambda_x=g\,R_v^{\,d}\;,
\]
and the homogeneous space $\cY=\ASL_d(K_v)/\ASL_d(R_v)$ with the space
of unimodular grids by the equivariant homeomorphism
\begin{equation}\label{eq:defiwtLambdasuby}
y=\begin{pmatrix} g & u\\ 0 & 1\end{pmatrix}\Ga\mapsto \wt\Lambda_y=
g\,R_v^{\,d} + u\;.
\end{equation}
We denote by $\pi:\cY\to \cX$ the natural projection map (forgetting
the translation factor), which is a proper map. Note that the fibers of
$\pi$ are exactly the orbits of $U$ in $\cY$, and in particular each
orbit under $U^+$ in $\cY$ is contained in some fiber of $\pi$ (see
Lemma \ref{lem:isomphiA} for a precise understanding of the
$U^+$-orbits).

\medskip
For every $N\in\NN-\{0\}$, we denote by $d_{\SL_N(K_v)}$ the
right-invariant distance on $\SL_N(K_v)$ defined by
\[
\forall\;g,h\in \SL_N(K_v),\;\;\;d_{\SL_N(K_v)}(g,h)=
\max\{\,\ln(1+\interleave \,gh^{-1}-\id\interleave\,),\;
\ln(1+\interleave \,hg^{-1}-\id\interleave\,)\,\}\;,
\]
where $\interleave\;\;\interleave$ is the operator norm on $\cM_N(K_v)$
defined by the sup norm $\|\;\|$ on $K_v^{\,N}$.  We endow every
closed subgroup $H$ of $G$ with the right-invariant distance $d_H$ on $H$,
which is the restriction to $H$ of the distance $d_{\SL_{d+1}(K_v)}$.
For instance, identifying the additive group $K_v^{\,m}$ with $U^+$ by
the map $w\mapsto \wh w=\begin{pmatrix}I_m & 0 & w\\ 0 & I_n & 0\\ 0 &
0 & 1\end{pmatrix}$, we have 
\begin{equation}\label{eq:dUplus1}
\forall\;w,w'\in K_v^{\,m},\;\;\;d_{U^+}(\wh w,\wh {w'})=\ln(1+\|w-w'\|)\;.
\end{equation}
We also consider the distance $d_{U^+, \|\;\|}$ on $U^+$ such that the map
$w\mapsto \wh w$ is an isometry between $K_v^m$ endowed with the
distance induced by the norm $\|\cdot\|$ and $(U^+,d_{U^+,\|\;\|})$~:
\begin{equation}\label{eq:dUplus}
\forall\;w,w'\in K_v^{\,m},\;\;\;d_{U^+,\|\;\|}(\wh w,\wh {w'})=\|w-w'\|\;.
\end{equation}
Observe that the distances $d_{U^+,\|\;\|}$ and $d_{U^+}$ on $U^+$
are locally equivalent. More precisely, for all $w,w'\in K_v^{\,m}$,
\begin{equation}\label{eq:locLip}
  {\rm if}\quad d_{U^+}(\wh w,\wh {w'})\leq\ln 2 \quad{\rm then}\quad
  \frac{1}{2}\|w-w'\| \leq d_{U^+}(\wh w,\wh {w'})\leq \|w-w'\|.
\end{equation} 

We endow $\cY=G/\Ga$ with the quotient distance $d_\cY$ of the
distance $d_G$ on $G$, defined by
\begin{equation}\label{eq:defidcY}
\forall\,y,y'\in\cY,\;\;\; d_\cY(y,y')=
\min_{\ga\in\Ga}\;d_G(\,\wt{y}\,\ga ,\,\wt{y}\,')
\end{equation} 
for any representatives $\wt{y}$ and $\wt{y}\,'$ of the classes $y$
and $y'$ in $G/\Ga$ respectively. This is a well defined distance
since the canonical projection $G\ra \cY$ is a covering map and the
distance $d_G$ on $G$ is right-invariant. Given any closed subgroup
$H$ of $G$, we denote by $B_H(x,r)$ (respectively $B_\cY(x,r)$) the
open ball of center $x$ and radius $r>0$ for the distance $d_H$
(respectively $d_\cY$), and by $B_r^H$ the open ball
$B_H(\id,r)$. Note that for all $y\in\cY$ and $r>0$, we have (for the
left action of subsets of $G$ on $\cY$)
\[
B_\cY(y,r)=B_r^Gy\;.
\]
We denote by $B_r^{U^+,\|\;\|}$ the open ball of center $\id$ and
radius $r>0$ for the distance $d_{U^+,\|\;\|}$ on $U^+$. Remark that
$B_r^{U^+,\|\;\|}=B_{\ln(1+r)}^{U^+}$ for every $r>0$, by Equations
\eqref{eq:dUplus1} and \eqref{eq:dUplus}.

\blemm\label{lem:contractboulUplus} For all $\epsilon>0$ and
$k\in\bZ_{\geq 0}$, we have 
$$
\ta^{-k}B_\epsilon^{U^+,\|\;\|}\ta^k\subset
B_{\epsilon \,q_v^{-k\min\mb{r}}}^{U^+,\|\;\|}
\quad \text{and}\quad \ta^{-k}B_\epsilon^{U^+}\ta^k\subset
B_{\ln(1+(e^\epsilon-1) \,q_v^{-k\min\mb{r}})}^{U^+}\;.
$$
Similary, we have
$$
\ta^{k}B_\epsilon^{U^+,\|\;\|}\ta^{-k}\subset
B_{\epsilon \,q_v^{k\max\mb{r}}}^{U^+,\|\;\|}
\quad \text{and}\quad \ta^{k}B_\epsilon^{U^+}\ta^{-k}\subset
B_{\ln(1+(e^\epsilon-1)\,q_v^{k\max\mb{r}})}^{U^+}\;.
$$
\elemm

\dem
The proof of the second assertion being similar, we only prove
the first one. For every $w=(w_1,\dots,w_m)\in K_v^{\,m}$, we have
$\ta^{-k}\wh w\;\ta^k=\widehat{\ta_-^{-k}w}$ and
\[
\|\,\ta_-^{-k}w\,\| =
\max_{1\leq i\leq m}|\,\pi_v^{r_ik}w_i\,|\leq q_v^{-k\min\mb{r}}
\;\|\,w\,\|\;.
\]
The part on the left of the first assertion hence follows from
Equation \eqref{eq:dUplus}. The part on the right of the first
assertion follows from the remark just before the statement of Lemma
\ref{lem:contractboulUplus}.  \cqfd

\medskip
Given a point $x$ in $\cY$ (and similarly for $x$ in $\cX$), we define
the {\it injectivity radius} of $\cY$ at $x$ to be
\[
\inj(x)=\sup\big\{r>0: \forall\;\ga\in\Ga-\{\id\},\;\;
  B_G(\wt x,r)\cap B_G(\wt x\,\ga,r)=\emptyset\big\}\;,
\]
which does not depend on the choice of $\wt x\in G$ such that $x=\wt
x\,\Ga$, and is positive and finite since the canonical projection $G\ra
\cY$ is a nontrivial covering map. For every $r>0$, we denote the {\it
  $r$-thick part} of $\cY$ by
\[
\cY(r)=\{x\in\cY:\inj(x)\geq r\}\;.
\]
It follows from the finiteness of a Haar measure of the homogeneous
space $\cY$ that $\cY(r)$ is a compact subset of $\cY$ for every
$r>0$, and that the Haar measure of the {\it $r$-thin part}
$\cY-\cY(r)$ tends to $0$ as $r$ goes to $0$. For every compact subset
$K$ of $\cY$, there exists $r>0$ such that $K\subset \cY(r)$.

\subsection{Dani correspondence}
\label{subsec:Dani}

In this subsection, we give an interpretation of the property for a
matrix $A \in \cM_{m,n}(K_v)$ to be $(\mb{r},\mb{s})$-singular on
average in terms of dynamical properties of the action of the
one-parameter diagonal subgroup $(\ta^\ell)_{\ell\in\bZ}$ on the space
of unimodular lattices, as originally developped by Dani (see for
instance \cite[\S 4]{Kleinbock99}). For every $A \in \cM_{m,n}(K_v)$,
let $u_A = \begin{pmatrix} I_m & A \\ 0& I_n \end{pmatrix}\in G_0$.

\bprop\label{prop:critdynescapeonaver} A matrix $A \in \cM_{m,n}(K_v)$
is $(\mb{r},\mb{s})$-singular on average if and only
if the forward orbit $\{\ta^\ell u_A R_v ^{\,d} : \ell\in\bZ_{\geq 0}
\}$ in $\cX$ of the lattice $u_A R_v ^{\,d}$ under $\ta$ {\rm
  diverges on average} in $\cX$, that is, if and only if for any
compact subset $Q$ of $\cX$, we have
$$
\lim_{N \to \infty} \frac{1}{N} \;\card \{ \ell \in \{ 1 ,\cdots, N \}:
\ta^\ell \,u_A\, \Ga_0 \in Q \} =0\;.
$$
\eprop

\dem Let $Q$ be a compact subset of $\cX$. By Mahler's compactness
criterion (see for instance \cite[Theo.~1.1]{KleShiTom17}), there
exists $\varepsilon\in\;]0,1[$ such that $Q$ is contained in 
\[
\cX_{>\epsilon}= \{g\,R_v^{\,d}\in\cX : \forall \;(\bm{\theta},\bm{\xi})
\in g\,R_v^{\,d} - \{ 0\}\subset K_v^{\,m}\times K_v^{\,n},\;\;
\max \{\|\,\bm{\theta}\,\|_{\mb{r}},
\|\,\bm{\xi}\,\|_{\mb{s}}\} > \varepsilon\}\;,
\]
which is the subset of $\cX$ consisting of the unimodular lattices
with systole (for an appropriate quasinorm) larger than $\epsilon$.
Observe that by Equation \eqref{eq:dilatweight}, for all sufficiently
large $\ell \in \bZ_{\geq 1},$ there exists an element $\mb{y} \in
R_v^{\,n}-\{ 0 \}$ such that $\idist{A\,\mb{y}}_{\mb{r}} \leq \varepsilon
q_v^{-\ell}$ and $\| \,\mb{y}\, \|_{\mb{s}} \leq \varepsilon q_v^\ell$ if and
only if we have $\ta^\ell u_A R_v^{\,d}=
\begin{pmatrix} \ta_-^{\,\ell} &0 \\0 &\ta_+^{\,\ell} \end{pmatrix}
\begin{pmatrix} I_m & A\\0 & I_n \end{pmatrix}R_v^{\,d} \in
\cX- \cX_{> \varepsilon}$.

With $\ell_\epsilon =\lfloor-\log_{q_v}\varepsilon \rfloor$, it follows that
\begin{align*}
0\leq \;&
\card \{ \,\ell\in\{1,\cdots,N\}: \ta^\ell\, u_A\, R_v^{\,d} \in Q \,\} \\
\leq \;&
\card \{ \,\ell\in\{1,\cdots,N\}:
\ta^\ell\, u_A\, R_v^{\,d} \in \cX_{>\varepsilon} \,\}\\
= \;& 
\card \{\ell\in\{1,\cdots,N\}:
\nexists\; \mb{y} \in R_v^{\,n}-\{0\},\;\;
\idist{A\,\mb{y}}_{\mb{r}} \leq \varepsilon q_v^{-\ell} , \;\;
\|\,\mb{y}\,\|_{\mb{s}} \leq \varepsilon q_v^\ell \}\\
\leq \;& 
\card \{\ell\in\{1,\cdots,N\}:
\nexists\; \mb{y} \in R_v^{\,n}-\{0\},\;\;
\idist{A\,\mb{y}}_{\mb{r}} \leq \frac{\varepsilon^2}{q_v}
q_v^{-(\ell-\ell_\epsilon)} , \;\;
\|\,\mb{y}\,\|_{\mb{s}} \leq  q_v^{\ell-\ell_\epsilon} \}\\
\leq \;&  
\ell_\epsilon+\card \{\ell\in\{1,\cdots,N-\ell_\epsilon\}:
\nexists\; \mb{y} \in R_v^{\,n}-\{0\},\;\;
\idist{A\,\mb{y}}_{\mb{r}} \leq \frac{\varepsilon^2}{q_v}
q_v^{-\ell} , \;\;
\|\,\mb{y}\,\|_{\mb{s}} \leq  q_v^{\ell} \}\;.
\end{align*}
After dividing by $N$ (or equivalently by $N-\ell_\epsilon$) this last
expression, its limit as $N$ tends to $0$ exists and is equal to $0$
if $A$ is $(\mb{r},\mb{s})$-singular on average (see Equation
\eqref{eq:defisingonaver}). Hence we have $ \lim_{N \to \infty}
\frac{1}{N} \;\card \{ \ell \in \{ 1 ,\cdots, N \}: \ta^\ell \,u_A\,
\Ga_0 \in Q \} =0 $ by the above string of (in)equalities.

The converse implication follows similarly by taking for the compact set
$Q$ the subset $\cX_{> \varepsilon}$.
\cqfd

\medskip
We denote by $\|\;\;\|_{\mb{s},\mb{r}}$ the quasi-norm on
$K_v^{\,d}=K_v^{\,m}\times K_v^{\,n}$ defined by
$$
\|\,(\bm{\theta},\bm{\xi})\,\|_{\mb{r},\mb{s}}=
\max\big\{\|\,\bm{\theta}\,\|_{\mb{r}}^{\;\frac{d}{m}},\;
\|\,\bm{\xi}\,\|_{\mb{s}}^{\;\frac{d}{n}}\big\}\;.
$$
Let $\varepsilon>0$. We define the {\it $\epsilon$-compact part} of
$\cY$ to be
\begin{equation}\label{eq:defiLsubeps}
\cL_\varepsilon = \{y\in\cY: \forall\;u \in \wt\Lambda_y,\;\;
\| \,u\, \|_{\mb{r}, \mb{s}} \geq \varepsilon \}\;.
\end{equation}
By Mahler's compactness criterion (see for instance
\cite[Theo.~1.1]{KleShiTom17}) and since the natural projection
$\pi:\cY\ra\cX$ is proper, the subset $\cL_\varepsilon$ is compact.

For every $\bm{\theta}\in K_v^{\,m}$, we denote by $y_{A,\bm{\theta}}$
the unimodular grid $u_A R_v^{\,d} - \begin{psmallmatrix}\bm{\theta}
\\ 0\end{psmallmatrix}$.

\blemm\label{lem:isomphiA} For every $A\in\cM_{m,n}(K_v)$, the map
$K_v^{\,m}\to \cY$ defined by $\bm{\theta}\mapsto y_{A,\bm{\theta}}$
induces a bilipschitz homeomorphism $\phi_A$ from $\bT^m=
K_v^{\,m}/R_v^{\,m}$ endowed with the quotient distance $d_{\bT^m}$ of
the distance on $K_v^{\,m}$ defined by the standard norm $\|\;\|$, and
the $U^+$-orbit $U^+y_{A,0}$ endowed with the restriction of the
distance $d_\cY$ of $\cY$. Furthermore, the map $\phi_A$ is an
isometry between $(\bT^m,d_{\bT^m})$ and the $U^+$-orbit $U^+y_{A,0}$
endowed with the quotient distance $d_{U^+y_{A,0},\|\cdot\|}$ of the
distance $d_{U^+,\|\cdot\|}$ of $U^+$ defined in Equation
\eqref{eq:dUplus}.
\elemm

\dem The map $K_v^{\,m}\to \cY$ defined by $\bm{\theta}\mapsto
y_{A,\bm{\theta}}$ is clearly invariant under translations by
$R_v^{\,m}$, and induces a bijection
\begin{equation}\label{eq:defiphiA}
  \phi_A:\bm{\theta}\!\!\mod R_v^{\,m}\mapsto y_{A,\bm{\theta}}
\end{equation}
from $\bT^m=K_v^{\,m}/R_v^{\,m}$ to the orbit $U^+y_{A,0}$. This orbit
is contained in the fiber $\pi^{-1}(x_A)$ of $x_A=u_A\,R_v^{\,m}$ for
the natural projection $\pi:\cY\ra \cX$, as already seen.

For all $A \in \cM_{m,n}(K_v)$ and $\bm{\theta}\in K_v^{\,m}$, let
$u_{A,\bm{\theta}} = \begin{pmatrix} I_m & A & \bm{\theta} \\ 0& I_n &
  0\\0 & 0 & 1 \end{pmatrix}\in G$, so that we have $y_{A,\bm{\theta}}
=u_{A,\,-\bm{\theta}}\Ga$. For all $\bm{\theta}, \bm{\theta}' \in
\bT^m$, denoting lifts of them to $K_v^{\,m}$ by $\wt{\bm{\theta}},
\wt{\bm{\theta}}'$ respectively, identifying $K_v^{\,d}$ with
$K_v^{\,m}\times K_v^{\,n}$, and using Equations \eqref{eq:defidcY}
and \eqref{eq:prodG}, we have
\begin{align*}
  &d_\cY(\phi_A(\bm{\theta}), \phi_A(\bm{\theta}')) =
  \inf_{\ga\in\Ga} d_G\big(u_{A,\,-\wt{\bm{\theta}}}\;\ga,
  u_{A,\,-\wt{\bm{\theta}}'} \big)\\
  &=\;\inf_{\substack{g\in\Ga_0 \\ x\in R_v^{\,m},\; y\in R_v^{\,n}}}
  d_G\big((u_Ag, (x+Ay -\wt{\bm{\theta}},y)),
  (u_A, (-\wt{\bm{\theta}}', 0))\big)\\
  &=\;\inf_{x\in R_v^{\,m}}
  d_{U^+}\big((\id, (x -\wt{\bm{\theta}},0)),
  (\id, (-\wt{\bm{\theta}}', 0))\big)\\
  &=\; \inf_{x\in R_v^{\,m}}
  \ln(1+\|\,\wt{\bm{\theta}}-\wt{\bm{\theta}}' -x\,\|)=
  \ln(1+d_{\bT^m}(\bm{\theta},\bm{\theta}))\;.
\end{align*}
The last claim of Lemma \ref{lem:isomphiA} hence follows from Equation
\eqref{eq:dUplus1}. By the properties of the logarithm map and by the
compactness of $\bT^m$, which implies that the diameter $\delta$ of
$(\bT^m,d_{\bT^m})$ is finite, we have
\[
\frac{\ln(1+\delta)}{\delta} \;d_{\bT^m} (\bm{\theta}, \bm{\theta}')
\leq d_\cY(\phi_A(\bm{\theta}), \phi_A(\bm{\theta}'))
\leq d_{\bT^m} (\bm{\theta}, \bm{\theta}')
\]
This proves the first claim of Lemma \ref{lem:isomphiA}.
\cqfd

\bprop\label{prop:bad} Let $\varepsilon>0$. For every $(A,\bm{\theta})
\in \M_{m,n}(K_v)\times K_v^{\,m}$ such that $\bm{\theta}\in\Bad_A
  (\varepsilon)$, one of the following statements holds.
\begin{enumerate}
\item
There exists $\mb{y} \in R_v^{\,n}$ such that $\idist{A\,\mb{y}-
  \bm{\theta}}_{\bbr}=0$. Note that given $A$, there are only
countably many $\bm{\theta}$ satisfying this statement.
\item
The forward $\ta$-orbit of the point $y_{A,\bm{\theta}}$ is eventually
in the $\epsilon$-compact part $\cL_\varepsilon$, that is, there
exists $T \geq 0$ such that for every $\ell \geq T$, we have
$\ta^\ell\, y_{A,\bm{\theta}} \in \cL_\varepsilon$.
\end{enumerate}
\eprop

\dem Assume for a contradiction that both statements do not hold.
Then there exist infinitely many $\ell \in \bZ_{\geq 1} $ such that
$\ta^\ell y_{A,\bm{\theta}} \notin \cL_\varepsilon$, hence such that
there exists $\mb{y}_\ell \in R_v^{\,n}$ with $\idist{A\,\mb{y}_\ell-
  \bm{\theta}}_\bbr < q_v^{-\ell} \varepsilon^{\frac{m}{d}}$ and
$\|\,\mb{y}_\ell\,\|_\mb{s} < q_v^\ell \varepsilon^{\frac{n}{d}}$.
Since the statement (1) does not hold, the inequality
\[
\|\,\mb{y}\,\|_s \idist{A\,\mb{y}- \bm{\theta}}_\bbr  < \varepsilon
\]
has infinitely many solutions $\mb{y} \in R_v^{\,n}$, which
contradicts the assumption $\bm{\theta}\in\Bad_A (\varepsilon)$.
\cqfd

\subsection{Entropy, partition construction, and effective variational
  principle}
\label{subsec:effvarprin}

In this subsection, after recalling the basic definitions and
properties about entropy (using \cite{EinLinWar22} as a general
reference, and in particular its Chapter 2), we give the preliminary
constructions of $\sigma$-algebras and results on entropy that will be
needed in Section \ref{sec:upperbound}. In particular, we give an
effective and positive characteristic version of the variational
principle for conditional entropy of \cite[\S 7.55]{EinLin10}. See
\cite{KimKimLim21} for the real case.

\medskip
Let $(X, \B, \mu)$ be a standard Borel probability space. For every
set $E$ of subsets of $X$, we denote by $\sigma(E)$ the
$\sigma$-algebra of subsets of $X$ generated by $E$. Let $\P$ be a
(finite or) countable $\B$-measurable partition of $X$. Let $\A$, $\C$
and $\C'$ be sub-$\sigma$-algebras of $\B$. Suppose that $\C$ and
$\C'$ are countably generated.

For every $x\in X$, we denote by $[x]_\P$ the {\it atom} of $x$ for
$\P$, which is the element of the partition $\P$ containing $x$. We
denote by $[x]_\C$ the {\it atom} of $x$ for $\C$, which is the
intersection of all elements of $\C$ containing $x$. Note that
$[x]_{\sigma(\P)}= [x]_\P$. We denote by $(\mu_x^\A)_{x\in X}$ an
$\A$-measurable family of (Borel probability) {\it conditional
  measures} of $\mu$ with respect to $\A$ on $X$, given for instance
by \cite[Theo. 5.9]{EinLin10}.

Using the standard convention $0 \log_{q_v} 0 =0$ and using
$\log_{q_v}$ instead of $\log$ for computational purposes in the field
$K_v$, the {\it entropy} of the partition $\P$ with respect to $\mu$
is defined by
\[
H_\mu(\P) = - \sum_{P\in \P} \mu(P) \log_{q_v} \mu(P) \;\in [0, \infty]\;.
\]
Recall the (logarithmic) cardinality majoration
\begin{equation}\label{eq:majocard}
  H_\mu(\P)\leq \log_{q_v}(\card\P)\;.
\end{equation}
The {\it information function} of $\C$ given $\A$ with respect to
$\mu$ is the measurable map $I_\mu (\C | \A):X\ra [0, \infty]$ defined by
\[
\forall\;x\in X,\;\;\;\;
I_\mu(\C | \A)(x) = - \log_{q_v} \mu_x^{\A}([x]_\C)\;.
\]
The {\it conditional entropy} of $\C$ given $\A$ with respect to $\mu$
is defined by
\begin{equation}\label{eq:deficondentrop}
H_\mu(\C| \A) = \int_X I_\mu (\C | \A) \;d\mu\;.
\end{equation}
Recall the additivity property $H_\mu(\C\vee\C' \,|\,\A )=H_\mu(\C
\,|\,\C'\vee\A )+H_\mu(\C' \,|\,\A)$ (see for instance
\cite[Prop.~2.13]{EinLinWar22}) so that if $\A\subset\C'\subset\C$,
we have
\begin{equation}\label{eq:chainrule}
  H_\mu(\C\,|\,\A)= H_\mu(\C\,|\,\C')+H_\mu(\C'\,|\,\A)\;.
\end{equation}

Let $T:(X, \B, \mu)\ra(X, \B, \mu)$ be a measure-preserving
transformation. Assume that the $\sigma$-algebra $\A$ is strictly
$T$-invariant, i.e., that $T^{-1}\A = \A$. If the partition $\cP$ has
finite entropy with respect to $\mu$, let
\[
h_{\mu}(T,\P|\A)=\lim_{n\to\infty}\frac{1}{n}H_{\mu}
\Big(\bigvee_{i=0}^{n-1}T^{-i}\P|\A\Big)
=\inf_{n\geq 1}\frac{1}{n}H_{\mu}\Big(\bigvee_{i=0}^{n-1}T^{-i}\P|\A\Big)\;.
\] 
The {\it conditional (dynamical) entropy} of $T$ given $\A$ is 
\[
h_{\mu}(T|\A)=\sup_{\P}\;h_{\mu}(T,\P|\A)\;,
\]
where the upper bound is taken on all countable $\B$-measurable partitions
$\cP$ of $X$ with finite entropy with respect to $\mu$.

Let $\E=\{B\in\B:\mu(T^{-1}B\Delta B)=0\}$ be the sub-$\sigma$-algebra
of $T$-invariant elements of $\B$ modulo $\mu$. We denote by
$(\mu_{x}^{\E})_{x\in X}$ the associated family of conditional
measures. The following result is proven for instance in
\cite[Theo.~2.34]{EinLinWar22}.

\bprop[{\bf Entropy and ergodic decomposition}] \label{prop:ergDec} If
$\A$ is strictly $T$-invariant, then we have
\[
h_{\mu}(T | \A)=\int_{X}h_{\mu_{x}^{\E}}(T |\A)\;d\mu(x)\;.\;\;\;\Box
\]
\eprop

\medskip
We now work in the standard Borel space $\cY$ of unimodular grids,
endowed with the distance $d_\cY$ (see Subsection
\ref{subsec:homogdyn}). Let $\delta>0$. For every subset $B$ of $\cY$,
we define the {\it $\delta$-boundary} $\partial_\delta B$ of $B$ by
\[
\partial_\delta B=\big\{\,y\in\cY:
\inf_{y'\,\in\, B} \;d_\cY(y,y')+\inf_{y''\,\in\, \cY-B}\;
d_\cY(y,y'')<\delta\,\big\}
\]
if $B$ and $\cY-B$ are nonempty, and $\partial_\delta B=\emptyset$
otherwise. Note that for all subsets $B$ and $B'$ of $\cY$, we have
\begin{equation}\label{eq:proprideltneigh}
  \partial_\delta (B\cup B')\subset \partial_\delta B\cup \partial_\delta B'
  \;\;\;{\rm and}\;\;\;\partial_\delta (B- B'\cap B)
  \subset \partial_\delta B\cup \partial_\delta B'\;.    
\end{equation}
We also have $\partial_\delta B\subset \partial_{\delta'} B$ if
$\delta\leq\delta'$.  Given any set $\cP$ of subsets of $\cY$, we
define the {\it $\delta$-boundary} $\partial_\delta \cP$ of $\cP$ by
\[
\partial_{\delta}\cP=\bigcup_{B\in\cP}\;\partial_\delta B\;.
\]

\blemm\label{partconst} For every $r>0$, there exist $\delta_r\in
\;]0,r]$ and a finite measurable partition $\cP= \{P_1, \dots,P_N,
P_\infty\}$ by closed and open subsets of $\cY$ such that
\begin{enumerate}
\item\label{partconst:prop1}
  the subset $P_\infty$ is contained in the $r$-thin part $\cY-\cY(r)$,
\item\label{partconst:prop2}
  for every $i\in\{1,\dots, N\}$, there exists $y_i\in \cY(r)$
  such that $B^G_{\frac{r}{2}}y_i\subset P_i\subset B^G_r y_i$,
\item\label{partconst:prop3}
  the set $\partial_{\delta_r}\P$ is empty.
\end{enumerate}
\elemm

\dem Choose a finite maximal $r$-separated subset $\{y_1,\dots,y_N\}$
of $\cY(r)$ for the distance $d_\cY$, which exists by the compactness of
$\cY(r)$. By induction on $i=1,\dots,N$, we define a Borel subset
$P_i$ of $\cY$ by
\[
P_i=B^G_r y_i-\Big(\bigcup_{j=1}^{i-1}P_j\cup
\bigcup_{j=i+1}^N B^G_{\frac{r}{2}}y_j\Big)\;.
\]
Define $P_\infty =\cY-\bigcup_{j=1}^{N}P_j$, which is also a Borel
subset of $\cY$.

By construction, we have $P_i\subset B^G_r y_i$. Since the set
$\{y_1,\dots,y_N\}$ is $r$-separated, the intersection of open
balls $B^G_{\frac{r}{2}}y_i\cap B^G_{\frac{r}{2}}y_j=
B_\cY(y_i,\frac{r}{2}) \cap B_\cY(y_j,\frac{r}{2})$ is empty if $j>i$. By
construction, the intersection $B^G_{\frac{r}{2}}y_i\cap P_j$ is empty
if $j<i$. Therefore $P_i$ contains $B^G_{\frac{r}{2}}y_i$, and
Assertion (ii) follows.

By construction, we have $\bigcup_{j=1}^{N}P_j\subset
\bigcup_{j=1}^{N} B^G_ry_j=\bigcup_{j=1}^{N} B_\cY(y_j,r)$, and the
later union contains $\cY(r)$, since the $r$-separated set
$\{y_1,\dots,y_N\}$ is maximal. Assertion (i) follows.

\medskip
For every $s>0$, let $n_s=\Big\lceil\frac{\ln(e^s-1)}{\ln q_v}
\Big\rceil\in\bZ$ and $\delta'_s=
\ln\big(\frac{1+q_v^{n_s}}{1+q_v^{n_s-1}}\big)>0$.  For all $\delta>0$
and $y\in\cY$, assume that there exists a point $z\in\partial_\delta
B_\cY(y,s)$.  Let $z'\in B_\cY(y,s)$ and $z'' \notin B_\cY(y,s)$ be such that
$d_\cY(z,z')+d_\cY(z,z'')<\delta$.  Since the operator norm on
$\cM_{d+1}(K_v)$ has values in $\{0\}\cup q_v^\bZ$, the set
$\{d_\cY(y,y'):y,y'\in\cY\}$ of values of the distance function
$d_\cY$ on $\cY$ is contained in $\{0\}\cup \{\ln(1+q_v^n):
n\in\bZ\}$.  Since $s\in\;]\ln(1+q_v^{n_s-1}), \ln(1+q_v^{n_s})]$, we
hence have $d_\cY(y,z')\leq\ln(1+q_v^{n_s-1})$ since $z'\in B_\cY(y,s)$
and $d_\cY(y,z'')\geq \ln(1+q_v^{n_s})$ since $z''\notin B_\cY(y,s)$.
Therefore by the triangle inequality and the inverse triangle
inequality, we have
\begin{align*}
\delta &> d_\cY(z,z')+d_\cY(z,z'')\geq d_\cY(z',z'')\geq
d_\cY(y,z'')-d_\cY(y,z')\\ & \geq 
\ln(1+q_v^{n_s})-\ln(1+q_v^{n_s-1})=\delta'_s\;.
\end{align*}
Hence $\partial_{\delta} B_\cY(y,s)$ is empty for every $\delta\in\;
  ]0,\delta'_s]$.

By Equation \eqref{eq:proprideltneigh}, for every $\delta>0$, we have
\[
\partial_{\delta}\cP\subset \bigcup_{j=1}^{N}\partial_{\delta}(B^G_ry_j)\cup
\bigcup_{j=1}^{N}\partial_{\delta}(B^G_{\frac{r}{2}}y_j)\;.
\]
Hence Assertion (iii) follows with $\delta_r=
\min\{\delta'_{\frac{r}{2}},r\}$.

Note that since the distance $d_G$ has values in $\{0\}\cup
\{\ln(1+q_v^n): n\in\bZ\}$, the open balls in $G$ are open and
compact, and since the canonical projection $G\ra\cY$ is open and
continuous, the subsets $P_i$ of $\cY$ are by construction open and
compact, and $P_\infty$ is closed and open.
\cqfd

\medskip
Let $\C$ be a countably generated $\sigma$-algebra of subsets of
$\cY$. Note that for every $j\in\bZ$, the $\sigma$-algebra $\ta^j\C$
is also countably generated and
\begin{equation}\label{eq:equivatom}
[y]_{a^j\C}=\ta^j\,[\ta^{-j}y]_{\C}\;.
\end{equation}
We say that $\C$ is {\it $\ta^{-1}$-descending} if $\ta\C$ is
contained in $\C$. If $\C$ is $\ta^{-1}$-descending, for all $y\in\cY$
and $j\in\bZ_{\geq 0}$, we have
\[
[y]_{\C}\subset[y]_{\ta^j\C}\;.
\]
Given a Borel probability measure $\mu$ on $\cY$ and a closed subgroup
$H$ of $G$, we say that $\C$ is {\it $H$-subordinate modulo $\mu$}
if for $\mu$-almost every $y\in\cY$, there exists $r=r_y\in\;]0,1]$
such that we have
\[
B_r^{H}y\subset [y]_{\C}\subset B_{1/r}^{H} y\;.
\]
If $\C$ is $U^+$-subordinate modulo $\mu$ and if furthermore $\mu$ is
$\ta$-invariant, since $\ta$ normalises $U^+$ and by Lemma
\ref{lem:contractboulUplus}, for every $j\in\bZ$, the $\sigma$-algebra
$\ta^j\C$ is also $U^+$-subordinate modulo $\mu$.

For every $\sigma$-algebra $\A$ of subsets of $\cY$, for all $a,b$ in
$\bZ\cup\{\pm\infty\}$ with $a<b$, we define a $\sigma$-algebra
$\A_a^b$ of subsets of $\cY$ by
\[
\A_a^b =\bigvee_{i=a}^{b}\;\ta^i \A=
\sigma\Big(\;\bigcup_{a\leq i\leq b}\ta^i \A\;\Big) \;.
\]
Note that if $\A$ is countably generated, then so is $\A_a^b$.

\bprop\label{prop:Apartconst}
For every $r\in\;]0,1[\,$, there exists a countably generated
sub-$\sigma$-algebra $\A^{U^+}$ of the Borel $\sigma$-algebra of $\cY$
such that
\begin{enumerate}
\item\label{Apartconst:prop1} the countably generated $\sigma$-algebra
  $\A^{U^+}$ is $\ta^{-1}$-descending,
\item\label{Apartconst:prop2} for every $y\in \cY(r)$, we have
  $[y]_{\A^{U^+}}\subset B_{r}^{U^+}y$,
\item\label{Apartconst:prop3} for every $y\in \cY$, we have
  $B_{\delta_r}^{U^+} y \subset [y]_{\A^{U^+}}$, where
  $\delta_r\in\;]0,r]$ is as in Lemma \ref{partconst}.
\end{enumerate}

\medskip
\noindent Let $\mu$ be a Borel $\ta$-invariant ergodic probability
measure on $\cY$ with $\mu(\cY(r))>0$.  Then $\A^{U^+}$ is
$U^+$-subordinate modulo $\mu$.
\eprop

\dem Fix $r\in\;]0,1[\,$. Let $\P=\{P_1,\dots,P_N,P_\infty\}$ be a
partition given by Lemma \ref{partconst} for this $r$.  We prove a
preliminary result on the countably generated sub-$\sigma$-algebra
$\sigma(\P)_0^\infty$.

\blemm\label{lem:calP0infty} For every $y\in \cY$, we have
$B_{\delta_r}^{U^+} y \subset [y]_{\sigma(\P)_0^\infty}$.
\elemm

\dem Let $h\in B_{\delta_r}^{U^+}$. Assume for a contradiction that
$hy \notin [y]_{\sigma(\P)_0^\infty}$. Then there exists $k\in
\bZ_{\geq 0}$ such that $\ta^{-k}hy$ and $\ta^{-k}y$ belong to
different atoms of the partition $\P$. Let $\alpha=\min\mb{r}>0$. By
Lemma \ref{lem:contractboulUplus}, we have
\begin{equation}\label{eqApartconst}
  d_{\cY} (\ta^{-k}hy, \ta^{-k}y) \leq d_G (\ta^{-k}h\ta^k, \id)=
  d_{U^+} (\ta^{-k}h\ta^k, \id)
  < \ln(1+(e^{\delta_r}-1)\,q_v^{-k\alpha})\leq \delta_r\leq r\;.
\end{equation}
It follows that both $\ta^{-k}hy$ and $\ta^{-k}y$ belong to the
$\delta_r$-boundary $\partial_{\delta_r}\P$ of $\P$.  But the set
$\partial_{\delta_r}\P$ is empty by Lemma \ref{partconst}
\eqref{partconst:prop3}, which gives a contradiction.
\cqfd

\medskip
By Lemma \ref{partconst}, for every $i\in\{1,\dots,N\}$, there exist
$y_i\in\cY(r)$ and a Borel subset $V_i$ of $G$ contained in $B_r^G$
such that $P_i=V_i y_i$. Let $\P^{U^+}$ be the sub-$\sigma$-algebra of
the Borel $\sigma$-algebra of $\cY$ generated by the subsets $P_\infty
\cap \pi^{-1}(W)$, where $W$ is a Borel subset of $\cX$, and the
subsets $((U^+B)\cap V_i)y_i$, where $i\in\{1,\dots,N\}$ and $B$ is a
Borel subset of $G$.  Then $\P^{U^+}$ is countably generated, since the
Borel $\sigma$-algebra of $\cX$ is countably generated and $U^+$ is a
closed subgroup of $G$. For every $y\in\cY$, the atom of $y$ for
$\P^{U^+}$ is equal to
\begin{equation}\label{eq:atomPUplus}
  [y]_{\P^{U^+}}=
  \left\{\begin{array}{ll} Uy& {\rm if}\;\;y\in P_\infty\\
  P_i\cap (B^{U^+}_ry)& {\rm if}\;\;
  \exists \;i\in\{1,\dots,N\},\; y\in P_i\;.\end{array}\right.
\end{equation}

Let us now define
\begin{equation}\label{eq:deficalAUplus}
  \A^{U^+}=(\P^{U^+})_0^\infty\;,
\end{equation}
which is a countably generated sub-$\sigma$-algebra of the Borel
$\sigma$-algebra of $\cY$, since so is $\P^{U^+}$. Note that
$\ta\,\A^{U^+}=(\P^{U^+})_1^\infty \subset \A^{U^+}$, which proves
Assertion \eqref{Apartconst:prop1}.

For every $y\in \cY(r)$, since $P_\infty\subset \cY-\cY(r)$ by Lemma
\ref{partconst} \eqref{partconst:prop1} and by Equation
\eqref{eq:atomPUplus}, we have $[y]_{\A^{U^+}}\subset [y]_{\P^{U^+}}
\subset B_r^{U^+}y$, which proves Assertion \eqref{Apartconst:prop2}.

In order to prove the last Assertion \eqref{Apartconst:prop3}, let us
take $y\in\cY$ and $h \in B_{\delta_r}^{U^+}$ and let us prove that
$hy \in [y]_{\A^{U^+}}$.  Since we have $hy \in
[y]_{\sigma(\P)_0^\infty}$ by Lemma \ref{lem:calP0infty}, for every
$k\geq 0$, there exists $i\in\{1,\dots,N,\infty\}$ such that the
points $\ta^{-k}y$ and $\ta^{-k}hy=\ta^{-k}h\ta^{k} (\ta^{-k}y)$ both
belong to $P_i \in \P$. If $i=\infty$, then by Equation
\eqref{eq:atomPUplus}, the points $\ta^{-k}y$ and $\ta^{-k}hy$ lie in
the same atom $[\ta^{-k}y]_{\P^{U^+}}=U\ta^{-k}y$ since
$\ta^{-k}h\ta^{k}\in U^+$.  Assume that $1\leq i\leq N$. Since $h \in
B_{\delta_r}^{U^+}$, it follows from Equation \eqref{eqApartconst}
that $\ta^{-k}h\ta^{k} \in B_r^{U^+}$.  Hence by Equation
\eqref{eq:atomPUplus}, the points $\ta^{-k} y$ and $\ta^{-k} hy$ lie in
the same atom $[\ta^{-k}y]_{\P^{U^+}}= P_i\cap (B_r^{U^+}\ta^{-k}y)$
of $\P^{U^+}$. This proves Assertion \eqref{Apartconst:prop3}.

Now let $\mu$ be an $\ta$-invariant ergodic probability measure on
$\cY$ with $\mu(\cY(r))>0$.  By ergodicity, for $\mu$-almost every
$y\in \cY$, there exists $k\in\bZ_{\geq 1}$ such that $\ta^{-k} y \in
\cY(r)$. Since $\ta^k\A^{U^+} \subset \A^{U^+}$, by Assertion (1) and
by Lemma \ref{lem:contractboulUplus}, we have
\[
[y]_{\A^{U^+}} \subset [y]_{\ta^k\A^{U^+} } = \ta^k[\ta^{-k} y]_{\A^{U^+}}
\subset \ta^k B_r^{U^+}\ta^{-k} y \subset
B_{\ln(1+(e^r-1)\,q_v^{k\max\mb{r}})}^{U^+} y\;.
\]
With Assertion (3), this proves that $\A^{U^+}$ is $U^+$-subordinate
modulo $\mu$.  \cqfd

\bigskip
Let us introduce some material before stating and proving our next
Lemma \ref{lem:algexiA}. The map $d_{K_v^{\,m},\mb{r}}:K_v^{\,m}\times
K_v^{\,m}\ra[0,+\infty[$ defined by
\begin{equation}\label{eq:defidistbfr}
\forall\;\bm{\xi},\bm{\xi}'\in K_v^{\,m},\;\;\;\;
d_{K_v^{\,m},\mb{r}}(\bm{\xi},\bm{\xi}')= \|\,\bm{\xi}-\bm{\xi}'\,\|_{\mb{r}}
\end{equation}
is an ultrametric distance on $K_v^{\,m}$, since the
$\mb{r}$-pseudonorm $\|\;\|_{\mb{r}}$ satisfies the ultrametric
inequality : for all $\bm{\xi},\bm{\xi}'\in K_v^{\,m}$, we have
\begin{equation}\label{eq:ultraineqpsudonorm}
  \|\,\bm{\xi}+\bm{\xi}'\,\|_{\mb{r}}\leq
  \max\{\|\,\bm{\xi}\,\|_{\mb{r}},\;\|\,\bm{\xi}'\,\|_{\mb{r}}\}\;,
\end{equation}
with equality if $\|\,\bm{\xi}\,\|_{\mb{r}}\neq
\|\,\bm{\xi}'\,\|_{\mb{r}}$.  Note that the map similar to
$d_{K_v^{\,m},\mb{r}}$ in the real case of \cite{KimKimLim21} is not a
distance if $m\geq 2$ for general $\mb{r}$. For every $\epsilon>0$, we
denote by $B^{K_v^{\,m},\mb{r}}_\epsilon$ the open ball of center $0$
and radius $\epsilon$ in $K_v^{\,m}$ for $d_{K_v^{\,m},\mb{r}}$.  Note
that the distance $d_{K_v^{\,m},\mb{r}}$ is bihölder equivalent to the
standard one: For all $\bm{\xi},\bm{\xi}'\in K_v^{\,m}$ such that
$\|\,\bm{\xi}-\bm{\xi}'\,\|\leq 1$, we have
\begin{equation}\label{eq:holderdistwithr}
\|\,\bm{\xi}-\bm{\xi}'\,\|^{\frac{1}{\min\mb{r}}}\leq
d_{K_v^{\,m},\mb{r}}(\bm{\xi},\bm{\xi}') \leq
\|\,\bm{\xi}-\bm{\xi}'\,\|^{\frac{1}{\max\mb{r}}}\;.
\end{equation}

We also endow the quotient space $\bT^m=K_v^{\,m}/R_v^{\,m}$ with the
quotient distance $d_{\bT^m,\mb{r}}$ of the distance
$d_{K_v^{\,m},\mb{r}}$ on $K_v^{\,m}$ defined by Equation
\eqref{eq:defidistbfr}. For every $A \in \cM_{m,n}(K_v)$, we denote by
$d_{U^+y_{A,0},\mb{r}}$ the distance on the orbit $U^+y_{A,0}=
\phi_A(\bT^m)$ defined by requiring the homeomorphism $\phi_A$ given
in Lemma \ref{lem:isomphiA} to be an isometry for the distances
$d_{\bT^m,\mb{r}}$ and $d_{U^+y_{A,0},\mb{r}}$, that is,
\[
\forall\;\bm{\theta},\bm{\theta}'\in\bT^m,\quad
d_{U^+y_{A,0},\mb{r}}(\phi_A(\bm{\theta}), \phi_A(\bm{\theta}'))=
d_{\bT^m,\mb{r}}(\bm{\theta},\bm{\theta}'))\;.
\]

Using the identification $w\mapsto \wh w$ between $K_v^{\,m}$ and
$U^+$ (see Subsection \ref{subsec:homogdyn}), for every $\epsilon>0$,
we denote by $B_\epsilon^{U^+,\mb{r}}$ the open ball of radius
$\epsilon$ in $U^+$ centered at the identity element for the distance
$d_{U^+,\mb{r}}$ on $U^+$ corresponding to the distance
$d_{K_v^m,\mb{r}}$ on $K_v^m$.  The map $u\mapsto u\,y_{A,0}$ from
$U^+$ onto $U^+y_{A,0}$ is $1$-Lipschitz and locally isometric for the
distances $d_{U^+,\mb{r}}$ and $d_{U^+y_{A,0},\mb{r}}$. Improving
Lemma \ref{lem:contractboulUplus}, for all $\epsilon>0$ and $k\in\bZ$,
we have
\begin{equation}\label{eq:contracboulUplusexact}
  \ta^{-k}B_\epsilon^{U^+,\mb{r}}\ta^k= B_{\epsilon \,q_v^{-k}}^{U^+,\mb{r}}\;.
\end{equation}

Again using the (locally compact) topological group identification
$w\mapsto \wh w$ between $(K_v^{\,m},+)$ and $U^+$, we endow $U^+$
with the Haar measure $m_{U^+}$ which corresponds to the normalized
Haar measure $\vol_v^m$ of $K_v^{\,m}$ (see Section
\ref{subsec:funcfield}).  For every $j\in\bZ$, the Jacobian
$\operatorname{Jac}_j$ with respect to the measure $m_{U^+}$ of the
homeomorphism $\varphi_j:u\mapsto \ta^j \;u\;\ta^{-j}$ from $U^+$
to $U^+$ (which is constant since $\varphi_j$ is a group automorphism
and $m_{U^+}$ is bi-invariant) is easy to compute: we have
\begin{equation}\label{eq:caljac}
    \operatorname{Jac}_j=q_v^{\;j\;|\mb{r}|}\;.
\end{equation}

\medskip
Using Equation \eqref{eq:deficalAUplus}, we consider the following
tail $\sigma$-algebra:
\begin{equation}\label{eq:tailalg}
  \A^{U^+}_\infty =
  \bigcap_{k=1}^{\infty} \bigvee_{i=k}^{\infty} \ta^{i}\A^{U^+}
  = \bigcap_{k=1}^{\infty} (\A^{U^+})_k^\infty
= \bigcap_{k=1}^{\infty} (\P^{U^+})_k^\infty.
\end{equation}
This $\sigma$-algebra may not be countably generated, but it is
strictly $\ta$-invariant, i.e., we have $\ta\A^{U^+}_\infty =
\A^{U^+}_\infty =\ta^{-1}\A^{U^+}_\infty$. We will use this
$\sigma$-algebra in order to observe the relative entropy of $\ta^{-1}$.

\blemm\label{lem:algexiA} For every $r\in\;]0,1[\,$, let $\A^{U^+}$ be
    as in Proposition \ref{prop:Apartconst} and $\A^{U^+}_\infty$ be
    as in Equation \eqref{eq:tailalg}. Let $\mu$ be an $\ta$-invariant
    ergodic probability measure on $\cY$. Then
\[
h_\mu(\ta^{-1}|\,\A^{U^+}_\infty)\leq |\mb{r}|\;.
\]
Furthermore, if $\mu(\cY(r))>0$, then
\[
h_\mu(\ta^{-1}|\,\A^{U^+}_\infty) = H_\mu(\A^{U^+}|\,\ta\,\A^{U^+})\;.
\]
\elemm

\dem Let us prove the first assertion. By \cite[Prop. 7.44]{EinLin10},
there exists a countable Borel-measurable partition $\G$ with finite
entropy which is a generator for the invertible transformation $\ta$
modulo $\mu$, such that $\sigma(\G)_0^\infty$ is $\ta^{-1}$-descending
and $G^+$-subordinate modulo $\mu$. It follows from
\cite[Kolmogorov-Sinai Theo.~2.20 and Prop.~2.19 (8)]{EinLinWar22}
that
\[
h_{\mu}(\ta^{-1}|\,\A^{U^+}_\infty)=h_{\mu}(\ta^{-1},\G\,|\,\A^{U^+}_\infty)=
H_{\mu}(\sigma(\G)\,|\,\sigma(\G)_{1}^{\infty}\vee \A^{U^+}_\infty)\;.
\]
Using the continuity and monotonicity of entropy \cite[Prop.~2.12 and
  Prop.~2.13 (2)]{EinLinWar22}, we have
\[
\begin{split}
H_{\mu}(\sigma(\G)\,|\,\sigma(\G)_{1}^{\infty}\vee \A^{U^+}_\infty)
&=\lim_{\ell\to\infty} H_{\mu}(\sigma(\G)\,|\,
\sigma(\G)_{1}^{\infty}\vee (\A^{U^+})_\ell^\infty)\\ &\leq \lim_{\ell\to\infty}
H_{\mu}(\sigma(\G)_{0}^{\infty}\vee (\A^{U^+})_\ell^\infty\,|\,
\ta(\sigma(\G)_{0}^{\infty}\vee (\A^{U^+})_\ell^\infty))\;.
\end{split}
\]
Note that for every integer $\ell\geq 1$ the $\sigma$-algebra
$\sigma(\G)_0^\infty\vee(\A^{U^+})_\ell^\infty$ is countably generated,
$\ta^{-1}$-descending, and $U^+$-subordinate since
$[y]_{(\A^{U^+})_\ell^\infty}\subset Uy$ for all $y\in\cY$ and since
$\sigma(\G)_0^\infty$ is $G^+$-subordinate.  Thus by
\cite[Prop. 7.34]{EinLin10} (recalling that we are using logarithms
with base $q_v$), we have
\[
H_\mu(\sigma(\G)_0^\infty\vee(\A^{U^+})_\ell^\infty|
\,\ta\,(\sigma(\G)_0^\infty\vee(\A^{U^+})_\ell^\infty))
=\lim_{k\to \infty}
\frac{\log_{q_v}\mu_{x}^{U^+}(\ta^{k}B_{1}^{U^+}\ta^{-k})}{k}\;,
\]
where $\mu_{x}^{U^+}$ is the leafwise measure of $\mu$ at $x\in \cY$
with respect to $U^+$ as defined in \cite[Theo.~6.3]{EinLin10}.  By
\cite[Theo.~6.30]{EinLin10} (which applies since $U^+$ is abelian,
hence unimodular) and by Equation \eqref{eq:caljac} (see also \cite[\S
  7.42]{EinLin10}), we have
\[
\limsup_{k\to\infty}
\frac{\mu_{x}^{U^+}(\ta^{k}B_{1}^{U^+}\ta^{-k})}{k^2 \,q_v^{k|\mb{r}|}}=0\;.
\]
Hence we have
\[
\lim_{k\to \infty}
\frac{\log_{q_v}\mu_{x}^{U^+}(\ta^{k}B_{1}^{U^+}\ta^{-k})}{k} \leq |\mb{r}|\;.
\]
This proves the first assertion of the lemma.

In order to prove the second assertion, let us take an increasing
sequence of finite partitions $(\P_k^{U^+})_{k\geq 1}$ of $\cY$ such
that the union $\bigcup_{k\geq 1}\P_k^{U^+}$ generates $\P^{U^+}$,
which is possible since the $\sigma$-algebrba $\P^{U^+}$ is countably
generated.  Since $\mu$ is ergodic and $\mu(\cY(r))>0$, for
$\mu$-almost every $y\in \cY$, there exists an increasing sequence of
positive integers $(k_i)_{i\geq 1}$ such that $\ta^{k_i}y\in \cY(r)$.
By Proposition \ref{prop:Apartconst} \eqref{Apartconst:prop2}, we have
$[\ta^{k_i}y]_{\A^{U^+}}\subset B_r^{U^+} \ta^{k_i}y$ for all $i\geq
1$. Hence it follows from Equation \eqref{eq:equivatom} and Lemma
\ref{lem:contractboulUplus} that
\[
  [y]_{(\P^{U^+})_{-k_i}^{\infty}}=\ta^{-k_i}[\ta^{k_i}y]_{(\P^{U^+})_{0}^{\infty}}
  \subset \ta^{-k_i} B_r^{U^+} \ta^{k_i}y \subset
  B_{\ln(1+(e^r-1)\,q_v^{-k_i \min \mb{r}})}^{U^+} y.
\]
Taking $i\to \infty$, we have $[y]_{(\P^{U^+})_{-\infty}^{\infty}}=
\{y\}$ for $\mu$-almost every $y\in \cY$. This means that
$(\P^{U^+})_{-\infty}^\infty = \B_{\cY} $ modulo $\mu$, where
$\B_{\cY}$ is the Borel $\sigma$-algebra of $\cY$.  It follows that
$\bigvee_{k=1}^{\infty}\sigma(\P_k^{U^+})_{-\infty}^{\infty}=
(\P^{U^+})_{-\infty}^\infty=\B_{\cY}$ modulo $\mu$, and
$\sigma(\P_k^{U^+})_{-\infty}^{\infty}\subseteq
\sigma(\P_{k+1}^{U^+})_{-\infty}^\infty$ for each $k\geq 1$.  Using
respectively the Kolmogorov-Sinai theorem for sequences
\cite[Theo.~2.20]{EinLinWar22}, the Future formula \cite[Prop.~2.19
  (8)]{EinLinWar22} and the continuity of entropy
\cite[Prop.~2.12]{EinLinWar22}, we have
\[\begin{split}
h_{\mu}(\ta^{-1}|\,\A^{U^+}_\infty)&=
\lim_{k\to\infty}h_{\mu}(\ta^{-1},\P_k^{U^+}|\,\A^{U^+}_\infty)
=\lim_{k\to\infty}
H_{\mu}(\P_k^{U^+}|\,\sigma(\P_k^{U^+})_1^\infty \vee \A_\infty^{U^+})\\
&=H_{\mu}(\P^{U^+}|\,(\P^{U^+})_1^\infty \vee \A_\infty^{U^+})
=H_{\mu}(\P^{U^+}|\,(\P^{U^+})_1^\infty)\\
&=H_{\mu}(\A^{U^+}|\,\ta\,\A^{U^+})\;.
\end{split}\]
This proves the second assertion of the lemma.
\cqfd

\medskip
Let us introduce some more material before stating and proving our
final Proposition \ref{prop:effEL} of Subsection
\ref{subsec:effvarprin}. Let $\cA$ be a countably generated
sub-$\sigma$-algebra of the Borel $\sigma$-algebra of $\cY$.  For all
$j\in\bZ_{\geq 0}$ and $y\in \cY$, let
\[
r_{y,j}=\inf\big\{r>0:
[y]_{\ta^j\cA} \subset B^{U^+,\mb{r}}_r \,y\big\}\in[0,+\infty]\;,
\]
with the usual convention that $\inf \emptyset =+\infty$. Note that
$r_{y,j}$ depends measurably on $y$ and that $r_{\ta^{-j}y,0}=
r_{y,j}\,q_v^{-j}$ by Equation \eqref{eq:contracboulUplusexact}. We
now define
\begin{equation}\label{eq:defishape}
  V_y^{\ta^j\cA}= \{u \in \overline{B}^{U^+,\mb{r}}_{r_{y,j}}:
    u\, y\in [y]_{\ta^j\cA}\}\;,
\end{equation}
which is a Borel subset of $U^+$, called the {\it $U^+$-shape of the
  atom $[y]_{\ta^j\cA}$}. Note that for every $j\in\bZ_{\geq 0}$, we
have
\[
V_y^{\ta^j\cA}= \ta^j\;V_{\ta^{-j}y}^{\cA}\;\ta^{-j}\;.
\]
Let us define a Borel-measurable family $\big(\tau_y^{\ta^j\cA}
\big)_{y\in \cY}$ of Borel measures on $\cY$, that we call the {\it
  $U^+$-subordinate Haar measure of $\ta^j\cA$}, as follows:

$\bullet$~
if $m_{U^+}(V_y^{\ta^j\cA})$ is equal to $0$ or $\infty$,
we set $\tau_y^{\ta^j\cA}=0$,

$\bullet$~ otherwise,  $\tau_y^{\ta^j\cA}$ is the
push-forward of the normalized measure $\frac{1}{m_{U^+}(V_y^{\ta^j\cA})}
\;m_{U^+}\!\!\mid_{V_y^{\ta^j\cA}}$ by the map $u\mapsto u\, y$.

Now let $\mu$ be a Borel $\ta$-invariant probability measure on $\cY$,
such that $\cA$ is $U^+$-subordinate modulo $\mu$. In particular, for
$\mu$-almost every $y\in\cY$, the atom $V_y^{\ta^j\cA}$ has positive
and finite $m_{U^+}$-measure, hence the measure $\tau_y^{\ta^j\cA}$ is
a probability measure with support in $[y]_{\ta^j\cA}$. Furthermore,
if $z\in [y]_{\ta^j\cA}$ then there exists $u\in U^+$ such that $z=u\,
y$, $V_z^{\ta^j\cA} =V_y^{\ta^j\cA}\,u^{-1}$, and $\tau_z^{\ta^j\cA}
=\tau_y^{\ta^j\cA}$, by the right-invariance of $m_{U^+}$.

\medskip
The following proposition is a function field analog of the effective
real case version \cite[Prop.~2.10, \S 2.4]{KimKimLim21} of \cite[\S
  7.55]{EinLin10}.

\bprop\label{prop:effEL} Let $\mu$ be a Borel $\ta$-invariant ergodic
probability measure on $\cY$ and let $\A$ be a countably generated
sub-$\sigma$-algebra of the Borel $\sigma$-algebra of $\cY$ which is
$\ta^{-1}$-descending and $U^{+}$-subordinate modulo $\mu$.  Fix
$j\in\bZ_{\geq 1}$ and a $U^{+}$-saturated Borel subset $K'$ of $\cY$.
Suppose that there exists $\epsilon>0$ such that $[z]_{\A}\subset
B_{\epsilon}^{U^+,\mb{r}} z$ for every $z\in K'$.  Then we have
\[
H_\mu(\A|\ta^{j}\A)\leq j\,|\mb{r}|+\int_\cY \log\tau_y^{\ta^{j}\A}
((\cY- K')\cup B^{U^+,\mb{r}}_{\epsilon}\Supp\mu)\;d\mu(y).
\]
\eprop

\dem We fix $\mu$, $\cA$, $j$, $K'$ and $\epsilon$ as in the
statement.  By for instance \cite[Theo.~5.9]{EinLin10}, let
$\big(\mu_y^{\ta^j\cA}\big)_{y\in \cY}$ be a measurable family of
conditional measures of $\mu$ with respect to $\ta^j\cA$, so that for
$\mu$-almost every $y\in\cY$, the measure $\mu_y^{\ta^j\cA}$ is a
probability measure on $\cY$ giving full measure to the atom
$[y]_{\ta^j\cA}$, with $\mu_z^{\ta^j\cA} =\mu_y^{\ta^j\cA}$ if
$z\in[y]_{\ta^j\cA}$, and such that the following {\it disintegration
  formula} holds true:
\begin{equation}\label{eq:disinting}
  \mu=\int_{y\in  \cY}\mu_y^{\ta^j\cA}\;d\mu(y)\;.
\end{equation}

Let $p_\mu:y\mapsto \mu_y^{\ta^j\cA}([y]_{\cA})$ and $p_\tau: y\mapsto
\tau_y^{\ta^j\cA}([y]_{\cA})$, which are nonnegative and measurable
functions on $\cY$. Since $\cA$ is $\ta^{-1}$-descending and
$U^+$-subordinate modulo $\mu$, the atom $[y]_{\cA}$ contains an open
neighborhood of $y$ in the atom $[y]_{\ta^j\cA}$ for $\mu$-almost
every $y\in\cY$. In particular, the function $p_\tau$ is $\mu$-almost
everywhere positive.

Since $\cA$ is countably generated and $\ta^{-1}$-descending, for every
$y\in\cY$, the atom of $y$ for $\ta^j\cA$ is countably partitioned
into atoms for $\cA$ up to measure $0$, that is, there exist a finite
or countable subset $I_y$ of $[y]_{\ta^j\cA}$ and a
$\mu_y^{\ta^j\cA}$-measure zero subset $N_y$ of $[y]_{\ta^j\cA}$ such
that
\begin{equation}\label{eq:partitionatom}
[y]_{\ta^j\cA}=N_y\sqcup\bigsqcup_{x\in I_y}\;[x]_{\cA}\;.
\end{equation}
Let $I'_y=\{x\in I_y:[x]_{\cA}\cap \operatorname{Supp}\mu\neq
\emptyset\}$.

\blemm\label{lem:atomZZZ} Let $x\in I_y$
\begin{enumerate}
\item[(1)] If $x\notin I'_y$, then $\mu_y^{\ta^j\cA}([x]_{\cA})= 0$.
\item[(2)] If $x\in I'_y$, then $[x]_{\cA}$ is contained in $(\cY
  -K')\cup B_\epsilon^{U^+,\mb{r}}\operatorname{Supp}\mu $.
\end{enumerate}
\elemm

\dem
(1) This follows since $\operatorname{Supp} \mu_y^{\ta^j\cA}$ is
contained in $\operatorname{Supp}\mu$.

\medskip
(2) If $x\in I'_y$, there exists $z\in [x]_{\cA}\cap
\operatorname{Supp}\mu$. For every $z'\in [x]_{\cA}$, we have either
$z'\in \cY -K'$ or $z'\in K'$. In the second case, since $\cA$ is
$U^+$-subordinate and $K'$ is $U^+$-saturated, we have $z\in
[x]_{\cA}=[z']_{\cA}\subset U^+z'\subset K'$. Hence by the
assumption of Proposition \ref{prop:effEL}, we have $z'\in
[x]_{\cA}=[z]_{\cA}\subset B_\epsilon^{U^+,\mb{r}} z\subset
B_\epsilon^{U^+,\mb{r}}\operatorname{Supp}\mu $, which proves the
result.
\cqfd

\medskip
By the definition of the $U^+$-subordinate Haar measure of $\ta^j\cA$,
for $\mu$-almost every $y\in\cY$, we have
$$
p_\tau(y)=\frac{m_{U^+}(V_y^\cA)}{m_{U^+}(V_y^{\ta^j\cA})}=
\frac{m_{U^+}(V_y^\cA)}{m_{U^+}(\ta^j\;V_{\ta^{-j}y}^{\cA}\;\ta^{-j})}
=\frac{m_{U^+}(V_y^\cA)}{\operatorname{Jac}_j\;m_{U^+}(V_{\ta^{-j}y}^{\cA})}\;.
$$
Hence, by the $\ta$-invariance of $\mu$ and by Equation
\eqref{eq:caljac}, we have
$$
\int_{z\in\cY} \log_{q_v} p_\tau(z)\;d\mu(z)=
- \log_{q_v}\operatorname{Jac}_j=-\,j\,|\mb{r}|\;.
$$
Respectively

$\bullet$~ by the definition of the conditional entropy in Equation
\eqref{eq:deficondentrop},

$\bullet$~ by the disintegration formula \eqref{eq:disinting},

$\bullet$~ since $\mu_y^{\ta^j\cA}$ gives full measure to
$[y]_{\ta^j\cA}$ which is partitionned as in Equation
\eqref{eq:partitionatom}, and by Lemma \ref{lem:atomZZZ} (1),

$\bullet$~ since when $z$ varies in $[x]_{\cA}\subset [y]_{\ta^j\cA}$,
the values $p_\mu(z)=\mu_z^{\ta^j\cA}([z]_{\cA})= \mu_y^{\ta^j\cA}
([x]_{\cA})$ and $p_\tau(z)= \tau_z^{\ta^j\cA}([z]_{\cA}) =
\tau_y^{\ta^j\cA}([x]_{\cA})$ are constant,

$\bullet$~ by the concavity property of the logarithm,

$\bullet$~ by Lemma \ref{lem:atomZZZ} (2),

\noindent we hence have
\begin{align*}
  &H_\mu(\cA\,|\,\ta^j\cA)\;-j\;|\mb{r}|\\=\;& -\int_{z\in\cY}
  \big(\log_{q_v} p_\mu(z)- \log_{q_v} p_\tau(z)\big)\;d\mu(z) \\ =\;&
  \int_{y\in\cY}\int_{z\in\cY} \big(\log_{q_v} p_\tau(z)- \log_{q_v}
  p_\mu(z)\big)\;d\mu_y^{\ta^j\cA}(z)\;d\mu(y) \\ =\;&
  \int_{y\in\cY}\sum_{x\in I'_y}\int_{z\in[x]_{\cA}} \big(\log_{q_v}
  p_\tau(z)- \log_{q_v} p_\mu(z)\big)\;d\mu_y^{\ta^j\cA}(z)
  \;d\mu(y)\\ =\;&\int_{y\in\cY}\sum_{x\in I'_y} \log_{q_v}
  \frac{\tau_y^{\ta^j\cA}([x]_{\cA})}{\mu_y^{\ta^j\cA} ([x]_{\cA})}
  \;\mu_y^{\ta^j\cA}([x]_{\cA}) \;d\mu(y)
  \\ \leq\;&\int_{y\in\cY}\log_{q_v} \Big(\sum_{x\in I'_y}
  \tau_y^{\ta^j\cA}([x]_{\cA})\Big) \;d\mu(y)\\ \leq\;& \int_{y\in\cY}
  \log_{q_v}\big(\tau_y^{\ta^j\cA}((\cY-K')\cup
  B_\epsilon^{U^+,\mb{r}}\Supp\mu)\big)\;d\mu(y)\;.
\end{align*}
This proves the result. \cqfd


\section{Upper bound on the Hausdorff dimension of $\Bad_A(\epsilon)$}
\label{sec:upperbound}

\subsection{Constructing measures with large entropy}
\label{subsec:consmeas}

In this subsection, we construct, as in \cite[Prop.~4.1]{KimKimLim21}
in the real case, an $\ta$-invariant probability measure on $\cY$
giving an appropriate lower bound on the conditional entropy of $\ta$
relative to the $\sigma$-algebra $\A^{U^+}_\infty$ defined in Equation
\eqref{eq:tailalg} using the $\sigma$-algebra $\A^{U^+}$ constructed
in Proposition \ref{prop:Apartconst} (see Equation
\eqref{eq:deficalAUplus}).

For any point $x$ in a measurable space, we denote by $\Delta_x$ the
unit Dirac measure at $x$. We denote by
$\stackrel{*}{\rightharpoonup}$ the weak-star convergence of Borel
measures on any locally compact space.

\medskip
Let us denote by $\ov{\cX}=\cX\cup\{\infty_\cX\}$ and $\ov{\cY}=
\cY\cup\{\infty_\cY\}$ the one-point compactifications of $\cX$ and
$\cY$, respectively.  We denote by $\ov{\pi}:\ov{\cY}\ra \ov{\cX}$ the
unique continuous extension of the natural projection $\pi:\cY\ra\cX$,
mapping $\infty_\cY$ to $\infty_\cX$. The left actions of $\ta$ on
$\cX$ and $\cY$ continuously extend to actions on $\ov{\cX}$ and
$\ov{\cY}$ fixing the points at infinity $\infty_\cX$ and
$\infty_\cY$. For every countably generated $\sigma$-algebra $\A$ of
subsets of $\cX$ or $\cY$, we denote by $\ov{\A}$ the countably
generated $\sigma$-algebra of subsets of $\ov{\cX}$ or $\ov{\cY}$
generated by $\cA$ and its point at infinity.  For a finite partition
$\Q= \{Q_1,\dots,Q_N,Q_\infty \}$ of $\cY$ with only one unbounded
atom $Q_\infty$, we denote by $\ov{\Q}$ the finite partition
$\{Q_1,\dots, Q_N, \ov{Q}_{\infty}= Q_\infty\cup \{\infty_\cY\}\}$ of
$\ov{\cY}$.  Note that $\overline{\bigvee_{i=a}^{b}\ta^{-i}\Q}=
\bigvee_{i=a}^{b} \ta^{-i}\,\ov{\Q}$ for all $a,b$ in $\bZ$ with
$a<b$.

For every $\eta\in[0,1]$, we say that an element $x\in\cX$ has {\it
  $\eta$-escape of mass on average} under the action of $\ta$ if for
every compact subset $Q$ of $\cX$,
\[
\liminf_{N\to\infty} \frac{1}{N} \;\card\,\big\{\ell\in \{1, \cdots, N\} :
\ta^\ell x \notin Q \big\} \geq \eta\;.
\]
When $\eta=1$, as defined in the Introduction and in Proposition
\ref{prop:critdynescapeonaver}, we say that $x$ diverges on average in
$\cX$ under the action of $\ta$. For every $A\in \cM_{m,n}(K_v)$, we
denote by $x_A=u_A R^{\,m}_v\in \cX$ its associated unimodular lattice
(see Section \ref{subsec:Dani}), and by $\eta_A\in[0,1]$ the upper
bound of the elements $\eta\in[0,1]$ such that $x_A$ has $\eta$-escape
of mass on average. Note that this upper bound is actually a maximum.

\bprop\label{prop:measconst} For every $A\in \cM_{m,n}(K_v)$, there
exists a Borel probability measure $\mu_A$ on $\ov{\cX}$
with $\mu_A(\cX)=1-\eta_A$ such that for every $\epsilon>0$, there
exists an $\ta$-invariant Borel probability measure $\ov{\mu}$
on $\ov{\cY}$ satisfying the following properties.
\begin{enumerate}
\item\label{supp} The support of $\ov{\mu}$ is contained in
  $\cL_\epsilon \cup \{\infty_\cY\}$, where the $\epsilon$-compact
  part $\cL_\epsilon$ is defined in Equation \eqref{eq:defiLsubeps}.
\item\label{cusp} We have $\ov{\pi}_{*}\ov{\mu}=\mu_A$. In particular,
  there exists an $\ta$-invariant Borel probability measure $\mu$ on
  $\cY$ such that
$$
\ov{\mu}=(1-\eta_A)\mu + \eta_A \Delta_{\infty_{\cY}}.
$$
\item\label{entropy} For every $r\in\;]0,1[\,$, let $\A^{U^+}$ be the
    $\sigma$-algebra of subsets of $\cY$ constructed in Proposition
    \ref{prop:Apartconst} and let $\A^{U^+}_\infty$ be as in Equation
    \eqref{eq:tailalg}. Then
\[
h_{\ov{\mu}}\big(\ta^{-1}|\,\ov{{\A^{U^+}_\infty}}\,\big)
=h_{\ov{\mu}}\big(\ta\,|\,\ov{{\A^{U^+}_\infty}}\,\big)
\geq |\mb{r}| (1-
\eta_A) - \max\mb{r}\;(m -\dimH \mb{Bad}_A(\epsilon))\;.
\]
\end{enumerate}
\end{prop}

\dem Since $x_A$ has $\eta_A$-escape of mass on average but does not
have $(\eta_A+\delta)$-escape of mass on average for any $\delta>0$,
there exists an increasing sequence of positive integers
$(k_i)_{i\in\bZ_{\geq 1}}$ such that, for the weak-star convergence of
Borel probability measures on the compact space $\ov{\cX}$, as
$i\ra+\infty$, we have
\begin{equation}\label{eq:definmuA}
\frac{1}{k_i}\sum_{k=0}^{k_i-1}\Delta_{\,\ta^k x_A}
\;\stackrel{*}{\rightharpoonup}\;\mu_A\;,
\end{equation}
and $\mu_A$ is a Borel probability measure on $\ov{\cX}$ with
$\mu_{A}(\cX)=1-\eta_A$. This is equivalent to $\mu_A(\{\infty_\cX\})
=\eta_A$.

Let $\epsilon>0$. For every $T\in\bZ_{\geq 0}$, with the notation of
Subsection \ref{subsec:Dani} (see in particular Equations
\eqref{eq:defiLsubeps} and \eqref{eq:defiphiA}), let
\[
R_{T}=\{\bm{\theta}\in\bT^m:
\forall k \geq T,\ta^k \phi_{A}(\bm{\theta})\in
\cL_\epsilon\}\cap\Bad_A(\epsilon)\;.
\]
By Proposition \ref{prop:bad}, since a countable subset of
$K_{\nu}^{\,m}$ has Hausdorff dimension $0$, we have $\dimH
\big(\bigcup_{T=1}^\infty R_{T}\big) = \dimH \Bad_A(\epsilon)$.  Thus,
for every $j\in \bZ_{\geq 1}$, there exists $T_j \in \bZ_{\geq 0}$
satisfying
\[
\dimH R_{T_j}\geq \dimH \Bad_A(\epsilon)-\frac{1}{j}\;.
\]

For all $i,j\in \bZ_{\geq 1}$ such that $k_i\ge T_j$, let $S_{i,j}$ be
a maximal $q_v^{-k_i}$-separated subset of $R_{T_j}$ for the distance
$d_{\bT^m,\mb{r}}$ defined after Equation \eqref{eq:holderdistwithr}.
Then $R_{T_j}$ can be covered by $\card\; S_{i,j}$ open balls of
radius $q_v^{-k_i}$ for $d_{\bT^m,\mb{r}}$. Each open ball of radius
$q_v^{-k_i}$ for $d_{\bT^m,\mb{r}}$ can be covered by $\prod_{j=1}^m
q_v^{-k_i r_j} / q_v^{-k_i \max\mb{r}} =q_v^{k_i (m \max\mb{r} -
  |\mb{r}|)}$ open balls of radius $q_v^{-k_i \max\mb{r}}$ with
respect to the standard distance $d_{\bT^m}$ (defining the Hausdorff
dimension of subsets of $\bT^m$). Since the lower Minskowski dimension
is at least equal to the Hausdorff dimension, we have
\[
\liminf_{i\to\infty}\frac{\log_{q_v} \big(q_v^{k_i (m\max\mb{r} -
    |\mb{r}|)}\card\;S_{i,j}\big)}{-\log_{q_v} \big(q_v^{-k_i \max
    \mb{r}}\big)}\geq \dimH R_{T_j} \geq \dimH
\Bad_A(\epsilon)-\frac{1}{j}\;,
\]
which implies that
\begin{equation}\label{ineqDim}
  \liminf_{ i \to \infty} \frac{ \log_{q_v} \card\; S_{i,j}}{k_i}
  \geq |\mb{r}|-\max\mb{r}\;\big(m + \frac{1}{j} -
  \dimH \Bad_A(\epsilon)\big)\;.
\end{equation}

Let us define the Borel probability measures
\[
\nu_{i,j}=\frac{1}{\card\; S_{i,j}} \sum_{\bm{\theta}\in S_{i,j}}
\Delta_{\phi_A(\bm{\theta})}\;,
\]
which is the normalized counting measure on the finite subset
$\phi_A(S_{i,j})$ of the $U^+$-orbit $\phi_A(\TT^m)= U^+y_{A,0}
\subset\pi^{-1}(x_A)$, and
\[
\wt\nu_{i,j}= \frac{1}{k_i} \sum_{0\leq k\leq k_i-1}
\ta_*^k\nu_{i,j}\;,
\]
which is the average of the previous one on the first $k_i$ points of
the $\ta$-orbit. Since $\ov{\cY}$ is compact, extracting diagonally a
subsequence if necessary, we may assume that $\wt\nu_{i,j}$ weak-star
converges as $i\ra+\infty$ towards an $\ta$-invariant Borel
probability measure $\wt\mu_j$, and that $\wt\mu_j$ weak-star
converges as $j\ra+\infty$ towards an $\ta$-invariant Borel
probability measure $\ov{\mu}$.  Let us prove that $\ov{\mu}$
satisfies the three assertions of Proposition \ref{prop:measconst}.

\medskip\noindent \eqref{supp} For all $k\geq T_j$ and $\bm{\theta}\in
S_{i,j}\subset R_{T_j}$, we have $\ta^k \phi_A(\bm{\theta})\in
\cL_\epsilon$ by the definition of $R_{T_j}$. Since
$\ta^k_{*}\nu_{i,j}$ is a probability measure, we hence have
\[
  \wt\nu_{i,j}(\cY-\cL_\epsilon)=\frac{1}{k_i}
  \sum_{k=0}^{k_i-1}\ta^k_{*}\nu_{i,j}(\cY-\cL_\epsilon)= \frac{1}{k_i}
  \sum_{k=0}^{T_j}\ta^k_{*}\nu_{i,j}(\cY-\cL_\epsilon)
\leq \frac{T_j}{k_i}\;.
\]
Since $\cL_\epsilon\cup\{\infty_\cY\}$ is closed in $\ov{\cY}$ and by
taking limits first as $i\ra+\infty$ then as $j\ra+\infty$, we
therefore have $\ov{\mu}(\cY- \cL_\epsilon) =0$. This proves Assertion
\eqref{supp}.

\medskip\noindent \eqref{cusp} Since $\phi_A(S_{i,j})$ is contained in
the fiber above $x_A$ of $\ov{\pi}$ and since $\nu_{i,j}$ is a probability
measure, we have $\ov{\pi}_*\nu_{i,j}=\Delta_{x_A}$. By the linearity and
equivariance of $\ov{\pi}_*$, we hence have
\[
\ov{\pi}_*\wt\nu_{i,j}=
\frac{1}{k_i} \sum_{0\leq k\leq k_i-1} \ta_*^k\;\ov{\pi}_*\,\nu_{i,j}=
\frac{1}{k_i} \sum_{0\leq k\leq k_i-1} \Delta_{\ta^kx_A}\;.
\]
By the weak-star continuity of $\ov{\pi}_*$ and Equation
\eqref{eq:definmuA}, we thus have
\[
\ov{\pi}_*\ov{\mu}=
\lim_{j\ra+\infty}\lim_{i\ra+\infty}\ov{\pi}_*\wt\nu_{i,j}=
\lim_{j\ra+\infty}\mu_A=\mu_A\;.
\]
Note that the point at infinity $\infty_\cY$ is an isolated point in
the support of $\ov{\mu}$ by Assertion (1), since $\cL_\epsilon$ is
compact. We hence have
\begin{equation}\label{eqcusp}
  \ov{\mu}(\{\infty_\cY\})=\ov{\mu}(\,\ov{\pi}^{-1}(\{\infty_\cX\}))=
  \mu_A(\{\infty_\cX\})=\eta_A\;.
\end{equation}
This proves Assertion \eqref{cusp}.

\medskip\noindent \eqref{entropy}
Suppose that $\Q$ is any finite Borel-measurable partition of $\cY$
satisfying
\begin{enumerate}
\item[(i)]
  the partition $\Q$ contains an atom $Q_\infty$ of the form
  $\pi^{-1}(Q_\infty^*)$, where $\cX- Q_\infty^*$ has
  compact closure,
\item[(ii)] there exists $\ell_0 \in\bZ_{\geq 1}$ such that for every atom
  $Q\in \Q$ different from $Q_\infty$ and for every $y'\in Q$, we have
  $\diam\; (U^{+}y' \cap Q) < q_v^{-\ell_0\max\mb{r}}$ for the distance
  $d_{U^+y',\|\;\|}$ defined in Lemma \ref{lem:isomphiA}.
\item[(iii)] for all $Q\in\Q$ and $j\in\bZ_{\geq 1}$, we have
  $\wt\mu_j(\partial Q)=0$ and $\ov{\mu}(\partial Q)=0$.
\end{enumerate}
We first prove the following entropy bound: For every $M\in\bZ_{\geq 1}$,
\begin{equation}\label{staticentropybound}
  \frac{1}{M}H_{\ov\mu}\big(\,\overline{\sigma(\Q^{(M)})}\,|
  \,\overline{\A^{U^+}_\infty}\,\big)\geq 
  |\mb{r}|(1-\ov{\mu}(\,\ov{Q}_\infty))-
  \max\mb{r}\;(m-\dimH\Bad_A(\epsilon))\;,
\end{equation}
where $\Q^{(M)}=\bigvee_{k=0}^{M-1}\ta^{-k}\Q$.
Since Equation \eqref{staticentropybound} is clear if $\ov{\mu}
(\,\ov{Q}_\infty)=1$, we may assume that $\ov{\mu}
(\,\overline{Q}_\infty) <1$, hence that $\wt\mu_{j}
(\,\overline{Q}_\infty) < 1$ for all large enough $j\geq 1$.  Now, we
fix such a $j\geq 1$.

Take $\rho>0$ small enough so that $\wt\mu_j (\,\overline{Q}_\infty) +
\rho <1$ and let
\begin{equation}\label{eq:defibeta}
  \beta=\wt\mu_j(\,\overline{Q}_\infty) + \rho \;.
\end{equation}
Then for all large enough $i\in\bZ_{\geq 1}$, since
$\phi_A(S_{i,j})\subset \pi^{-1}(x_A)$ and $Q_\infty=
\pi^{-1}(Q_\infty^*)$ by Property (i) of $\Q$, we have
\[\begin{split}
\beta=\wt\mu_j (\,\overline{Q}_\infty) + \rho > \wt\nu_{i,j}(Q_\infty)
&
=\frac{1}{k_i \card\; S_{i,j}}\sum_{k=0}^{k_i-1}
\sum_{\bm{\theta}\in S_{i,j}} \Delta_{\ta^k \phi_{A}(\bm{\theta})}(Q_\infty)
\\&
=\frac{1}{k_i}\sum_{k=0}^{k_i-1}\Delta_{\ta^k x_A}(Q_{\infty}^{*})\;.
\end{split}
\]
Thus, for every $\bm{\theta}\in \TT^m$, since  $\ta^k
\phi_{A}(\bm{\theta}) \in Q_\infty$ implies that $\ta^k x_A\in
Q^*_\infty$ by Property (i) of $\Q$, we have
\begin{equation}\label{Icount}
  \card\{k\in\{0,\dots,k_i-1\}:
  \ta^k \phi_{A}(\bm{\theta})\in Q_\infty\}<\beta\, k_i\;.
\end{equation}

Let us prove the following counting lemma inspired by
\cite[Lem.~4.5]{EinLinMicVen12} and \cite[Lem.~2.4]{LimSaxSha18},
where $\ell_0$ is given by Property (ii) of $\Q$.

\blemm \label{lem:CovCount} There exists a constant $C>0$ depending
only on $\mb{r}$ and $\ell_0$ such that for all $A\in \cM_{m,n}(K_v)$,
$\bm{\theta}\in \bT^m$ and $T\in\bZ_{\geq 0}$, defining
$y=\phi_{A}(\bm{\theta})$, $I=\{k\in\bZ_{\geq 0}: \ta^k y \in
Q_{\infty}\}$, and
\[
E_{y,T}=\{z\in U^+ y : \forall\; k\in\{0,\dots,T\} - I, \;\;
d_{U^+\ta^k y,\|\cdot\|}(\ta^k y, \ta^k z)< q_v^{-\ell_0 \max\mb{r}}\;\}\;,
\]
the set $E_{y,T}$ can be covered by
$C\,q_v^{\,|\mb{r}|\card(I\cap\{0,\dots,T\})}$ closed balls of radius
$q_v^{-(\ell_0+T)}$ for the distance $d_{U^+y,\,\mb{r}}$.
\elemm

\dem As in the proof of \cite[Lemma 2.4]{LimSaxSha18}, we proceed by
induction on $T$. 

By the compactness of $\TT^m$, there exists a constant $C\in\bZ_{\geq
  1}$ depending only on $\mb{r}$ and $\ell_0$ such that the metric
space $(\TT^m,d_{\TT^m,\mb{r}})$ can be covered by $C$ closed balls of
radius $q_v^{-\ell_0}$.  Since $\phi_A:\TT^m\ra U^+y$ is an isometry
for the distances $d_{\TT^m,\mb{r}}$ and $d_{U^+y,\,\mb{r}}$, the
orbit $U^{+}y$ can be covered by $C$ closed balls for
$d_{U^+y,\,\mb{r}}$ of radius $q_v^{-\ell_0}$. Thus the lemma holds
for $T=0$. Let $N_{T}=C\,q_v^{\,|\mb{r}|\card(I\cap\{0,\dots,T\})}$.

Assume by induction that $E_{y,T-1}$ can be covered by $N_{T-1}$ balls
for $d_{U^+y,\,\mb{r}}$ of radius $q_v^{-(\ell_0+T-1)}$. Note that for
every $k\in\bZ$, since $\pi_v^k\cO_v/(\pi_v^{k+1}\cO_v)$ has order
$q_v$, every closed ball in $K_v$ of radius $q_v^{-k}$ is the disjoint
union of $q_v$ closed ball of radius $q_v^{-k-1}$. Hence every closed
ball for $d_{U^+y,\,\mb{r}}$ of radius $q_v^{-(\ell_0 +T-1)}$ in
$U^+y$ can be covered by $q_v^{\,|\mb{r}|}$ closed balls for
$d_{U^+y,\,\mb{r}}$ of radius $q_v^{-(\ell_0 +T)}$.  Therefore, if
$T\in I$, then $E_{y,T}=E_{y,T-1}$ can be covered by $N_T=
q_v^{|\mb{r}|}N_{T-1}$ closed balls for $d_{U^+y,\,\mb{r}}$ of radius
$q_v^{-(\ell_0+T)}$.

Suppose conversely that $T\notin I$, so that in particular $N_T=
N_{T-1}$.  Denote the above covering of $E_{y,T-1}$ by $\{B_i :
i=1,\dots,N_{T-1}\}$. Since we have $E_{y,T}\subset E_{y,T-1}$, the
set $\{E_{y,T}\cap B_i : i=1,\dots,N_{T-1}\}$ is a covering of
$E_{y,T}$.

\medskip\noindent{\bf Claim. } For all $i=1,\dots,N_{T-1}$
and $z_1,z_2 \in E_{y,T}\cap B_i$, we have
$d_{U^+y,\,\mb{r}}(z_1,z_2)\leq q_v^{-(\ell_0+T)}$.

\medskip
\dem
Since $T\notin I$, we have $d_{U^+\ta^T y,\|\cdot\|}(\ta^T y,\ta^T z_j)
< q_v^{-\ell_0 \max\mb{r}}$ for each $j=1,2$. Thus we have
$d_{U^+\ta^T y,\|\cdot\|} (\ta^T z_1,\ta^T z_2)<q_v^{-\ell_0\max\mb{r}}$ by
the ultrametric inequality property of $\|\cdot\|$.  Note that since
$z_1,z_2 \in U^{+} y= U^{+} y_{A,\bm{\theta}}$, there exist
$\bm{\theta}_1 = (\theta_{1,1},\dots, \theta_{1,m})$ and
$\bm{\theta}_2 = (\theta_{2,1},\dots,\theta_{2,m})$ in $\TT^m$ such
that (denoting in the same way lifts of $\bm{\theta}_1$ and
$\bm{\theta}_2$ to $K_v^{\,m}$) we have $z_1=y_{A,\bm{\theta}_1}$ and
$z_2= y_{A,\bm{\theta}_2}$. With $|\idist{\;}|$ the map defined after
Equation \eqref{eq:RcapO} and by Lemma \ref{lem:isomphiA}, it follows
that we have
\begin{align*}
\max_{1\leq i\leq m} q_v^{r_i T}\,|\idist{\theta_{1,i}-\theta_{2,i}}|=
d_{\TT^m} (\ta_-^T \bm{\theta}_1,\ta_-^T \bm{\theta}_2)&=
d_{U^+\ta^T y,\|\cdot\|} (\ta^T y_{A,\bm{\theta}_1},\ta^T y_{A,\bm{\theta}_2})\\ & =
d_{U^+\ta^Ty,\|\cdot\|}(\ta^T z_1,\ta^T z_2)<q_v^{-\ell_0\max\mb{r}}\;.
\end{align*}
Hence, we have
\[
d_{U^+y,\,\mb{r}}(z_1,z_2)=d_{\TT^m,\mb{r}}(\bm{\theta}_1, \bm{\theta}_2)=
\max_{1\leq i\leq m}|\idist{\theta_{1,i}-\theta_{2,i}}|^{\frac{1}{r_i}}
< q_v^{-(\ell_0+T)}\;,
\]
which concludes the claim. \cqfd

\medskip
By the above claim, the intersection $E_{y,T}\cap B_i$ is contained in a
single ball for $d_{U^+y,\,\mb{r}}$ of radius $q_v^{-(\ell_0+T)}$ for
each $i=1,\dots,N_{T-1}$.  Thus $E_{y,T}$ can be covered by $N_T=
N_{T-1}$ balls for $d_{U^+y,\,\mb{r}}$ of radius $q_v^{-(\ell_0+T)}$.
\cqfd

\medskip
Recall that as constructed in the proof of Proposition
\ref{prop:Apartconst}, there exist a Borel-measura\-ble partition
$\P=\{P_1,\dots, P_N,P_\infty\}$ of $\cY$ with $N+1$ elements, and a
countably generated Borel-measurable $\sigma$-algebra $\P^{U^+}$ of
subsets of $\cY$, with $[y]_{\P^{U^+}}=[y]_\P\cap B^{U^+}_ry$ for
every $y\in \cY(r)$ by Equation \eqref{eq:atomPUplus}, such that we
have $\A^{U^+}= (\P^{U^+})_0^{\infty}$ as defined in Equation
\eqref{eq:deficalAUplus}.  We now consider the sequence of
$\sigma$-algebras $\big((\P^{U^+})_{\ell}^\infty\big)_{\ell\geq 1}$,
which is a decreasing sequence of $\sigma$-algebras whose intersection
is $\A^{U^+}_\infty$ by Equation \eqref{eq:tailalg}.  Note that for
every $\ell\in\bZ_{\geq 1}$, the $\sigma$-algebra
$(\P^{U^+})_\ell^\infty$ is countably generated.

Let $Q$ be an atom of the finite partition $\Q^{(k_i)}=
\bigvee_{k=0}^{k_i-1}\ta^{-k}\Q$ of $\cY$. First assume that
$\phi_A(S_{i,j}) \cap Q$ is nonempty. Fix any $y\in\phi_A(S_{i,j})\cap
Q$ and for every $k=0,\dots,k_i-1$, let $Q_k\in \Q$ be such that $Q=
\bigcap_{k=0}^{k_i-1} \ta^{-k}Q_k$. For all $k\in\{0,\dots,k_i-1\}-
\{k\in\bZ_{\geq 0}: \ta^k y \in Q_{\infty}\}$ and $z\in\phi_A(S_{i,j})
\cap Q$, we have $z\in U^+y$ and $\ta^k y, \ta^k z\in \ta^k(U^+ y\cap
\ta^{-k}Q_k)=U^+\ta^k y\cap Q_k$ since $\ta$ normalizes $U^+$.  By
Property (ii) of $\Q$, we hence have
\[
d_{U^+\ta^k y,\|\cdot\|}(\ta^k y, \ta^k z)\leq
\diam_{d_{U^+\ta^k y,\|\cdot\|}}\; Q_k < q_v^{-\ell_0\max\mb{r}}\;.
\]
Therefore the intersection $\phi_A(S_{i,j}) \cap Q$ is contained in
the set $E_{y,k_i-1}$ defined in Lemma \ref{lem:CovCount}.  It follows
from Lemma \ref{lem:CovCount} and Equation \eqref{Icount} that
$\phi_A(S_{i,j}) \cap Q$ can be covered by $C\,q_v^{\,|\mb{r}|\beta
  k_i}$ closed balls for $d_{U^+y,\,\mb{r}}$ of radius
$q_v^{-(\ell_0+k_i-1)}=q_v^{-\ell_0+1}q_v^{-k_i}$, where $C$ depends
only on $\mb{r}$ and $\ell_0$. Since $S_{i,j}$ is
$q_v^{-k_i}$-separated (hence $q_v^{-\ell_0+1}q_v^{-k_i}$-separated
since $\ell_0\geq 1$) with respect to $d_{\TT^m,\mb{r}}$, and since
$\phi_A:(\TT^n,d_{\TT^m,\mb{r}})\ra (U^+y, d_{U^+y,\,\mb{r}})$ is an
isometry, we have
\[
\card(\phi_A(S_{i,j})\cap Q)\leq C\,q_v^{\,|\mb{r}|\beta k_i}\;.
\]
This formula is also satisfied if $\phi_A(S_{i,j}) \cap Q = \emptyset$.

\medskip\noindent{\bf Claim. }  For every large enough $\ell \in
\bZ_{\geq 1}$, we have $H_{\nu_{i,j}}(\sigma(\Q^{(k_i)})
\,|\,(\P^{U^+})_\ell^\infty)= H_{\nu_{i,j}}(\Q^{(k_i)})$.

\medskip
\dem Observe that using Proposition \ref{prop:Apartconst}
\eqref{Apartconst:prop3}, for every $\ell \in \bZ_{\geq 1}$, we have
\[
[y_{A,0}]_{(\P^{U^+})_\ell^\infty} =\ta^\ell [\ta^{-\ell}y_{A,0}]_{(\P^{U^+})_0^\infty}
=\ta^\ell [\ta^{-\ell}y_{A,0}]_{\A^{U^+}}
\supset \ta^{\ell} B_{\delta_r}^{U^+} \ta^{-\ell} y_{A,0}\;.
\]
As in the proof of Lemma \ref{lem:contractboulUplus}, we have
\[
\ta^{\ell} B_{\delta_r}^{U^+} \ta^{-\ell} \supset
B_{\ln (1+(e^{\delta_r} -1)q_v^{\ell \min \mb{r}})}^{U^+}\;,
\]
which implies that
\[
[y_{A,0}]_{(\P^{U^+})_\ell^\infty}  \supset
B_{\ln (1+(e^{\delta_r} -1)q_v^{\ell \min \mb{r}})}^{U^+} y_{A,0}\;.
\]
It follows that $[y_{A,0}]_{(\P^{U^+})_\ell^\infty} = U^+ y_{A,0}$ for
all large enough $\ell \in \bZ_{\geq 1}$ since the orbit $U^+ y_{A,0}$
is compact. Hence we have $H_{\nu_{i,j}}(\sigma(\Q^{(k_i)})\,|\,
(\P^{U^+})_\ell^\infty)= H_{\nu_{i,j}}(\Q^{(k_i)})$ for all large
enough $\ell \in\bZ_{\geq 1}$ since $\nu_{i,j}$ is supported on a
single atom of $(\P^{U^+})_\ell^\infty$. This proves the claim.
\cqfd

\medskip
Since the map $\Psi=-\log_{q_v}$ is nonincreasing, for every large
enough $\ell \in \bZ_{\geq 1}$, it hence follows that
\begin{align*}
  &H_{\nu_{i,j}}(\sigma(\Q^{(k_i)})\,|\,(\P^{U^+})_\ell^\infty)=
  H_{\nu_{i,j}}(\Q^{(k_i)})=\sum_{Q\in\Q^{(k_i)}}
  \nu_{i,j}(Q) \Psi\big(\nu_{i,j}(Q)\big)\\
  =\;&\sum_{Q\in\Q^{(k_i)}}\nu_{i,j}(Q) \Psi
  \Big(\frac{\card(\phi_A(S_{i,j})\cap Q)}{\card \;S_{i,j}}\Big)\\
  \geq\;& \Psi\Big(\frac{C\,q_v^{\,|\mb{r}|\beta k_i}}{\card \;S_{i,j}}\Big)
  \sum_{Q\in\Q^{(k_i)}}\nu_{i,j}(Q)
  =\log_{q_v}(\card\; S_{i,j}) -|\mb{r}|\,\beta\, k_i -\log_{q_v}C\;.
\end{align*}
By taking $\ell\to+\infty$, it follows from the continuity of entropy
that
\begin{equation}\label{entropybound1}
H_{\nu_{i,j}}(\sigma(\Q^{(k_i)})\,|\,\A^{U^+}_\infty)
\geq \log_{q_v}(\card\; S_{i,j}) -|\mb{r}|\,\beta\, k_i -\log_{q_v}C\;.
\end{equation}

Since $\A^{U^+}_\infty$ is strictly $\ta$-invariant, by the
subadditivity and concavity properties of the entropy as in the proof
of \cite[Eq.~(2.9)]{LimSaxSha18}, for every $M \in \bZ_{\geq 1}$, we
have
\begin{equation}\label{entropybound}
  \frac{1}{M}H_{\wt\nu_{i,j}}(\sigma(\Q^{(M)})\,|\,\A^{U^+}_\infty) \geq
  \frac{1}{k_i}H_{\nu_{i,j}}(\sigma(\Q^{(k_i)})\,|\,\A^{U^+}_\infty)-
  \frac{2\,M\,\log_{q_v} (\card\; \Q)}{k_i}\;.
\end{equation}
Therefore, since $\nu_{i,j}(\infty_\cY)=0$, it follows from Equations
\eqref{entropybound} and \eqref{entropybound1} that
\[
\begin{split}
  \frac{1}{M}&H_{\wt\nu_{i,j}}\big(\;\overline{\sigma(\Q^{(M)})}\,|\,
  \overline{\A^{U^+}_\infty} \big)=
  \frac{1}{M}H_{\wt\nu_{i,j}}\big(\sigma(\Q^{(M)})\,|\,\A^{U^+}_\infty\big)\\
  &\geq \frac{1}{k_i}\left(\log_{q_v}(\card\; S_{i,j}) -|\mb{r}|\beta k_i -
  \log_{q_v}C - 2\,M\,\log_{q_v} (\card\; \Q)\right).
\end{split}
\]
Now we can take $i\to\infty$ since the atoms $Q$ of the partition
$\ov{\Q}$ and hence of the partition $\overline{\Q^{(M)}}$, satisfy
$\wt\mu_j (\partial Q)=0$ by the property (iii) of $\Q$. Also, the
constants $C$ and $\card\; \Q$ are independent of $k_i$. Thus it
follows from Equation \eqref{ineqDim} that
\[
\frac{1}{M}H_{\wt\mu_j}\big(\;\overline{\sigma(\Q^{(M)})}\,|\,
  \overline{\A^{U^+}_\infty} \big) \geq |\mb{r}|(1-\beta) -
\max\mb{r}\; (m + \frac{1}{j} - \dimH \Bad_A(\epsilon))\;.
\]
By taking $\rho\to 0$ in Equation \eqref{eq:defibeta}, we have
\[
\frac{1}{M}H_{\wt\mu_j}\big(\;\overline{\sigma(\Q^{(M)})}\,|\,
  \overline{\A^{U^+}_\infty} \big) \geq
|\mb{r}|(1-\wt\mu_j(\,\ov{Q}_\infty))
- \max\mb{r}\; (m + \frac{1}{j} - \dimH \Bad_A(\epsilon))\;.
\]
Hence, it follows by taking $j\to \infty$ and by using the property
(iii) of $\Q$ that
\[
\frac{1}{M}H_{\ov{\mu}}\big(\;\overline{\sigma(\Q^{(M)})}\,|\,
  \overline{\A^{U^+}_\infty} \big) \geq
|\mb{r}|(1-\ov{\mu}(\,\ov{Q}_\infty))
- \max\mb{r}\; (m - \dimH \Bad_A(\epsilon))\;,
\]
which proves Equation \eqref{staticentropybound}.

Hence, by taking $M\to\infty$, we have
\[
h_{\ov\mu}(\ta^{-1}|\;\overline{\A^{U^+}_\infty}\,)
=h_{\ov\mu}(\ta\,|\;\overline{\A^{U^+}_\infty}\,)\geq
|\mb{r}| (1- \ov{\mu} (\ov{Q}_\infty)) -
\max\mb{r}\;(m -\dimH \mb{Bad}_A(\epsilon))\;,
\]
provided that we have a partition $\Q$ satisfying the above
requirements (i), (ii) and (iii).  After taking a sufficiently small
neighborhood of infinity $Q_\infty^*$ in $\cX$, so that if
$Q_\infty=\pi^{-1}(Q_\infty^*)$, then $\ov{\mu} (\ov{Q}_\infty)$ is
sufficiently close to $\ov\mu(\infty_\cY) = \eta_A$, we can
indeed construct a finite Borel-measurable partition $\cQ$ of $\cY$
satisfying Properties (i), (ii) and (iii), by following the procedure
in \cite[Proof of Theorem 4.2, Claim 2]{LimSaxSha18}.  This proves
Assertion \eqref{entropy}.  \cqfd

\subsection{Effective upper bound on $\dimH \mb{Bad}_{A}(\epsilon)$}

For every $\ell \in\bZ_{\leq 1}$, with $\lambda_1$ the shortest length
function of a nonzero vector of an $R_v$-lattice (see Subsection \ref
{subsec:geomnumb}), we define
\[
\cX^{\geq q_v^{\ell}}=\{x\in \cX : \lambda_1 (x) \geq q_v^{\ell} \}
\quad\text{and}\quad \cY^{\geq q_v^{\ell}} =\pi^{-1}(\cX^{\geq q_v^{\ell}})\;.
\]
Note that by Corollary \ref{coro:minko}, we have $\lambda_1 (x)\leq
q_v$ for all $x\in\cX$, thus $\cX=\bigcup_{\ell=-\infty}^1\cX^{\geq
  q_v^{\ell}}$.  By Mahler's compactness criterion (see for instance
\cite[Theo.~1.1]{KleShiTom17}), the subsets $\cX^{\geq q_v^{\ell}}$
and $\cY^{\geq q_v^{\ell}}$ are compact.

\blemm\label{lem:condmeas} Let $\mu'$ be an $\ta$-invariant Borel
probability measure on $\cY$ and let $\A$ be a countably generated
sub-$\sigma$-algebra of the Borel $\sigma$-algebra of $\cY$ which is
$\ta^{-1}$-descending and $U^+$-subordinate modulo $\mu'$.  For all
$r'\geq \delta'>0$, $\epsilon\in\;]0,1]$ and $\ell\in\bZ_{\leq 0}$, let
$j_1,j_2$ be integers satisfying
$$j_1>
\frac{d-(d-1)\ell}{\min\mb{r}} -\log_{q_v}\delta'
\quad \text{and}\quad j_2 >
\frac{d-(d-1)\ell}{\min\mb{s}} -\frac{n}{d}\log_{q_v}\epsilon\:.
$$
If $y\in\cY^{\geq q_v^{\ell}}$ satisfies
$B_{\delta'}^{U^+,\mb{r}} \ta^{-j_1}y \subset
[\ta^{-j_1}y]_\A \subset B_{r'}^{U^+,\mb{r}} \ta^{-j_1}y$,
then we have
\[
\tau_y^{\ta^{j_1} \A}(\ta^{-j_2}\cL_\epsilon)\leq 1-\frac{1}{q_v^m}
\left(q_v^{-(j_1+j_2)}(r')^{-1}\epsilon^{\frac{m}{d}}\right)^{|\mb{r}|}\;.
\]
\elemm

\dem Let $x=\pi(y)$, which belongs to $\cX^{\geq q_v^{\ell}}$. Since
$x$ is a unimodular $R_v$-lattice, by Minkowski's theorem
\ref{theo:Mink2}, we hence have
\[
q_v^{(d-1)\ell}\lambda_d (x)\leq (\lambda_1 (x))^{d-1}\lambda_d (x)
\leq \lambda_1 (x)\lambda_2 (x)\cdots \lambda_d (x) \leq q_v^d\;,
\]
therefore $\lambda_d (x)\leq q_v^{d-(d-1)\ell}$. There are linearly
independent vectors $v_1,\dots,v_d$ in the $R_v$-lattice $x$ such that
$\|v_i\|\leq q_v^{d-(d-1)\ell}$. Let $\Delta$ be the parallelepiped in
$K_v^{\,d}$ generated by $v_1,\dots,v_d$, that is,
\[
\Delta=\{t_1v_1+\cdots+t_dv_d\in K_v^{\,d} :\;
\forall\; i=1,\dots,d, \;\;|\,t_i\,| \leq 1\}\;.
\]
We identify $K_v^{\,d}$ with $K_v^{\,m}\times K_v^{\,n}$. Then for every
$\mb{b}=(\mb{b}^-,\mb{b}^+)\in \Delta$ with $\mb{b}^-\in K_v^{\,m}$
and $\mb{b}^+\in K_v^{\,n}$, we have $\|\,\mb{b}\,\|\leq
q_v^{d-(d-1)\ell}$, hence $\|\,\mb{b}^-\,\|_{\mb{r}}\leq
q_v^{\frac{d-(d-1)\ell}{\min\mb{r}}}$ and $\|\,\mb{b}^+\,\|_{\mb{s}}\leq
q_v^{\frac{d-(d-1)\ell}{\min\mb{s}}}$ since $\ell\leq 0$. Note that
the fiber $\pi^{-1}(x)$ can be parametrized as follows: Fixing $g\in
G_0$ with $x=g\Gamma_0$, since $\Delta$ is a fondamental domain for
the action of $R_v^{\,d}$ on $K_v^{\,d}$, we have
\[
\pi^{-1}(x)=\{w(\mb{b})g\Gamma : \mb{b}\in \Delta\}, \text{ where }
w(\mb{b})= \begin{pmatrix} I_d & \mb{b} \\  0 & 1 \end{pmatrix}\;.
\]
In particular, there exists $\mb{b}_0= (\mb{b}_0^-, \mb{b}_0^+) \in
\Delta$ such that $y=w(\mb{b}_0)g\Gamma$.

With a slightly simplified notation, let $V_y$ be the $U^+$-shape of
the atom $[y]_{\ta^{j_1}\A}$ (see Equation \eqref{eq:defishape}), so
that we have $V_y y= [y]_{\ta^{j_1}\A}$. Let $\Xi=\{\bm{\theta}\in
K_v^{\,m}:w(\bm{\theta},0)\in V_y\}$ be the Borel set corresponding to
$V_y$ by the canonical bijection $\bm{\theta}\mapsto w(\bm{\theta},0)$
(see above Equation \eqref{eq:dUplus}) between $K_v^{\,m}$ and
$U^+$. Note that $0\in\Xi$ as $I_{d+1}\in V_y$. Since $\ta_-^{j_1}$
expands the $\mb{r}$-quasinorm on $K_v^{\,m}$ with ratio exactly
$q_v^{j_1}$ (see Equation \eqref{eq:dilatweight}), and by the
assumption on $y$ in the statement of Lemma \ref{lem:condmeas}, we
have $B_{q_v^{j_1}\delta'}^{U^+,\mb{r}} y \subset [y]_{\ta^{j_1}
  \A}\subset B_{q_v^{j_1}r'}^{U^+,\mb{r}} y$. Hence by the definition
of the $U^+$-shapes of atoms in Equation \eqref{eq:defishape}, we have
\begin{equation}\label{eq:controlXi}
B_{q_v^{j_1}\delta'}^{K_v^{\,m},\mb{r}} \subset \Xi \subset
\overline{B}_{q_v^{j_1}r'}^{K_v^{\,m},\mb{r}}\;.
\end{equation}

The atom $[y]_{\ta^{j_1}\A}$ can be parametrized by
\[
  [y]_{\ta^{j_1}\A} = \big\{w(\mb{b})g\Gamma:\;\exists\;\mb{b}^-
  \in \mb{b}_0^- +\Xi,\;\; \mb{b}=(\mb{b}^-,\mb{b}_0^+) \big\}\;,
\]
and $\tau_y^{\ta^{j_1}\A}$ is the pushforward measure of the
normalized Haar measure on the Borel set (with positive measure)
$\mb{b}_0^- +\Xi$ of $K_v^{\,m}$.

Let us consider the sets
\[
\Theta^- =\{\mb{b}^- \in K_v^{\,m} :
\|\,\mb{b}^-\,\|_\mb{r}< q_v^{-j_2}\epsilon^{\frac{m}{d}}\}\;
\text{ and }\; \Theta^+ =\{\mb{b}^+ \in K_v^{\,n} :
\|\,\mb{b}^+\,\|_\mb{s}< q_v^{j_2}\epsilon^{\frac{n}{d}}\}\;.
\]
If $\mb{b}=(\mb{b}^-,\mb{b}^+)\in \Theta^- \times \Theta^+$, then
$\|\ta_{-}^{j_2}\mb{b}^-\|_{\mb{r}}< \epsilon^{\frac{m}{d}}$ and
$\|\ta_{+}^{j_2}\mb{b}^+\|_{\mb{s}}< \epsilon^{\frac{n}{d}}$ by
Equation \eqref{eq:dilatweight}. By the definition of $\cL_\epsilon$
in Equation \eqref{eq:defiLsubeps}, and since the grid
$\ta^{j_2}gR_v^{\,m}+(\ta_{-}^{j_2}\mb{b}^-, \ta_{+}^{j_2}\mb{b}^+)$
contains the vector $(\ta_{-}^{j_2}\mb{b}^-, \ta_{+}^{j_2}\mb{b}^+)$,
we have
\[
\ta^{j_2}w(\mb{b})g\Gamma=w(\ta_{-}^{j_2}\mb{b}^-, \ta_{+}^{j_2}\mb{b}^+)
\ta^{j_2}g\Gamma \notin \cL_\epsilon\;.
\]
Hence we have $w(\mb{b})g\Gamma \notin \ta^{-j_2}\cL_\epsilon$, so that
\begin{equation}\label{eq:poidsomgabgga}
  [y]_{\ta^{j_1}\A}-\ta^{-j_2}\cL_\epsilon \supset w\big(\,\big((\mb{b}_0^-+\Xi)
  \times \{\mb{b}_0^+\}\big) \cap (\Theta^- \times \Theta^+)\big)g\Ga\;.
\end{equation}

\noindent{\bf Claim. } We have the inclusion $\Theta^- \times
\{\mb{b}_0^+\} \subset \big((\mb{b}_0^- +\Xi ) \times
\{\mb{b}_0^+\}\big) \cap (\Theta^- \times \Theta^+)$.

\medskip
\dem We only have to prove that $\mb{b}_0^+ \in \Theta^+$ and that
$\Theta^-\subset \mb{b}_0^- + \Xi$. Since $(\mb{b}_0^-, \mb{b}_0^+) \in
\Delta$, we have $\|\,\mb{b}_0^+\,\|_{\mb{s}}\leq
q_v^{\frac{d-(d-1)\ell}{\min\mb{s}}}$, hence the former assertion follows
from the assumption that $j_2 > \frac{d-(d-1)\ell}{\min\mb{s}}
-\frac{n}{d}\log_{q_v}\epsilon$.

In order to prove the latter assertion, let us fix $\mb{b}^- \in
\Theta^-$. Recall that the $\mb{r}$-quasinorm $\|\cdot\|_\mb{r}$
satisfies the ultrametric inequality property, see Equation
\eqref{eq:ultraineqpsudonorm}.  Hence, it follows from the assumptions
$j_2 > \frac{d-(d-1)\ell}{\min\mb{s}} -\frac{n}{d}
\log_{q_v}\epsilon$ and $j_1> \frac{d-(d-1)\ell}{\min\mb{r}}
-\log_{q_v}\delta'$, since $\epsilon\leq 1$, that
\[
\begin{split}
  \|\,\mb{b}^- -\mb{b}_0^-\, \|_\mb{r}&
  \leq \max \{\|\,\mb{b}^-\,\|_\mb{r},\|\,\mb{b}_0^-\,\|_\mb{r}\} 
  \leq \max \big\{q_v^{-j_2}\epsilon^{\frac{m}{d}},
  q_v^{\frac{d-(d-1)\ell}{\min\mb{r}}}\big\}\\
  &\leq \max\big\{q_v^{-\frac{d-(d-1)\ell}{\min\mb{s}}} \,\epsilon\,,
  q_v^{\frac{d-(d-1)\ell}{\min\mb{r}}}\big\}
=q_v^{\frac{d-(d-1)\ell}{\min\mb{r}}} < q_v^{\,j_1}\delta'\;.
\end{split}
\]
Hence by the left inclusion in Equation \eqref{eq:controlXi}, we have
$\mb{b}^-\in\mb{b}_0^- + B_{q_v^{\,j_1}\delta'}^{K_v^{\,m},\mb{r}}
\subset \mb{b}_0^- + \Xi$, which concludes the latter assertion. \cqfd

\medskip
Now by Equation \eqref{eq:poidsomgabgga} and by the above claim, since
the canonical isomorphism $w\mapsto \wh w$ (see the line above
Equation \eqref{eq:dUplus1}) from $K_v^m$ to $U^+$ maps the Haar
measure $\vol_v^m$ of $K_v^m$ to the Haar measure $m_{U^+}$ of $U^+$,
we have
\begin{equation}\label{eq:correc13june1}
  1-\tau_y^{\ta^{\,j_1}\A}(\ta^{-j_2}\cL_\epsilon)=
  \tau_y^{\ta^{\,j_1}\A}([y]_{\ta^{j_1}\A}- \ta^{-j_2}\cL_\epsilon) 
\geq \frac{\vol_v^{\,m}(\Theta^-)}{\vol_v^{\,m}(\mb{b}_0^- +\Xi)}\;.
\end{equation}
On one hand, an open ball in $K_v$ of radius $\eta>0$ having Haar
measure $q_v^{\lceil\log_{q_v}\eta\rceil-1}$, we have
\begin{align}
  \nonumber\vol_v^m (\Theta^- ) &
  = \vol_v^m\big(B^{K_v^m,\mb{r}}_{q_v^{-j_2}\epsilon^{\frac{m}{d}}}\big)
  = \prod_{i=1}^m \vol_v \big(B^{K_v}_{(q_v^{-j_2}\epsilon^{\frac{m}{d}})^{r_i}}\big)\\
&\geq \prod_{i=1}^m \frac{1}{q_v} (q_v^{-j_2}\epsilon^{\frac{m}{d}})^{r_i} 
  =\frac{1}{q_v^m}(q_v^{-j_2}\epsilon^{\frac{m}{d}})^{| \mathbf{r} |}\;.
  \label{eq:correc13june2}
\end{align}
On the other hand, by the right inclusion in Equation
\eqref{eq:controlXi}, we have
\begin{align}
 \nonumber\vol_v^m(\mb{b}_0^- +\Xi) &= \vol_v^m(\Xi) 
\leq  \vol_v^m\big(\,\overline{B}^{K_{v, \mb{r}}^m}_{q_v^{j_1}r'}\big) \\
&= \prod_{i=1}^m  \vol_v\big(\,\overline{B}^{K_v}_{(q_v^{j_1} r')^{r_i}}\big)
\leq \prod_{i=1}^m (q_v^{j_1} r')^{r_i} = (q_v^{j_1} r')^{|\mb{r}|}\;.
\label{eq:correc13june3}
\end{align}
Lemma \ref{lem:condmeas} follows by combining the inequalities
\eqref{eq:correc13june1}, \eqref{eq:correc13june2} and
\eqref{eq:correc13june3}.
\cqfd

\bigskip\noindent{\bf Proof of Theorem \ref{theo:introEff1}. } We fix
a matrix $A\in \cM_{m,n}(K_v)$ which is not $(\mb{r},\mb{s})$-singular
on average, or equivalently by Proposition
\ref{prop:critdynescapeonaver} and the definition of $\eta_A$ just
before Lemma \ref{prop:measconst}, we assume that $\eta_A<1$. We also fix
$\epsilon\in\;]0,1]$.

By Proposition \ref{prop:measconst}, there exist an $\ta$-invariant
Borel probability measure $\ov{\mu}$ on $\ov{\cY}$ (depending on
$\epsilon$) and an $\ta$-invariant Borel probability measure $\mu$ on
$\cY$ (unique since $\eta_A<1$) such that
\[
\Supp{\ov{\mu}} \subset \cL_\epsilon \cup \{\infty_{\cY}\},\;\; 
\ov{\pi}_* \ov{\mu}=\mu_A,\;\;  \text{and}\;\;\; 
\ov{\mu}=(1-\eta_A)\mu + \eta_A \Delta_{\infty_\cY}\;.
\]
Take a compact subset $K_0$ of $\cX$ such that $\mu_A(K_0)>0.99\,
\mu_A (\cX)= 0.99\,(1-\eta_A)$. Write $K=\pi^{-1}(K_0)$ and choose
$r\in\;]0,\ln 2[$ such that $K \subset \cY(r)$. Then $\mu(\cY(r))
    \geq\mu(K)>0.99$ since $\eta_A<1$.  Note that the choices of $K$
    and $r$ are independent of $\epsilon$ since the measure $\mu_A$
    depends only on $A$ (see Proposition \ref{prop:measconst} and
    Equation \eqref{eq:definmuA}).

For such an $r>0$, let $\A^{U^+}$ be the $\sigma$-algebra of subsets
of $\cY$ constructed in Proposition \ref{prop:Apartconst}.
Proposition \ref{prop:measconst} \eqref{entropy} gives the inequality
\[
h_{\ov{\mu}}(\ta^{-1}|\;\ov{{\A^{U^+}_\infty}}\,)\geq 
|\mb{r}| (1- \eta_A) - \max\mb{r}\;(m -\dimH \mb{Bad}_A(\epsilon))\;.
\]
By the linearity of entropy (and since the entropy of $\ta^{-1}$
vanishes on the fixed set $\{\infty_\cY\}$), we have
\begin{equation}\label{entropyEq1}
  h_{\mu}(\ta^{-1}|\;\A^{U^+}_\infty)=
  \frac{1}{1-\eta_A}h_{\ov{\mu}}(\ta^{-1}|\;\ov{{\A^{U^+}_\infty}}\,) \geq |\mb{r}|-
  \frac{\max\mb{r}}{1-\eta_A}\;(m -\dimH \mb{Bad}_A(\epsilon))\;.
\end{equation}

In order to use Lemma \ref{lem:algexiA} and Proposition
\ref{prop:effEL}, we need an ergodicity assumption on the measures
that appear in these statements. We will choose an appropriate ergodic
component of $\mu$. Using the notation of Proposition
\ref{prop:ergDec} with $X=\cY$ and $T=\ta^{-1}$, let us denote the
ergodic decomposition of $\mu$ (see for instance
\cite[Theo.~2.7]{EinLinWar22}) by
\[
\mu=\int_{y\in\cY} \mu_y^\E \;d\mu(y).
\]
Let $E=\{y\in\cY : \mu_{y}^{\E}(K)>0.9\}$. It follows from the
inequality $\mu(K)>0.99$ that
\[
0.99 < \int_\cY \mu_y^\E (K)\;d\mu(y) \leq
\mu(E) + 0.9\,\mu(\cY - E)=0.9+0.1\,\mu(E)\;,
\]
hence $\mu(E)>0.9$. By Proposition \ref{prop:ergDec} with $T=\ta^{-1}$
and $\A=\A^{U^+}_\infty$, and by Equation \eqref{entropyEq1}, we have
\[
\int_\cY h_{\mu_{y}^{\E}}(\ta^{-1}|\;\A^{U^+}_\infty)\; d\mu(y) =
h_{\mu}(\ta^{-1}|\;\A^{U^+}_\infty) \geq
|\mb{r}|- \frac{\max\mb{r}}{1-\eta_A}\; (m -\dimH \mb{Bad}_A(\epsilon))\;. 
\]
Since $h_{\mu_{y}^{\E}}(\ta^{-1}|\;\A^{U^+}_\infty)\leq |\mb{r}|$ 
for every $y\in\cY$ by Lemma \ref{lem:algexiA}, we have
\[
\int_{\cY- E} h_{\mu_{y}^{\E}}(\ta^{-1}|\;\A^{U^+}_\infty)
\;d\mu(y) \leq |\mb{r}|\;\mu(\cY- E)\;.
\]
Hence
\begin{align*}
\int_E h_{\mu_{y}^{\E}}(\ta^{-1}|\;\A^{U^+}_\infty) \;d\mu(y) 
&\geq |\mb{r}|\,\mu(E)-\frac{\max\mb{r}}{1-\eta_A}\;
(m -\dimH \mb{Bad}_A(\epsilon))\\
&\geq \mu(E)\Big(|\mb{r}| -\frac{\max\mb{r}}{0.9\,(1-\eta_A)}\;
(m -\dimH \mb{Bad}_A(\epsilon))\Big)\;.
\end{align*}
Therefore, there exists $z\in \cY$ such that $\mu_z^{\E}(K)>0.9$
and
\[
h_{\mu_z^{\E}}(\ta^{-1}|\;\A^{U^+}_\infty)\geq|\mb{r}|
-\frac{\max\mb{r}}{0.9\,(1-\eta_A)}\; (m -\dimH \mb{Bad}_A(\epsilon))\;.
\] 
We denote  $\lambda=\mu_z^{\E}$ for such a $z\in \cY$. Then $\lambda$ is an
$\ta$-invariant ergodic Borel probability measure on $\cY$ and
$\Supp{\lambda}\subset \Supp{\mu}\subset \cL_{\epsilon}$. By Lemma
\ref{lem:algexiA}, we have
\begin{equation}\label{entropyEq2}
  H_{\lambda}(\A^{U^+}|\;\ta \,\A^{U^+})\geq |\mb{r}|
  -\frac{\max\mb{r}}{0.9\,(1-\eta_A)}\; (m -\dimH \mb{Bad}_A(\epsilon))\;.
\end{equation}

We will apply Lemma \ref{lem:condmeas} with $\mu'=\lambda$ and
$\A=\ta^{-k}\A^{U^+}$ for some $k\geq 1$.  Take an integer $\ell \leq
0$ such that $K\subset \cY^{\geq q_v^{\ell}}$, which depends only on
$A$. Set
\begin{equation}\label{eq:defj1j2}
j_1 = \Big\lceil \frac{d-(d-1)\ell}{\min\mb{r}}
-\log_{q_v}\delta' \Big\rceil+1\; \text{ and }\; 
j_2 = \Big\lceil \frac{d-(d-1)\ell}{\min\mb{s}}
-\frac{n}{d}\log_{q_v}\epsilon  \Big\rceil+1\;,
\end{equation}
where $\delta'$ will be determined later on.

Let $k= \big\lceil \log_{q_v} \big(r^{\frac{1}{\max\mb{r}}}\,
\epsilon^{-\frac{m}{d}}\big) \big\rceil + j_2+1$ and $\A=
\ta^{-k}\A^{U^+}$.  By the properties of $\A^{U^+}$ given in
Proposition \ref{prop:Apartconst} and since $K\subset \cY(r)$, 
for every $y\in K$, we have
\[
B_{\delta_r}^{U^+}  y \subset [y]_{\A^{U^+}} \subset B_r^{U^+} y\;.
\]
It follows from Equations \eqref{eq:locLip} and
\eqref{eq:holderdistwithr} that since $\delta_r\leq r\leq \ln 2$, for
every $y\in K$, we have
\[
B_{(\delta_r/2)^{\frac{1}{\min{\mb{r}}}}}^{U^+,\mb{r}}  y \subset
B_{\delta_r/2}^{U^+,\|\;\|}  y
\subset [y]_{\A^{U^+}} \subset B_{r}^{U^+,\|\;\|}  y
\subset B_{r^{\frac{1}{\max{\mb{r}}}}}^{U^+,\mb{r}} y\;.
\]
Hence, by Equation \eqref{eq:contracboulUplusexact}, we have 
\[
B_{q_v^{-k}(\delta_r/2)^{\frac{1}{\min{\mb{r}}}}}^{U^+,\mb{r}}  \ta^{-k}y \subset
[\ta^{-k} y]_{\ta^{-k}\A^{U^+}} =[\ta^{-k} y]_{\A}
\subset B_{q_v^{-k}r^{\frac{1}{\max{\mb{r}}}}}^{U^+,\mb{r}} \ta^{-k}y\;.
\]
Thus for every $y \in \ta^k K$, we have
\begin{equation}\label{eq:controlfinatomBUplus}
B_{\delta'}^{U^+,\mb{r}}  y \subset [y]_{\A} \subset B_{r'}^{U^+,\mb{r}} y\;,
\end{equation}
where, by the definition of $k$, we take
\[
r'= q_v^{-j_2-1}\epsilon^{\frac{m}{d}} \quad \text{and}\quad \delta' =
q_v^{-1}\;r^{-\frac{1}{\max\mb{r}}}\;(\delta_r/2)^{\frac{1}{\min\mb{r}}}\;r'\;.
\]
Equation \eqref{eq:controlfinatomBUplus} implies that for every $y\in
\ta^{j_1+k}K$, we have
\begin{equation}\label{eq:lemmaassumption}
B_{\delta'}^{U^+,\mb{r}} \ta^{-j_1}y \subset [\ta^{-j_1}y]_{\A} 
\subset B_{r'}^{U^+,\mb{r}} \ta^{-j_1}y\;.
\end{equation}

Now, we will use Proposition \ref{prop:effEL} with $j=j_1$, $K'=\ta^k
K$ (which is $U^+$-saturated since so is $K$ and as $\ta$ normalizes
$U^+$), and $\epsilon=r'$ (which satisfies the assumption of
Proposition \ref{prop:effEL} by Equation
\eqref{eq:controlfinatomBUplus}). We claim that
\begin{equation}\label{eq:stabLepsBU}
B_{q_v^{-1}\epsilon^{\frac{m}{d}}}^{U^+,\mb{r}}\cL_\epsilon \subset \cL_\epsilon\;.
\end{equation}
Indeed, for all $y\in \cL_\epsilon$ and $\bm{\theta}\in K_v^{\,m}$
such that $\|\,\bm{\theta}\,\|_\mb{r}\leq q_v^{-1}
\epsilon^{\frac{m}{d}}$, for every vector $u=(u^-,u^+)$ in the grid
$w(\bm{\theta}, \mb{0})y$, we can write $u=v+ (\bm{\theta},\mb{0})$
for some $v=(v^-,v^+)$ in the grid $\wt\Lambda_y$ associated with $y$
(see Equation \eqref{eq:defiwtLambdasuby}).  Since $y\in
\cL_\varepsilon$, we have (see Equation \eqref{eq:defiLsubeps})
$\|\,v\, \|_{\mb{r}, \mb{s}}= \max\{\| \,v^-\,
\|_{\mb{r}}^{\;\frac{d}{m}}, \| \,v^+\, \|_{\mb{s}}^{\;\frac{d}{n}}\}
\geq \epsilon$. Since $u^+=v^+$, if $\| \,v^+\,
\|_{\mb{s}}^{\;\frac{d}{n}} \geq \epsilon$, then $w(\bm{\theta},
\mb{0})y \in\cL_\epsilon$. Otherwise $\| \,v^-
\,\|_{\mb{r}}^{\;\frac{d}{m}} \geq\epsilon$.  We then have
$\|\,\bm{\theta}\,\| _{\mb{r}}\leq q_v^{-1} \epsilon^{\frac{m}{d}}<
\epsilon^{\frac{m}{d}}\leq \| \,v^- \,\|_{\mb{r}}$.  It follows from
the equality case of the ultrametric inequality property of $\|\;
\,\|_{\mb{r}}$ that
\[
\| \,u^- \,\|_\mb{r}=\|\,\bm{\theta}+v^- \,\|_\mb{r}=
\max\big\{\|\,\bm{\theta}\,\|_\mb{r},\;\|\,v^- \,\|_\mb{r}\big\}=
\| \,v^- \,\|_\mb{r}\geq \epsilon^{\frac{m}{d}}\;.
\]
Hence $w(\bm{\theta}, \mb{0})y\in\cL_\epsilon$, which proves Equation
\eqref{eq:stabLepsBU}.

\medskip
By Proposition \ref{prop:Apartconst}, the $\sigma$-algebra $\A^{U^+}$
is $\ta^{-1}$-descending and $U^+$-subordinate modulo $\lambda$, and
so is $\A=\ta^{-k}\A^{U^+}$ since $\ta$ normalizes $U^+$.  Note that
$\Supp{\lambda} \subset \ta^{-j_2}\cL_\epsilon$ since $\lambda$ is
$\ta$-invariant.  By Equations \eqref{eq:contracboulUplusexact} and
\eqref{eq:stabLepsBU}, we have
\[
B_{r'}^{U^+,\mb{r}}\ta^{-j_2}\cL_\epsilon =
\ta^{-j_2}B_{q_v^{\,j_2}r'}^{U^+,\mb{r}}\cL_\epsilon
=\ta^{-j_2}B_{q_v^{-1}\epsilon^{\frac{m}{d}}}^{U^+,\mb{r}}\cL_\epsilon 
\subset \ta^{-j_2}\cL_\epsilon\;.
\]
Note that we have
\[
\tau_y^{\ta^{j_1}\A}(\cY-\ta^k K)=0
\]
for $\lambda$-almost every $y \in \ta^k K$, since then (see just above
Proposition \ref{prop:effEL}) the support $\Supp \tau_y^{\ta^{j_1}\A}$
is contained in $[y]_{\ta^{j_1}\A}$, which is contained in $U^+ y$, hence
in $\ta^k K$ since $\ta$ normalizes $U^+$ and $K=\pi^{-1}(K_0)$ is
$U^+$-saturated. Therefore, it follows from Proposition \ref{prop:effEL}
for the first line, from the fact that the integrated function is
nonpositive (hence its integral on a smaller domain is larger) for the
third line, that
\begin{align*}
H_{\lambda}(\A|\ta^{j_1}\A)&\leq j_1 |\mb{r}|+\int_{\cY}\log_{q_v}
\tau_y^{\ta^{j_1}\A}((\cY- \ta^k K)\cup B_{r'}^{U^+,\mb{r}}
\Supp{\lambda})\;d\lambda(y)\\
&\leq j_1 |\mb{r}|+\int_{\cY}\log_{q_v}\tau_y^{\ta^{j_1}\A}
((\cY- \ta^k K)\cup \ta^{-j_2}\cL_\epsilon)\;d\lambda(y)\\
&\leq j_1 |\mb{r}|+\int_{\ta^k K \cap \ta^{j_1+k}K \cap \cY^{\geq q_v^{\ell}}}
\log_{q_v}\tau_y^{\ta^{j_1}\A}((\cY- \ta^k K)\cup \ta^{-j_2}\cL_\epsilon)
\;d\lambda(y)\\&
= j_1 |\mb{r}|+\int_{\ta^k K \cap \ta^{j_1+k}K \cap \cY^{\geq q_v^{\ell}}}
\log_{q_v}\tau_y^{\ta^{j_1}\A}(\ta^{-j_2}\cL_\epsilon)\;d\lambda(y)\;.
\end{align*}

We now apply Lemma \ref{lem:condmeas} with as said above
$\mu'=\lambda$ and $\A=\ta^{-k}\A^{U^+}$, with $j_1$ and $j_2$ defined
in Equation \eqref{eq:defj1j2}, and with $y \in \ta^{j_1+k}K \cap
\cY^{\geq q_v^{\ell}}$ which satisfies the assumption of Lemma
\ref{lem:condmeas} by Equation \eqref{eq:lemmaassumption}. Thus
\[
\tau_y^{\ta^{j_1}\A} (\ta^{-j_2}\cL_\epsilon) \leq
1-\frac{1}{q_v^m}\left(q_v^{-(j_1+j_2)}{r'}^{-1}\epsilon^{\frac{m}{d}}\right)^{|\mb{r}|}
=1-q_v^{-(j_1-1)|\mb{r}|-m}\;.
\]
Hence 
\[
-\log_{q_v}\tau_y^{\ta^{j_1}\A} (\ta^{-j_2}\cL_\epsilon) \geq
-\log_{q_v}\big(1-q_v^{-(j_1-1)|\mb{r}|-m}\big) 
\geq \frac{q_v^{-(j_1-1)|\mb{r}|-m}}{\ln q_v}\;.
\]
Note that $\lambda(\ta^k K \cap \ta^{j_1+k}K \cap \cY^{\geq
  q_v^{\ell}}) \geq \frac{1}{2}$ since $\lambda$ is $\ta$-invariant,
$K\subset\cY^{\geq q_v^{\ell}}$ and $\lambda(K)>0.9$, so that the
three sets $\ta^k K$, $\ta^{j_1+k}K$ and $\cY^{\geq q_v^{\ell}}$ have
$\lambda$-measure $>0.9$, hence their pairwise intersections have
$\lambda$-measure $>2\times 0.9-1=0.8\,$, and their triple intersection
has $\lambda$-measure $>2\times 0.8-1=0.6\,$.  It follows from Equation
\eqref{eq:chainrule} and the invariance under $\ta$ of $\lambda$,
hence of the conditional entropy, that
\begin{align*}
  |\mb{r}|-H_{\lambda}(\A^{U^+}|\;\ta\,\A^{U^+})&=
  |\mb{r}|-\frac{1}{j_1}H_{\lambda}(\A^{U^+}|\;\ta^{j_1}\A^{U^+})
= |\mb{r}|- \frac{1}{j_1}H_{\lambda}(\A|\;\ta^{j_1}\A)\\
&\geq -\frac{1}{j_1}\int_{\ta^k K \cap \ta^{j_1+k}K \cap \cY^{\geq q_v^{\ell}}}
\log_{q_v} \tau_y^{\ta^{j_1}\A}(\ta^{-j_2}\cL_\epsilon)\;d\lambda(y)\\
 &\geq \frac{q_v^{|\mb{r}|-m}}{2\ln{q_v}}\frac{q_v^{-j_1|\mb{r}|}}{j_1}\;.
\end{align*}
Therefore, by Equation \eqref{entropyEq2}, we have
\[
\frac{\max\mb{r}}{0.9\,(1-\eta_A)}\;(m-\dimH\Bad_A(\epsilon)) \geq
|\mb{r}|-H_{\lambda}(\A^{U^+}|\;\ta\,\A^{U+})
\geq \frac{q_v^{|\mb{r}|-m}}{2\ln{q_v}}\frac{q_v^{-j_1|\mb{r}|}}{j_1}\;.
\]
Observe that 
\begin{align*}
  j_1 &= \Big\lceil\; \frac{d-(d-1)\ell}{\min\mb{r}} -
  \log_{q_v}\delta' \;\Big\rceil +1
  = \Big\lceil \;\frac{d-(d-1)\ell}{\min\mb{r}}-\log_{q_v}
  \Big(\;\frac{(\delta_r/2)^{\frac{1}{\min\mb{r}}}}{q_v^2 \;r^{\frac{1}{\max\mb{r}}}}
  q_v^{-j_2}\epsilon^{\frac{m}{d}}\;\Big) \Big\rceil+1\\
  &=\Big\lceil \;\frac{d-(d-1)\ell}{\min\mb{r}} + \Big\lceil\;
  \frac{d-(d-1)\ell}{\min\mb{s}} -\frac{n}{d}\log_{q_v}\epsilon\;
  \Big\rceil+1 -\frac{m}{d}\log_{q_v}\epsilon -\log_{q_v}
  \frac{(\delta_r/2)^{\frac{1}{\min\mb{r}}}}{q_v^2 \;r^{\frac{1}{\max\mb{r}}}}\;
  \Big\rceil+1\\ &\leq (\,d-(d-1)\ell\,)
  \Big(\frac{1}{\min\mb{r}}+\frac{1}{\max\mb{s}}\Big)-\log_{q_v}
  \frac{(\delta_r/2)^{\frac{1}{\min\mb{r}}}}{q_v^2 r^{\frac{1}{\max\mb{r}}}}
  +4-\log_{q_v}\epsilon\;.
\end{align*}
The constants $\eta_A$, $\ell$, $\delta_r$, and $r$ depend only on the
fixed matrix $A\in\cM_{m,n}(K_v)$. Hence there exists a constant
$c(A)>0$ depending only on $d$, $\mb{r}$, $\mb{s}$ and $A$ such that
\[
m-\dimH \Bad_A(\epsilon)\geq c(A)\;\frac{\epsilon^{\,|\mb{r}|}}{\log_{q_v}(1/\epsilon)}\;.
\]
This proves Theorem \ref{theo:introEff1}.
\cqfd

{\small \bibliography{../biblio}

\bigskip
\noindent \begin{tabular}{l} 
\\ 
Department of Mathematical sciences\\ Seoul National University\\
San 56-1, Sillim-Dong, Seoul 151-747, SOUTH KOREA.\\ {\it e-mail:
kimth@snu.ac.kr}
\end{tabular}
\medskip

\noindent \begin{tabular}{l}
Department of Mathematical sciences
and Research Institute of Mathematics\\ Seoul National University\\
San 56-1, Sillim-Dong, Seoul 151-747, SOUTH KOREA.\\ {\it e-mail:
slim@snu.ac.kr}
\end{tabular}
\medskip

\noindent \begin{tabular}{l}
Laboratoire de math\'ematique d'Orsay, UMR 8628 Univ.~Paris-Sud, CNRS\\
Universit\'e Paris-Saclay,
91405 ORSAY Cedex, FRANCE\\
{\it e-mail: frederic.paulin@universite-paris-saclay.fr}
\end{tabular}
}

\end{document}